\documentclass{article}
\usepackage[utf8]{inputenc}
\usepackage{amsmath,amssymb,amsfonts}
\usepackage{graphicx}
\usepackage{lineno,hyperref}
\usepackage{adjustbox}
\usepackage{url}
\usepackage{hhline}
\usepackage{float}
\usepackage{caption}
\usepackage{float}
\usepackage{subfig}
\usepackage{authblk}
\restylefloat{table}
\graphicspath{ {./images/}{./images_grf/} }
\usepackage{placeins}
\usepackage{color}
\usepackage{bm}

\newcommand{\eqn}[1]{
\begin{equation}
#1
\end{equation}
}

\newcommand{\tb}[1]{
{\bf #1}
}

\newcommand{\norm}[1]{
\left\lVert #1 \right\rVert
}

\newcommand{\ip}[1]{
\left\langle #1 \right\rangle
}

\title{Machine-learning custom-made basis functions for partial differential equations}
\author[c]{Brek Meuris}
\author[a]{Saad Qadeer}
\author[a,b]{Panos Stinis}

\affil[a]{Advanced Computing, Mathematics and Data Division, Pacific Northwest National Laboratory, Richland, WA 99354, USA}
\affil[b]{Department of Applied Mathematics, University of Washington, Seattle, WA 98195, USA}
\affil[c]{Department of Mechanical Engineering, University of Washington, Seattle, WA 98195, USA}

\begin{document}
\maketitle
\begin{abstract}
Spectral methods are an important part of scientific computing's arsenal for solving partial differential equations (PDEs). However, their applicability and effectiveness depend crucially on the choice of basis functions used to expand the solution of a PDE. The last decade has seen the emergence of deep learning as a strong contender in providing efficient representations of complex functions. In the current work, we present an approach for combining deep neural networks with spectral methods to solve PDEs. In particular, we use a deep learning technique known as the Deep Operator Network (DeepONet), to identify candidate functions on which to expand the solution of PDEs. We have devised an approach which uses the candidate functions provided by the DeepONet as a starting point to construct a set of functions which have the following properties: i) they constitute a basis, 2)
they are orthonormal, and 3) they are hierarchical i.e., akin to Fourier series or orthogonal polynomials. We have exploited the favorable properties of our custom-made basis functions to both study their approximation capability and use them to expand the solution of linear and nonlinear time-dependent PDEs. 
\end{abstract}

\section{Introduction}\label{sec:intro}

In the last seventy years, scientific computing has made tremendous advancements in developing methods for solving partial differential equations (PDEs) \cite{tadmor2012review,iserles2009first,li2002meshfree}. Spectral methods constitute a significant part of scientific computing's arsenal due to their inherent hierarchical structure, their connection  to approximation theory and their superior performance in various situations \cite{bernardi1997spectral,Boyd_2001,hesthaven2007spectral,canuto2012spectral}. Spectral methods proceed by expanding the solution of a PDE as a linear combination of basis functions, and then estimating the coefficients of the linear combination so that the underlying PDE is satisfied in an appropriate sense. Even though spectral methods can be powerful, their versatility depends strongly on the choice of basis functions, which is far from obvious for many real-world applications.

In the last decade, due to advancements in algorithmic and computational capacity, machine learning, and in particular deep learning, has appeared as a strong contender in providing efficient representations of complex functions \cite{Lecun_2015}. In addition, physics-informed deep learning holds the promise to become a viable approach for the numerical solution of PDEs (see e.g. \cite{karniadakis2021physics,alber2019integrating}). In the current work, we propose a way to combine deep learning and spectral methods in order to solve PDEs. In particular, we put forth the use of deep learning techniques to identify basis functions to expand the solution of a PDE. These basis functions are custom-made i.e., they are constructed specifically for a particular PDE, and they are represented through appropriately defined and trained neural networks. 

Our construction starts with candidate functions that are extracted from a recent deep learning technique for approximating the action of, in general nonlinear, operators, known as the Deep Operator Network (DeepONet) \cite{Lu_2021}. Due to the inherent structure of the DeepONet, those candidate functions are custom-made for a particular PDE (including a class of boundary and initial conditions).  Using those candidate functions as a starting point, we have devised an approach to constructing a set of functions which have the following properties: i) they constitute a basis, 2)
they are orthonormal, and 3) they are hierarchical i.e., akin to Fourier series or orthogonal polynomials. We have exploited the favorable properties of our custom-made basis functions to both study their approximation capability and use them to expand the solution of linear and nonlinear time-dependent PDEs.

The paper is organized as follows. In Section \ref{sec:motive} we offer a motivating example and an overview of our approach. In Section \ref{sec:basis_functions_operator_nets} we present our approach for determining custom-made basis functions using candidate functions provided by the DeepONet. In Section \ref{sec:approx_capability} we examine the approximation capability of the custom-made basis functions. In Section \ref{sec:numerical} we present results using our approach to solve a collection of linear and nonlinear time-dependent PDEs. We conclude in Section \ref{sec:discussion_future_work} with a discussion of our results and suggestions for future work.

\section{Motivation}\label{sec:motive}
The Universal Approximation Theorem (UAT) of Chen and Chen \cite{Chen_1995} guarantees that a two-layer neural network can approximate the action of a continuous nonlinear operator to arbitrary accuracy. In principle, the result suggests an approach for using neural networks to compute solutions of partial differential equations (PDEs) by setting up neural network approximants of operators that map the given data (e.g., initial conditions, boundary data, forcing terms, or diffusivity coefficients) to the solutions. In addition, it posits a linear structure for an approximant, to wit
\eqn{
\mathcal{G}_\text{NN}[g](\tb{y}) = \sum_{k = 1}^w b_k[g]\gamma_k(\tb{y}), \label{UATNN}
}
where $w$ is the network width, $\{b_k\}_{1 \leq k \leq w}$ and $\{\gamma_k\}_{1 \leq k \leq w}$ are branch and trunk nets, $g$ is the problem data, and $\tb{y}$ is an evaluation point. Despite being a powerful theoretical result, this particular approach is infeasible in practice due to the high costs associated with training a fully-connected network of large unknown width.

This methodology was further developed in \cite{Lu_2021}, where the shallow branch and trunk nets were replaced by deeper ones. The resulting architecture, named DeepONet, enables solutions of a large number of PDEs to be computed rapidly and robustly with remarkable accuracy. Strikingly, the technique is agnostic to the nature of the spatial domain and operates at a much lower computational cost than conventional numerical methods. In addition, complementary error analyses \cite{deng2021convergence,lanthaler2021error} provide upper bounds for the approximation error in terms of network size, operator type, and data regularity, while practical performance demonstrates the low generalization and optimization errors associated with this architecture.  

\begin{figure}[ht]
	\centering
	\subfloat[]{{\includegraphics[width=0.45\textwidth]{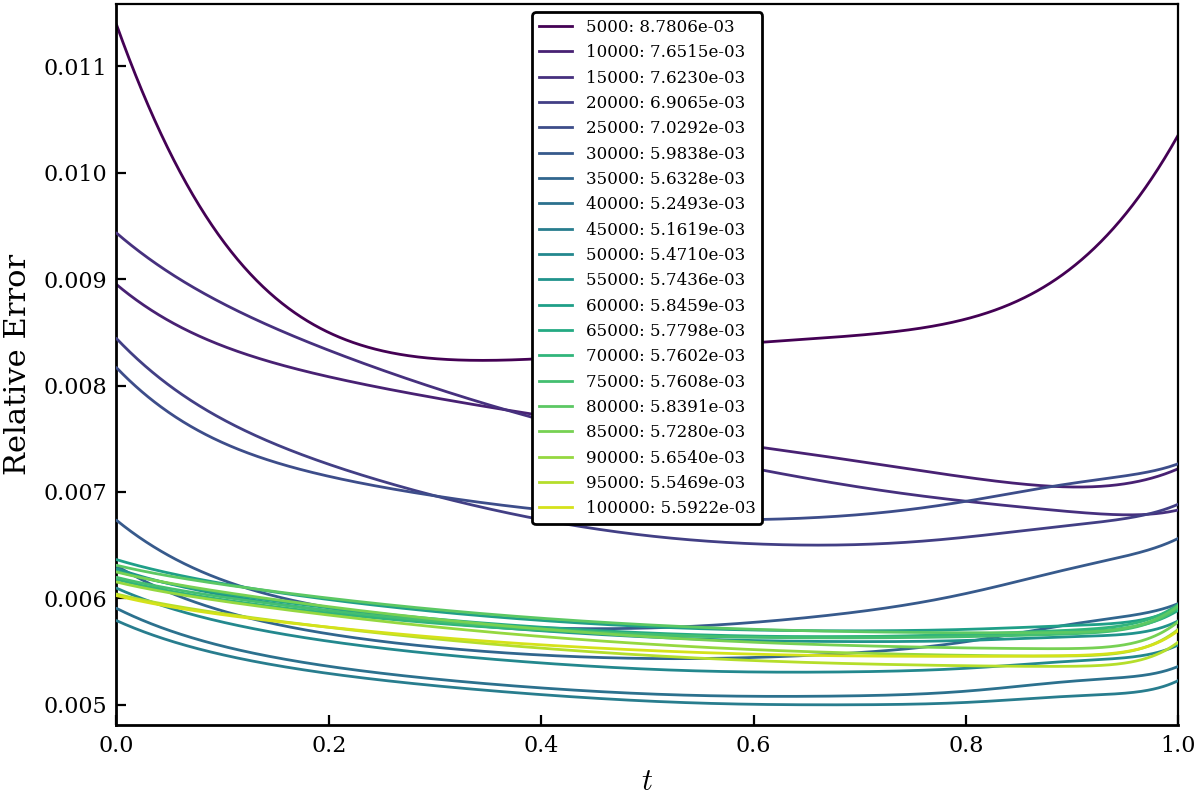} \label{OPNN1}}}
	\qquad
	\subfloat[]{{\includegraphics[width=0.45\textwidth]{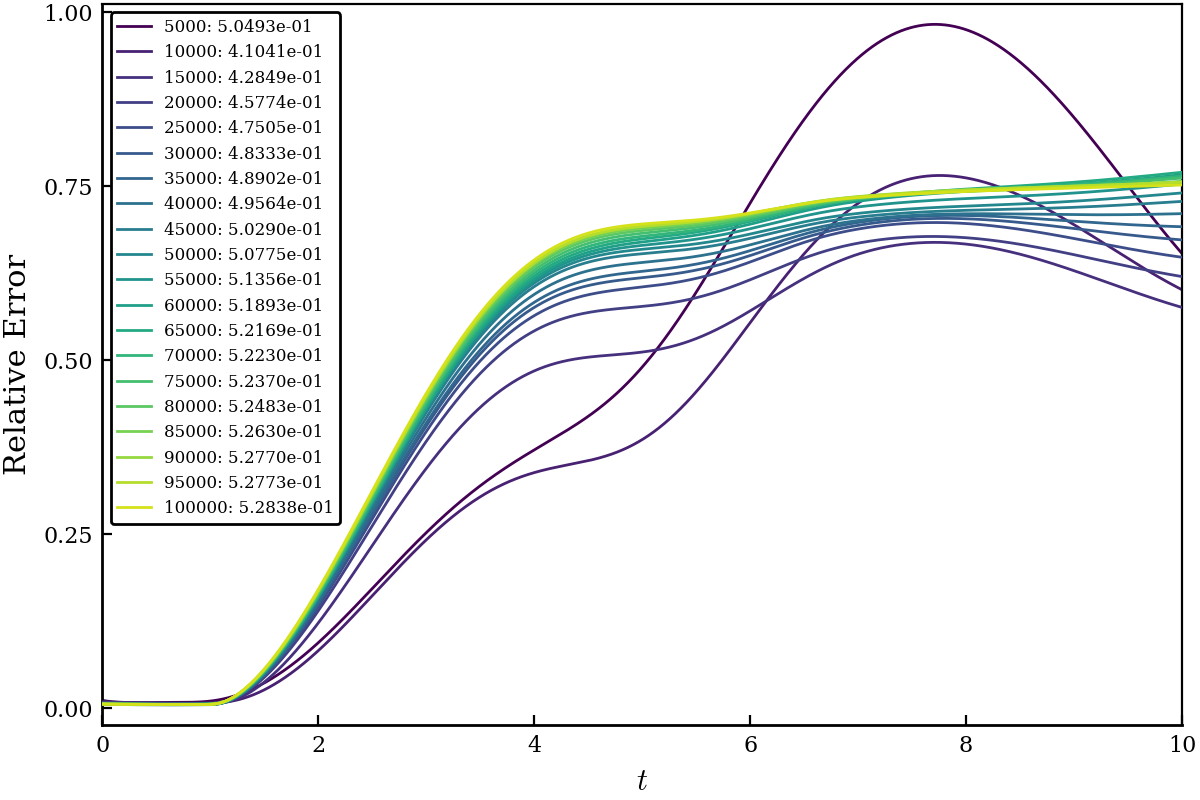}\label{OPNN10} }}
	\caption{Relative errors for the periodic advection problem on $[0,2\pi]$ using a DeepONet trained for $0 \leq t \leq 1$, shown for the initial condition $u_0(x) = \sin^2\left(x/2\right)$, with the number of training epochs going up to $10^5$. The errors are well under 1\% as long as we remain inside in the training interval, as seen in (a). Extrapolating beyond it, however, may lead to large errors, as shown in (b) for $0 \leq t \leq 10$.}
	\label{OPNNErrs}
\end{figure}

Figure \ref{OPNNErrs} shows the results for a DeepONet trained to solve the periodic advection problem $u_t + u_x = 0$ for $x \in [0,2\pi]$, applied to the initial condition $u_0(x) = \sin^2\left(x/2\right)$. The training was performed for $t \in [0,1]$ and the number of epochs increased up to $10^5$. While Figure \ref{OPNN1} shows that the errors are in check for time values in the training domain, the approximate solution quickly loses accuracy outside the training interval, as can be seen in Figure \ref{OPNN10}. This should not be seen as an indictment of the DeepONet approach as it clearly performs satisfactorily on the domain it is designed for. At the same time, it leaves room for the development of tools that can make use of a trained DeepONet to extrapolate solutions outside the training domain. 

In the current work, we present a procedure that harnesses the DeepONet machinery to compute solutions beyond the temporal training interval. Broadly speaking, our approach relies on extracting a hierarchical spatial basis from a trained DeepONet and employing it in a spectral method to solve the PDE of interest. By explicitly making use of the given problem, we expect to overcome the limitations associated with small input-output datasets. At the same time, the use of an operator-regression procedure to obtain the basis functions yields excellent representational capability, particularly on complex spatial domains, with the promise of overcoming the curse of dimensionality. 

Our prescription for obtaining a basis is similar in spirit to Proper Orthogonal Decomposition; however, instead of relying on solution snapshots at different times, our technique makes use of functions whose representation is provided by trained networks. At the same time, the procedure we propose can, in principle, complement any operator regression technique that can furnish high-quality spatial functions, e.g., \cite{li2020fourier,li2020neural}. Our technique of combining machine-learned basis functions for solving variational formulations can also be seen in the context of several important methodologies developed recently \cite{ainsworth2021galerkin,kharazmi2019variational,kharazmi2021hp,khodayi2020varnet}.

\section{Construction of custom-made basis functions using operator neural networks}\label{sec:basis_functions_operator_nets}
In this section, we provide the details of our procedure for obtaining machine-learned basis functions. Let $\mathcal{G}$ be the solution operator for a time-dependent problem on spatial domain $\Omega \subset \mathbb{R}^d$ that maps the initial condition to the solution at later times. A DeepONet $\mathcal{G}_\text{NN}$ trained to approximate $\mathcal{G}$ consists of the branch and trunk nets $\{b_k\}_{1 \leq k \leq w}$ and $\{\gamma_k\}_{1 \leq k \leq w}$ respectively, combined according to
\eqn{
	\mathcal{G}_\text{NN}[u_0](t,{\bf x}) = \sum_{k = 1}^w b_k[u_0]\gamma_k(t,{\bf x}).  \label{DONform}
}

Here, $(t,{\bf x}) \in [0,T] \times \Omega$, the training domain, and $u_0$ is an initial condition. The empirical effectiveness of this architecture and its linear structure suggest that the trunk net functions $\{\gamma_k\}$ possess good spatio-temporal representational capability on the training domain. By evaluating these functions at particular time values, we can therefore obtain a set of spatial functions that can be used to obtain an accurate spectral method. 

However, the trunk net functions may not be well-suited for numerical purposes in their raw form. Most obviously, there is no  particular reason for them to form an orthogonal set with respect to a natural inner product. In addition, they might not be hierarchically ordered, hence closing off an obvious avenue for truncation and reduced order modeling. Finally, they may not necessarily lead to well-conditioned linear systems which would severely limit their effectiveness. As a result, it is essential to develop methods for finding alternate bases for the span of trunk net functions that possess more favourable numerical properties.

Denote by $\ip{\cdot,\cdot}$ the $L^2$ inner product on $\Omega$:
\eqn{
\ip{h_1,h_2} = \int_{\Omega} \overline{h_1(x)} h_2(x) \ dx, \label{IPOmega}
}
and let $\{(x_i,\omega_i)\}_{1 \leq i \leq M}$ be a quadrature rule on $\Omega$ such that
\eqn{
\ip{h_1,h_2} \approx  \ip{\tb{h}_1,\tb{h}_2}_{(\tb{x},W)} := \tb{h}_1^*W \tb{h}_2, \label{IPDisc}
} 
where $\tb{x} = (x_i)_{1 \leq i \leq M}$, $W = \text{diag}(\omega_1,\hdots,\omega_M)$, and $\tb{h}_l = \left(h_l(x_i)\right)_{1 \leq i \leq M}$ for $l = 1,2$. 

We denote the collection of ``frozen-in-time'' trunk net functions by $\{\tau_k\}_{1 \leq k \leq p}$. For example, we can consider $\{\tau_k\}_{1 \leq k \leq p}$ to be the trunk net functions $\{\gamma_k\}_{1 \leq k \leq p}$ evaluated at $t=0.$ In a slight abuse of terminology, we henceforth refer to the $\{\tau_k\}$ as the trunk net functions. We assume that these are linearly independent and that have been normalized (i.e., $\norm{\tau_k} = 1$ for all $k$).  Set $\mathcal{S} = \text{span}\left(\{\tau_k\}_{1 \leq k \leq p}\right)$; recall then that for any $r \leq p$, there exists a (not necessarily unique) $r$-dimensional subspace $\mathcal{S}_r$ of $\mathcal{S}$ such that 
\eqn{
\sum_{k = 1}^p \min_{g_k \in \mathcal{S}_r} \norm{\tau_k - g_k}^2 \leq  \sum_{k = 1}^p \min_{v_k \in \mathcal{V}_r} \norm{\tau_k - v_k}^2,  \label{OptSub}
}
for any $r$-dimensional subspace $\mathcal{V}_r$ of $\mathcal{S}$. The optimal subspace $\mathcal{S}_r$ can be found by assembling the covariance operator
\eqn{
\mathcal{C} = \sum_{k = 1}^p \tau_k \otimes \tau_k = \sum_{k = 1}^p \tau_k\ip{\tau_k,\cdot}, \label{CovOp}
}
and taking the sum of the eigen-spaces associated with the $r$ largest eigenvalues. Note that the eigendecomposition coincides with the singular value decomposition due to the self-adjoint nature of $\mathcal{C}$.

Numerically, the computation of these eigenfunctions can be carried out in two ways. The more obvious approach relies on assembling a discretization of the covariance operator and performing its eigendecomposition. More precisely, defining the $M \times p$ matrix $A$ by $A_{ik} = \tau_k(x_i)$ allows us to build $C_M = AA^*W$, an approximation to $\mathcal{C}$. Its eigendecomposition then yields the eigenfunctions evaluated at the quadrature nodes $\{x_i\}_{1 \leq i \leq M}$. However, this approach is prohibitive in practice as it requires assembling a large square matrix whose dimensions scale with the number of quadrature points $M$, and hence the dimension $d$.

Defining the operators $\mathcal{A}: \mathbb{C}^p \to \mathcal{S}$ and $\mathcal{A}^*: \mathcal{S} \to \mathbb{C}^p$ by
\eqn{
\mathcal{A}{\bf b} = \sum_{k = 1}^p b_k\tau_k, \qquad \left(\mathcal{A}^*f\right)_k = \ip{\tau_k,f}, \ \text{for } 1 \leq k \leq p, \label{Adefn}
}
allows us to write $\mathcal{C} = \mathcal{A}\mathcal{A}^*$. Let 
\eqn{
\mathcal{A} = \sum_{k = 1}^p \sigma_k \phi_k \tb{v}_k^* \label{svdA}
}
be the singular value decomposition (SVD) of $\mathcal{A}$; here, $\sigma_1 \geq \sigma_2 \geq \hdots \geq \sigma_p > 0$ are the singular values, $\{\phi_k\}_{1 \leq k \leq p} \subset \mathcal{S}$ is a collection of functions orthonormal with respect to $\ip{\cdot,\cdot}$, and $\{\tb{v}_k\}_{1 \leq k \leq p} \subset \mathbb{C}^p$ is a set of orthonormal vectors with respect to the Euclidean inner product. This leads to
\eqn{
\mathcal{C} = \sum_{k = 1}^p \sigma_k^2 \left(\phi_k \otimes \phi_k\right), \label{svdC}
}
hence demonstrating that the desired eigenfunctions of $\mathcal{C}$ are simply the $\{\phi_k\}$. 

The alternative approach relies on the observation that knowledge of the singular values $\{\sigma_k\}$ and right singular vectors $\{\tb{v}_k\}$, together with the $\{\tau_k\}$, is sufficient for computing the $\{\phi_k\}$. This can be accomplished by setting $\mathcal{D} = \mathcal{A}^*\mathcal{A}$; this is a $p \times p$ matrix comprising the pairwise inner products $\mathcal{D}_{kl} = \ip{\tau_k,\tau_l}$. Note then that $\mathcal{D} = VS^2V^*$, where $V$ and $S$ are $p \times p$ matrices given by $V_{lk} = (\tb{v}_k)_l$ and $S = \text{diag}(\sigma_1,\hdots,\sigma_p)$. The discretization
\eqn{
\mathcal{D} \approx D_M := A^*WA \label{DAWform}
}
allows the calculation of (approximate) $V$ and $S$ via the SVD and hence the recovery of $\{\phi_k\}_{1 \leq k\leq p}$ by
\eqn{
\phi_k = \sigma_k^{-1}\sum_{l = 1}^p \left(\tb{v}_k\right)_l \tau_l. \label{uRecover}
}

\begin{figure}[ht]
    \centering
    \subfloat[The first few eigenfunctions]{{\includegraphics[width=0.45\textwidth]{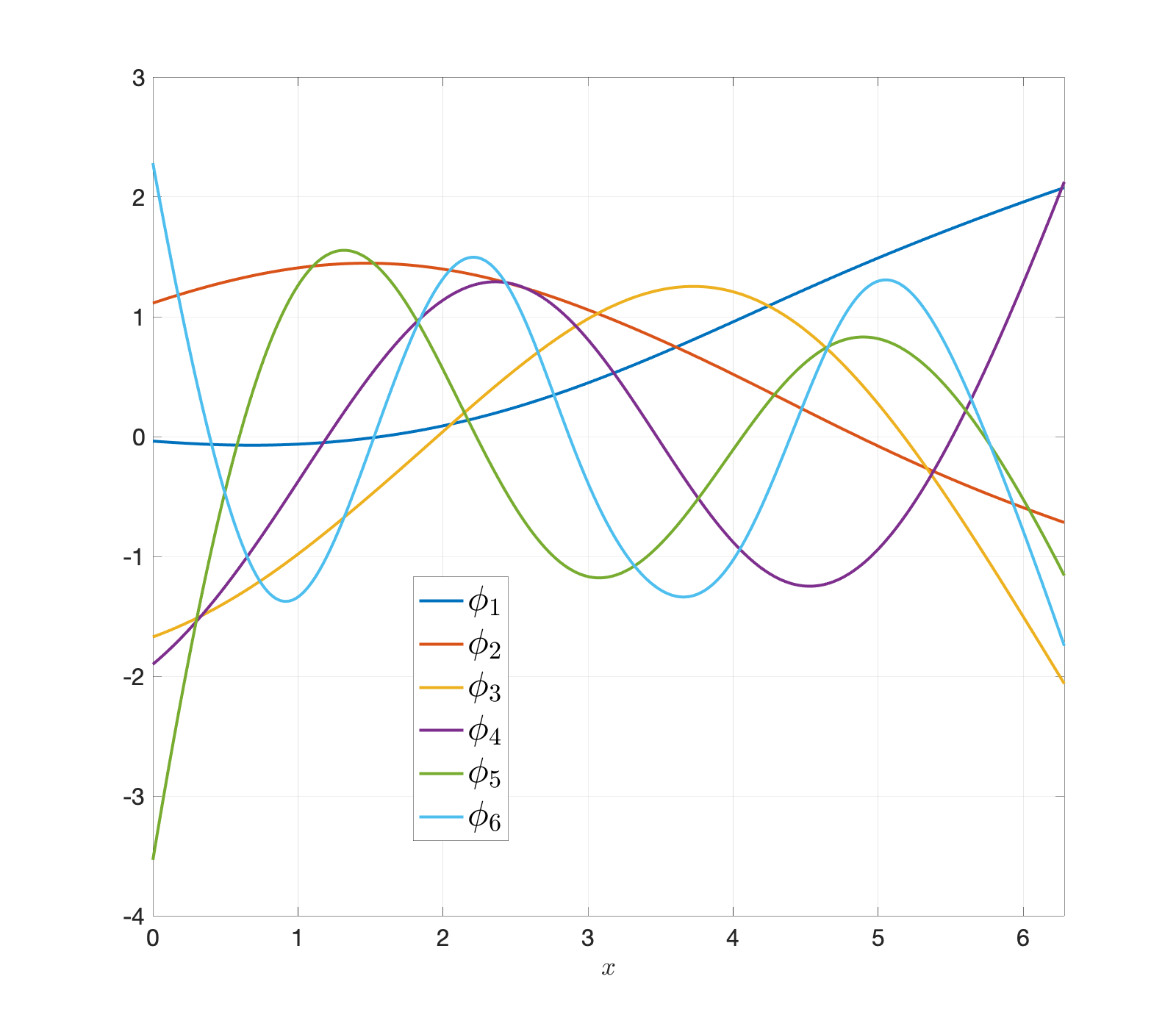} \label{phiKTanh}}}
    \qquad
    \subfloat[The singular values $\{\sigma_k\}_{1 \leq k \leq p}$]{{\includegraphics[width=0.45\textwidth]{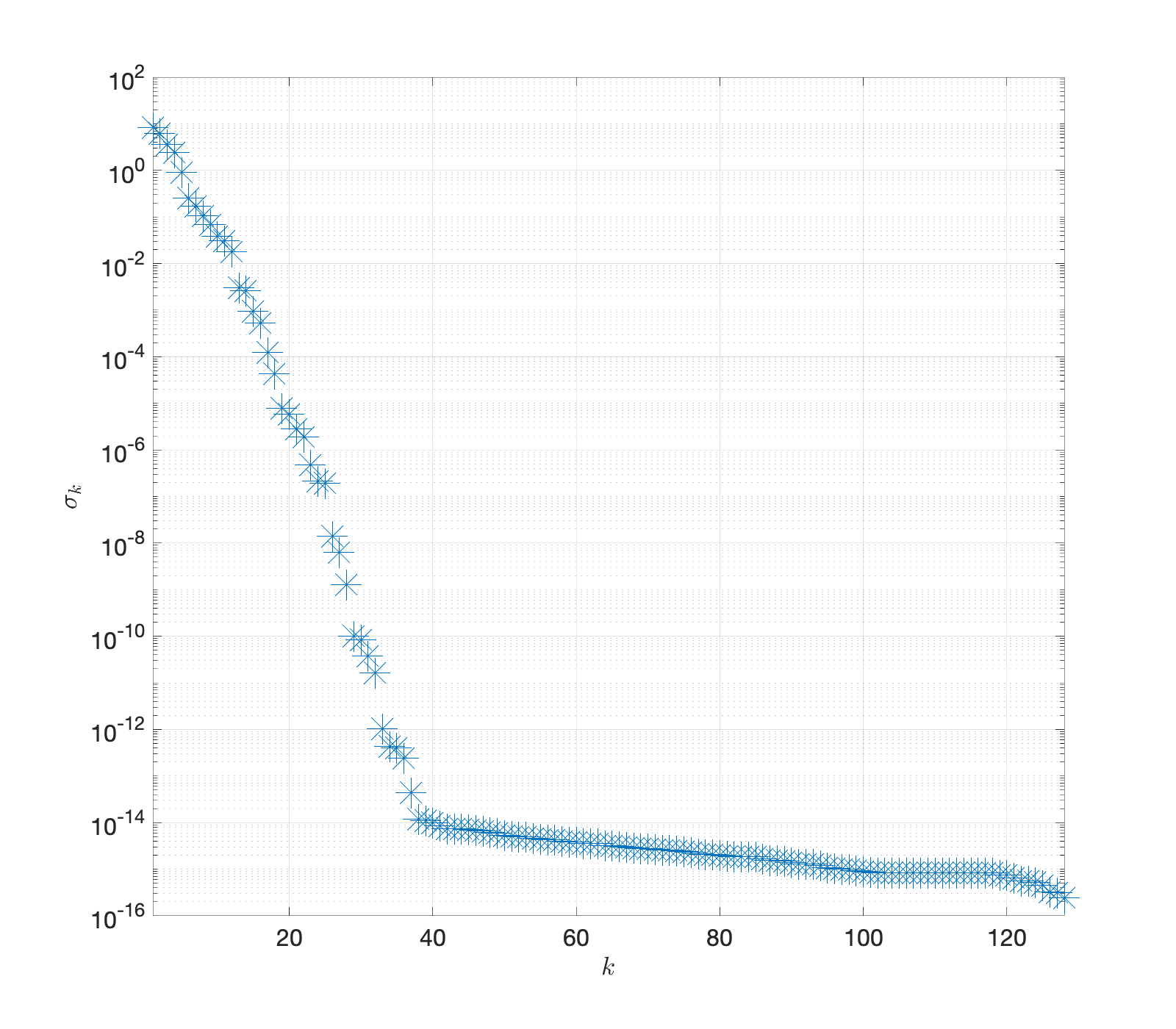}\label{sigKTanh} }}
    \qquad
    \subfloat[Expansion coefficients for $e^{\sin(x)}$]{{\includegraphics[width=0.45\textwidth]{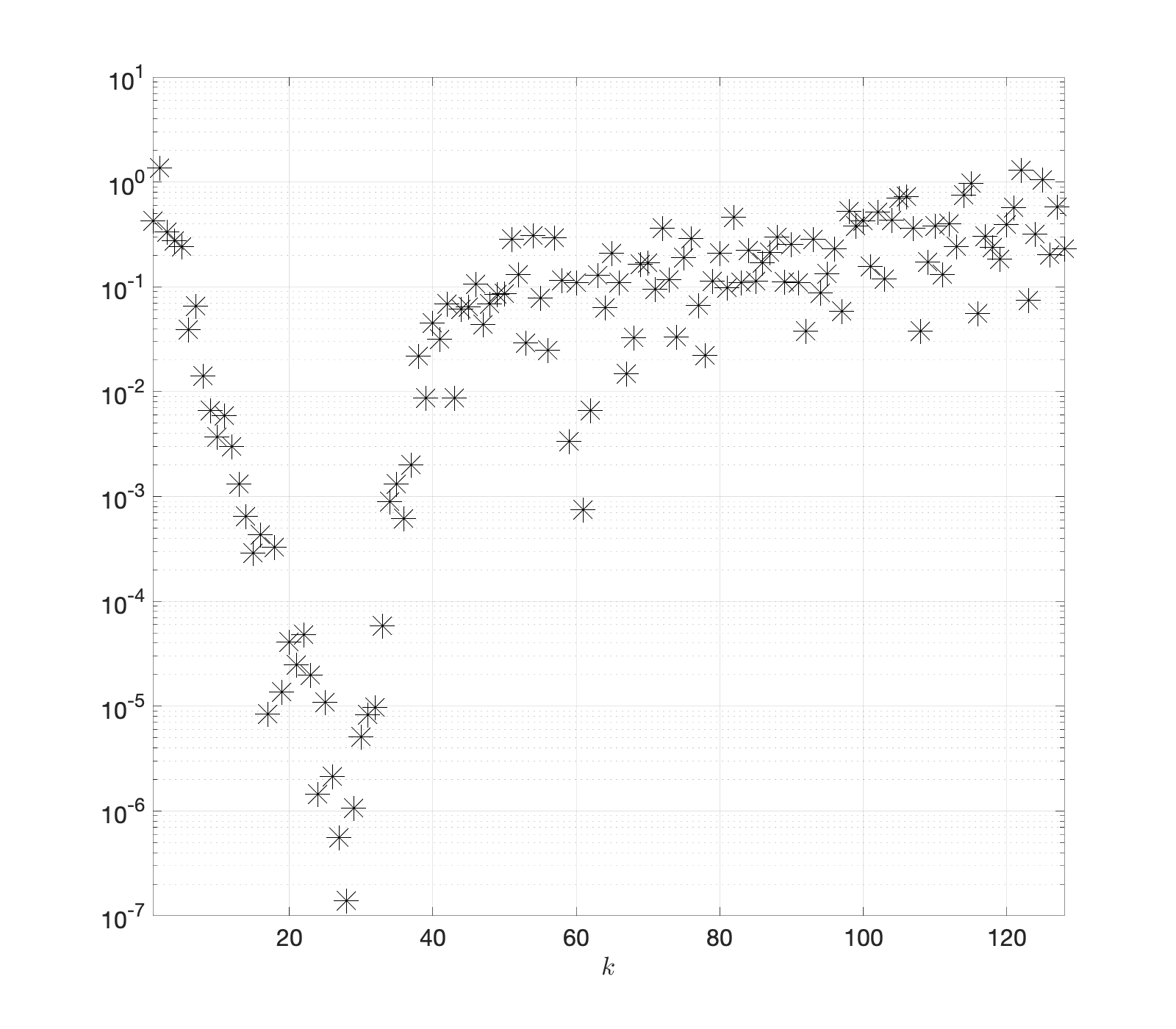} \label{CoeffsTanh}}}
    \qquad
    \caption{Results for the trunk net functions obtained from training the periodic advection problem on $\Omega = [0,2\pi]$ frozen at $t = 0$. We used $M = 2^{10}$ Gauss--Legendre quadrature points and used the ``square root'' formulation with \eqref{uRecover} to compute the eigenfunctions.}
    \label{TanhRes1}
\end{figure}

However, this prescription relies on division by the singular values that may decay rapidly, so the corresponding orthonormal basis calculations can suffer from large errors. In addition, assembling the matrix $\mathcal{D}$ explicitly ends up squaring the singular values, with the result that we can only compute them with accuracy on the order of the square root of machine precision. The latter shortcoming can be overcome by instead making use of the ``square root'' of $\mathcal{D}$, given by $B_M = W^{1/2}A$, and computing its SVD $B_M = QSV^*$.

Figure \ref{TanhRes1} shows the results of applying this procedure to the trunk net functions obtained from training the periodic advection problem on $\Omega = [0,2\pi]$. The trained DeepONet has width $w = 128$ and uses the Tanh activation function; further details about the network training are specified in Table \ref{tbl:network_params} but we draw attention to the choice of the activation function as the results are considerably worse if we employ ReLU instead. Evaluating the temporal trunk net functions at only $t = 0$ yields $p = w = 128$ spatial functions, with $M = 2^{10}$ Gauss--Legendre quadrature nodes used. From the graphs in Figure \ref{phiKTanh}, it appears that successive eigenfunctions undulate with increasing frequency, suggesting that this family possesses a natural frequency-based hierarchy. Also note that, barring the first pair, the zeros of successive eigenfunctions separate each other in the manner of orthogonal polynomials. In addition, the singular values allow us to gauge the contribution of each eigenfunction to $\mathcal{S}$: for singular values below a certain threshold, the corresponding eigenfunctions are essentially noise. As a result, the initial few eigenfunctions (till around $\phi_{36}$ in this example) are, for all practical purposes, sufficiently descriptive. On the flip side, as the ratio of the largest and smallest singular values gives the condition number, we deduce that $B_M$ (and hence $A$) has a large condition number. This hints that any procedure that directly makes use of the $\{\tau_k\}$ (i.e., up to a change of basis) is likely to be ill-conditioned.

One instance of this can be seen in Figure \ref{CoeffsTanh} where the magnitudes of the expansion coefficients $a_k = \ip{\phi_k,f}$ for $f(x) = e^{\sin(x)}$ are shown. The steep decay is interrupted at $k \approx 22$ and the values start increasing. As mentioned earlier, this is a consequence of using \eqref{uRecover}: the round-off error in the calculation of $\phi_k$, and hence $a_k$, is roughly $\epsilon \sigma_k^{-1}$, where $\epsilon$ is the machine precision, and swamps the expansion coefficient values whenever $\sigma_k$ goes below a certain value.

These results call for a technique that avoids making use of \eqref{uRecover} altogether. One way of going about this is to utilize the orthogonal matrix $Q$ obtained from the SVD of $B_M$. Note that the entries of $W^{-1/2}Q$ provide the values of the eigenfunctions at the quadrature points via
\eqn{
\phi_k(x_i) = (W^{-1/2}Q)_{ik} \ \text{ for } 1\leq i \leq M \text{ and } 1 \leq k \leq p. \label{UfromQ}
}

\begin{figure}[ht]
\centering
{\includegraphics[width=0.7\textwidth]{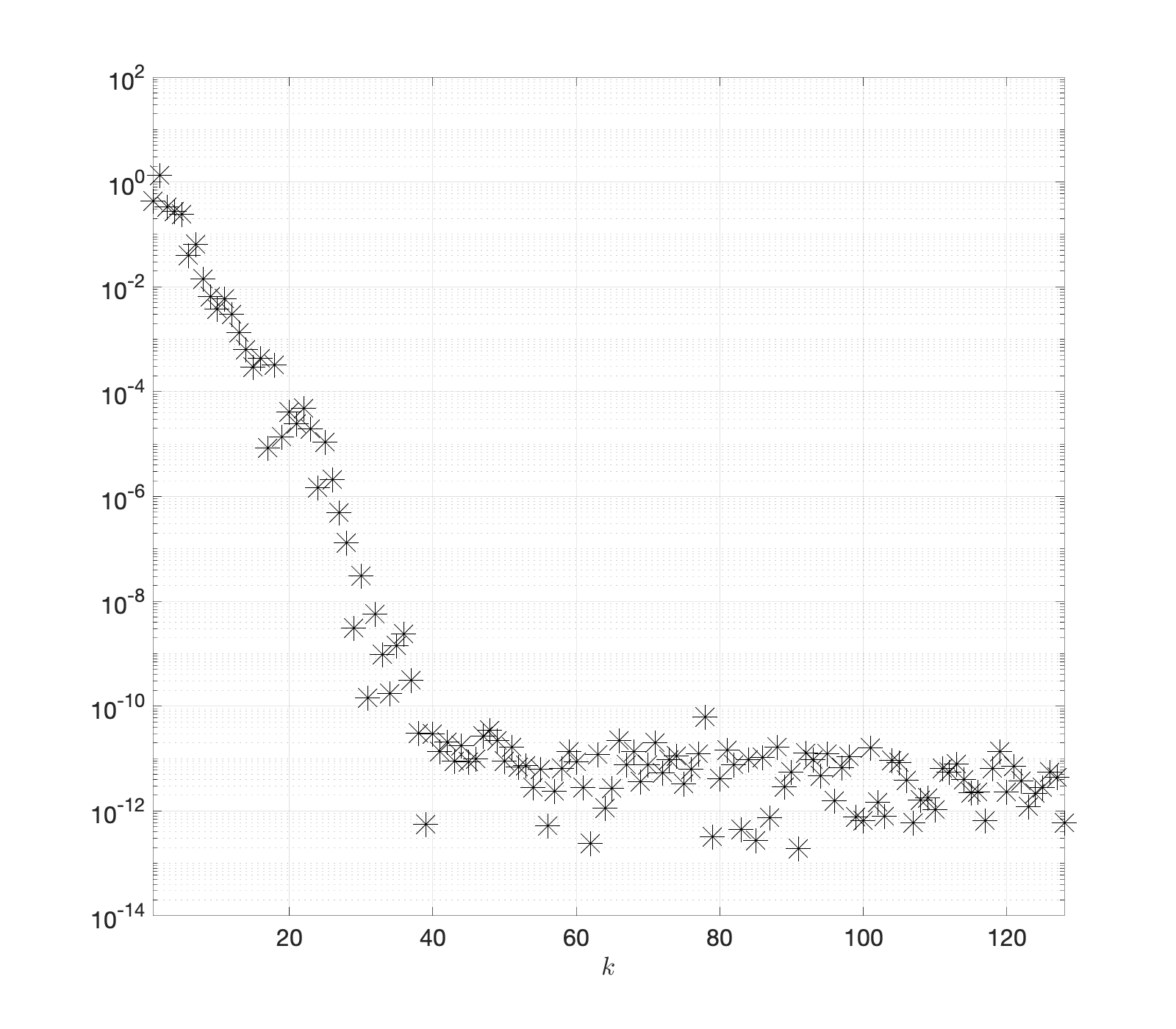}}
\caption{The expansion coefficients for $f(x) = e^{\sin(x)}$ in terms of the $\{\tilde \phi_k\}$ computed using \eqref{UfromQ} and \eqref{LegProj} with $L = 127$.}\label{CoeffTanh2}
\end{figure}

To allow the evaluation of the $\{\phi_k\}$ at points away from the quadrature grid, we can employ either a local interpolation procedure or an orthogonal polynomial expansion. The choice of Gauss--Legendre quadrature points as the $\{x_i\}$ makes the latter option more appealing. For any $L < M$, let $\{q_j\}_{0\leq j \leq L}$ be the orthonormal Legendre polynomials on $[0,2\pi]$ and define the functions $\{\tilde \phi_k\}_{1 \leq k \leq p}$ by
\eqn{
\tilde \phi_k = \sum_{j = 0}^{L} \left(\sum_{i = 1}^M q_j(x_i)\phi_k(x_i) \omega_i\right) q_j. \label{LegProj}
} 

Observe that for $L = M-1$, this is equivalent to a polynomial interpolation through $\{\phi_k(x_i)\}_{1 \leq i \leq M}$. As the $\{\phi_k\}$ are smooth, this polynomial expansion is guaranteed to be highly accurate. By choosing $L$ large enough, we can be fairly confident that the $\{\tilde \phi_k\}$ serve as good approximations to $\{\phi_k\}$. More significantly, the procedure only makes use of operations \eqref{UfromQ} and \eqref{LegProj}, both of which are well-conditioned. This is demonstrated in Figure \ref{CoeffTanh2}, where the magnitudes of the expansion coefficients $\tilde a_k = \ip{\tilde \phi_k,f}$ for $f(x) = e^{\sin(x)}$ can be seen to decay all the way to machine precision (using $L = 127$). This shows that SVD orthonormalization followed by a polynomial expansion forms an accurate and stable method to obtain a hierarchical basis from the trunk net functions.

Most of the ingredients used above generalize to higher dimensions and complex domains easily. However, one obvious pitfall is that higher dimensional analogues of Legendre expansions are not readily available on arbitrary domains. This poses a challenge to our ability to evaluate the orthonormal eigenfunctions away from the quadrature grid points if we forgo the use of \eqref{uRecover}. However, several approaches can be used to bridge this gap, chief among which is a local spline-based interpolation method. In a somewhat related spirit, one can also make use of the recent partition of unity networks proposed in \cite{lee2021partition,trask2021probabilistic} that solve the regression problem on complex geometries by computing a partition of unity coupled with high-order polynomial expansions supported on the different partitions. Yet another approach is to make use of an extension algorithm to extend the eigenfunctions to larger, more regular domains, where Legendre expansions are available. That such extensions exist in principle follows from the Whitney extension theorem \cite{whitney1934analytic,hormander2015analysis}, while several methods have been proposed to carry out this procedure in practice, e.g., \cite{boyd2002comparison,bruno2010high,stein2016immersed,qadeer2021smooth}. These methodologies will be further explored, and the results presented in a future publication.

\section{Assessing the approximation capability of the custom-made basis functions}\label{sec:approx_capability}
The sharp decay in the expansion coefficients in Figure \ref{CoeffTanh2} suggests that the orthonormal basis functions $\{\tilde \phi_k\}_{1 \leq k \leq p}$ are highly adept at approximating smooth functions. In this section, we investigate this capability more fully. From here on, we drop the tildes on the $\{\tilde \phi_k\}$ as, for all intents and purposes, they overlap with the $\{\phi_k\}$. In contrast with the DeepONet analyses \cite{deng2021convergence,lanthaler2021error}, we study the properties of a trained trunk net function space and, in doing so, account for estimation and optimization errors on top of the approximation errors.

For any function $f:[0,2\pi] \to \mathbb{C}$, we can define the orthogonal projection
\eqn{
\mathcal{P}f := \sum_{k = 1}^p \ip{\phi_k,f}\phi_k. \label{Pproj}
}

It follows that
\eqn{
\norm{f - \mathcal{P}f} = \min_{g \in \mathcal{S}} \norm{f - g}. \label{ProjMinimizer}
}

Let $\{q_j\}_{j \geq 0}$ be the orthonormal Legendre polynomials on $[0,2\pi]$ and, for any $r \geq 0$, let 
\eqn{
\mathcal{L}_rf = \sum_{j = 0}^{r} \ip{q_j,f}q_j \label{LegExp}
}
denote the Legendre expansion of $f$ up to $q_{r}$. As $\mathcal{P}\mathcal{L}_rf \in \mathcal{S}$, we have from \eqref{ProjMinimizer}
\eqn{
\norm{f - \mathcal{P}f} \leq \norm{f - \mathcal{P}\mathcal{L}_rf},  \label{AC1}
}
and hence
\eqn{
\norm{f - \mathcal{P}f} \leq \norm{f - \mathcal{L}_rf} + \norm{\mathcal{L}_rf - \mathcal{P}\mathcal{L}_rf}.  \label{AC2}
}

From Parseval's theorem, we have
\eqn{
\norm{f - \mathcal{L}_rf}^2 = \sum_{j \geq r+1} |\ip{q_j,f}|^2, \label{Pars}
}
while \eqref{LegExp} yields
\eqn{
\norm{\mathcal{L}_rf - \mathcal{P}\mathcal{L}_rf} = \norm{\sum_{j = 0}^{r} \ip{q_j,f}\left(q_j - \mathcal{P}q_j\right)}\leq \sum_{j = 0}^{r} |\ip{q_j,f}| \norm{q_j - \mathcal{P}q_j}. \label{LQdiff}
}

Plugging \eqref{Pars} and \eqref{LQdiff} into \eqref{AC2} yields
\eqn{
\norm{f - \mathcal{P}f} \leq \left(\sum_{j \geq r+1} |\ip{q_j,f}|^2\right)^{1/2} + \sum_{j = 0}^{r} |\ip{q_j,f}| \norm{q_j - \mathcal{P}q_j}. \label{AC3}
}

Thus, we can bound the approximation error of an arbitrary function in terms of the Legendre expansion coefficients and the errors made in the approximation of Legendre polynomials using our custom basis. Recall that for smooth $f$, the Legendre expansion coefficients $\ip{q_j,f}$ decay rapidly so, for sufficiently large $r$, we only need to focus on the second term in \eqref{AC3}.

\begin{figure}[htbp]
    \centering
    \subfloat[]{{\includegraphics[width=0.45\textwidth]{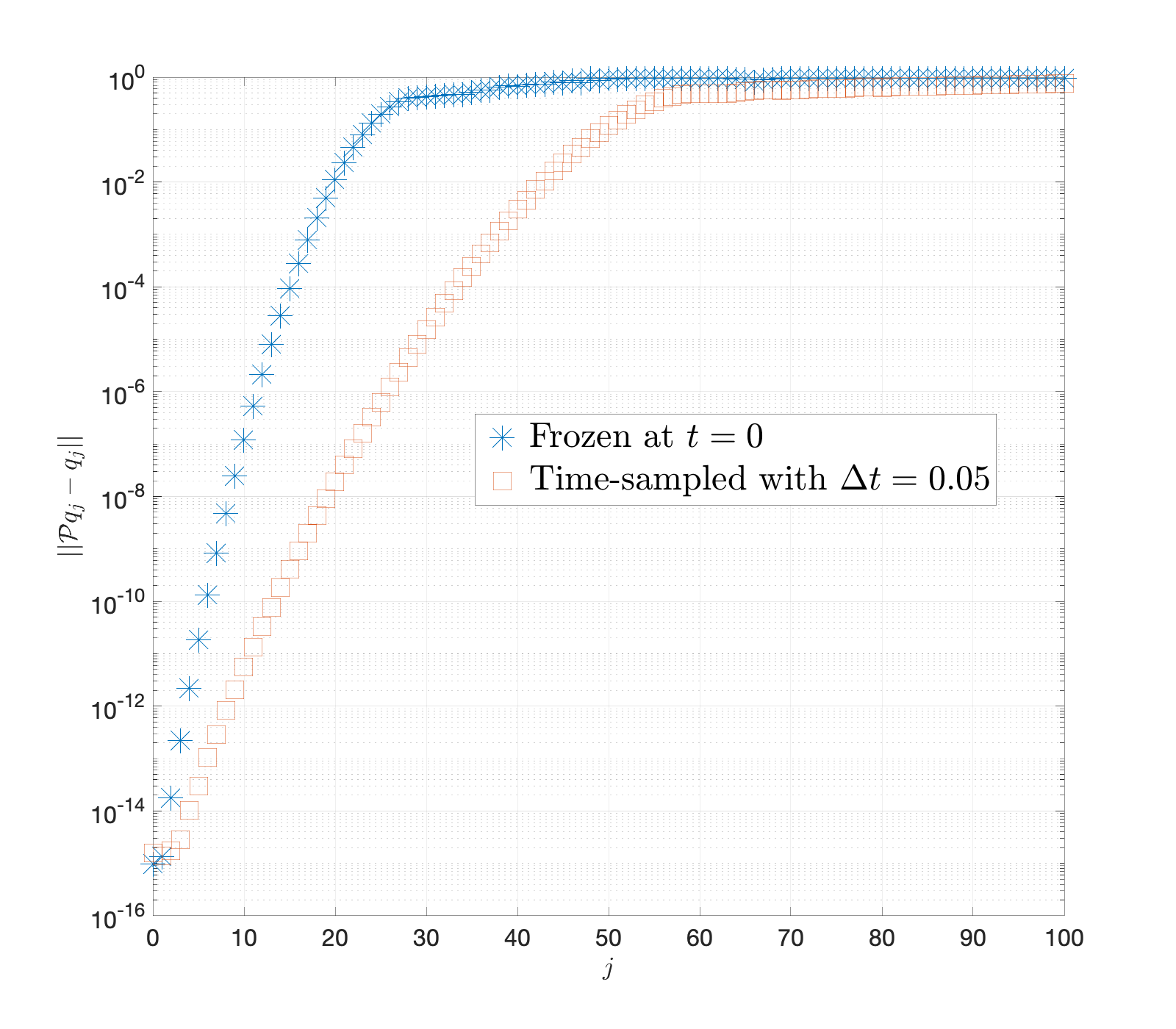} \label{qPsiErrs}}}
    \qquad
    \subfloat[]{{\includegraphics[width=0.45\textwidth]{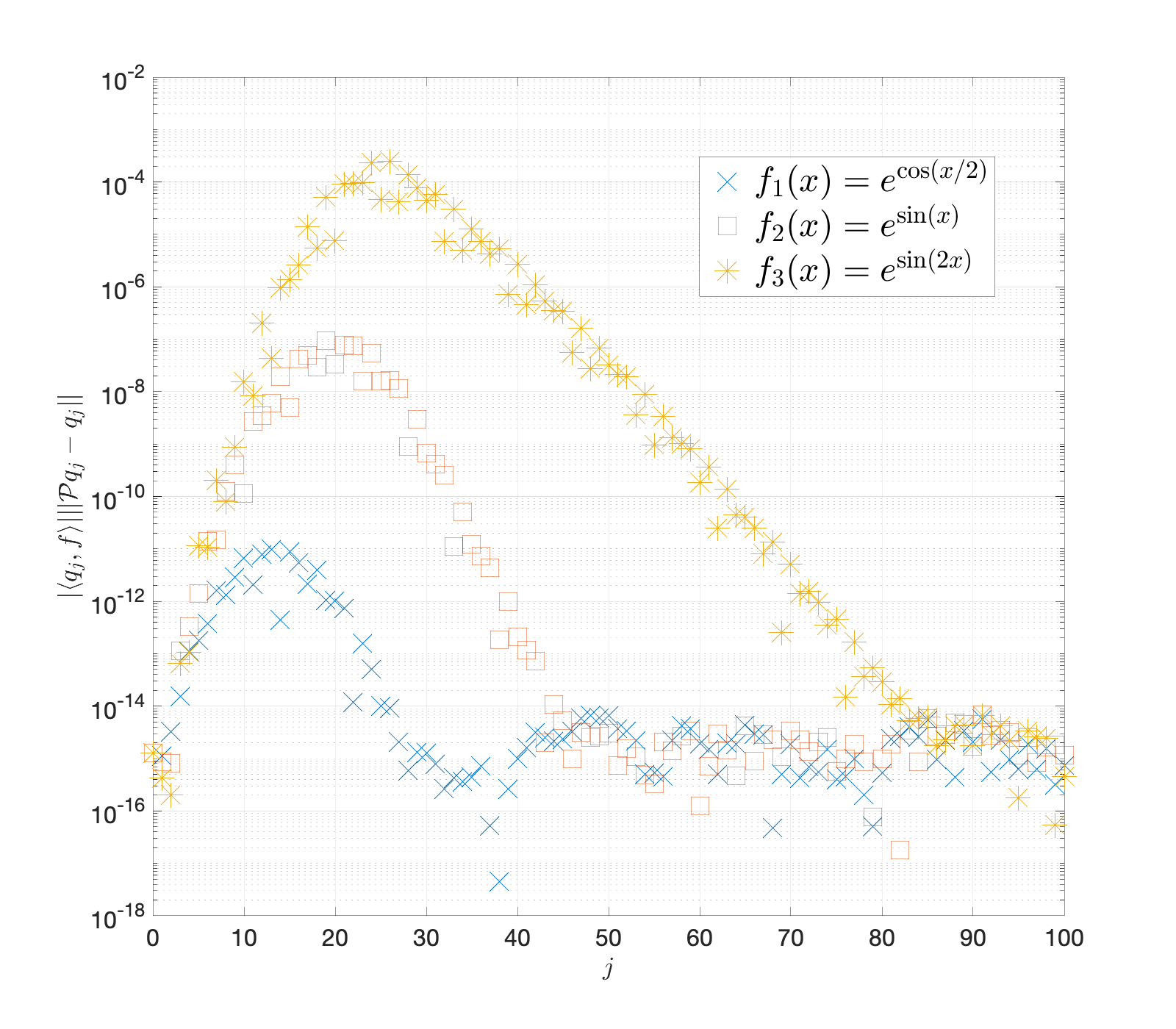}\label{QFPsiErrs} }}
    \qquad
    \subfloat[]{{\includegraphics[width=0.45\textwidth]{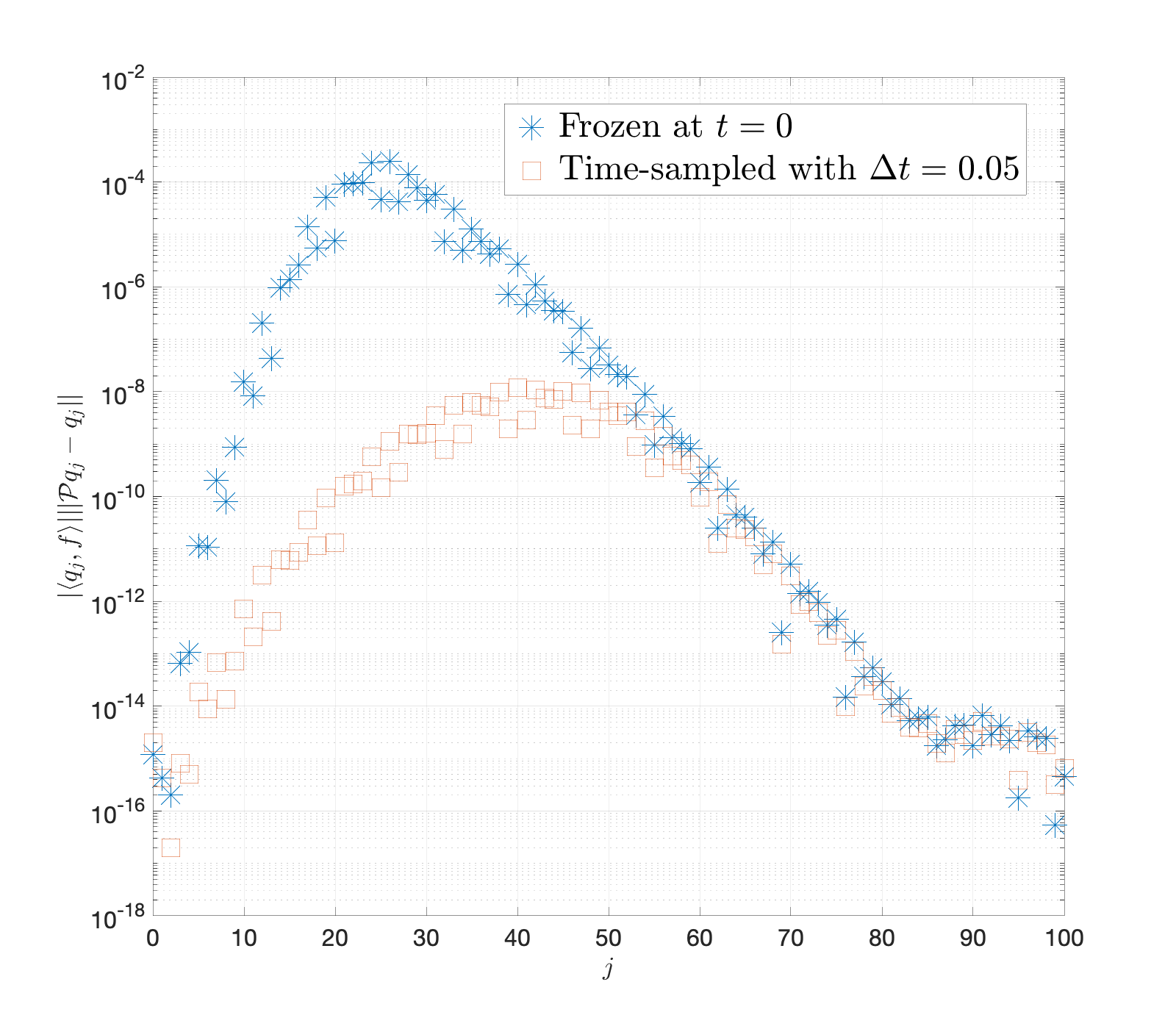}\label{QFPsiErrs3} }}
    \caption{(a) Errors in the approximation of the $j$th Legendre polynomial using the custom basis functions. The errors appear to increase exponentially, hinting at the difficulty in approximating high-frequency functions. The growth is markedly slower when the basis functions are drawn from the larger collection of time-sampled trunk net functions. (b) Terms in the upper bound \eqref{AC3} computed for functions with varying frequencies, using the $t = 0$ trunk net functions. Low frequency functions possess rapidly decaying Legendre coefficients and hence are able to damp the growth seen in (a); this becomes harder with increasing frequency. (c) The upper bound terms in \eqref{AC3} for $f(x) = e^{\sin(2 x)}$, shown for both cases considered in (a). The slower growth associated with the time-sampled trunk functions leads to much smaller error bounds, helping explain the lower errors in Table \ref{tab:cbf_errs}.}
    \label{LegErrs}
\end{figure} 

Figure \ref{qPsiErrs} shows that the errors $\norm{q_j - \mathcal{P}q_j}$ in approximating Legendre polynomials using the basis functions derived from the advection DeepONet increase exponentially with $j$ before levelling off. Along with the basis functions derived from the trunk net frozen at $t = 0$, we present the results from sampling the trunk net at equispaced values in the temporal domain (we call those time-sampled functions). Using $w = 128$ as before and choosing the time-step size $\Delta t = 0.05$, we obtain $p = (1 + 1/\Delta t)w = 2688$ trunk net functions from which the basis functions are drawn. For this case, we use $M = 2^{12}$ Gauss--Legendre quadrature nodes; to allow for a fair comparison with the $t = 0$ case, we employ 128 basis functions for both tests. The growth is markedly slower in the time-sampled case, suggesting that it is better equipped to approximate high frequency functions. The plateauing of the errors in both cases is an expected outcome of using an architecture based on the UAT for operators as it is only valid for compact subsets of $C(\Omega)$ and hence cannot account for arbitrarily large frequencies. Nevertheless, this diagram goes some way towards assessing the completeness properties of the custom basis functions by showing the errors in approximating successive elements of a complete polynomial basis.

As these errors appear in tandem with the decaying Legendre expansion coefficients in \eqref{AC3}, we can expect this steep growth to be damped. This is illustrated in Figure \ref{QFPsiErrs}, where the $|\ip{q_j,f}| \norm{q_j - \mathcal{P}q_j}$ have been computed for $f_1(x) = e^{\cos(x/2)}$, $f_2(x) = e^{\sin(x)}$, and $f_3(x) = e^{\sin(2x)}$. It can be seen that the low frequency $f_1$ possess rapidly decaying Legendre coefficients and hence is more capable at damping the growth seen in Figure \ref{qPsiErrs}; for increasing frequencies, this becomes a much more difficult task. Finally, in Figure \ref{QFPsiErrs3}, we compare the upper bound terms from the $t = 0$ trunk net functions and the time-sampled ones for $f_3$. Due to the superior capability of the latter at approximating higher order Legendre polynomials, the upper bound terms are much smaller, hinting that the actual errors $\norm{f - \mathcal{P}f}$ would be lower as well. That this is indeed the case is confirmed in Table \ref{tab:cbf_errs}.

\begin{table}[ht]
    \centering
    \begin{tabular}{|c||c|c|}\hline
       $f(x)$ & Frozen at $t = 0$ & Time-sampled  \\ \hline\hline
        $e^{\cos(x/2)}$  &  $5.69 \times 10^{-11}$ & $4.79 \times 10^{-15}$ \\ \hline
        $e^{\sin(x)}$  & $5.59 \times 10^{-7}$ & $4.83 \times 10^{-15}$ \\ \hline
        $e^{\sin(2x)}$  & $1.33 \times 10^{-3}$ & $6.99 \times 10^{-15}$ \\ \hline
    \end{tabular}
    \caption{The errors $\norm{f - \mathcal{P}f}$ computed for functions with increasing frequencies on $[0,2\pi]$, using the custom basis functions derived from trunk net functions frozen at $t = 0$ as well as from those sampled at time values spaced apart at $\Delta t = 0.05$. To account for the much larger dimension of the space in the latter case, we only make use of the first 128 basis functions. The time-sampled approach yields much smaller errors, indicating the scope of this approach.}
    \label{tab:cbf_errs}
\end{table}

\section{Numerical results}\label{sec:numerical}

Once we have obtained the custom-made basis functions, we use a standard Galerkin method to obtain equations for the evolution of the coefficients of the expansion (please see Appendix \ref{app:custom_basis_evolution} for details). We use four one-dimensional PDEs on periodic domains to validate our approach. We consider two linear PDEs, the advection and the advection-diffusion problems, and two nonlinear ones, the viscous and inviscid Burgers equations. To demonstrate the temporal interpolation and extrapolation capabilities of the presented procedure for the advection, advection-diffusion, and viscous Burgers equations, results are presented for times beyond the operator neural network training interval. For the inviscid Burgers equation the temporal interpolation capabilities are demonstrated, with commentary on extrapolation provided in Section \ref{sec:discussion_future_work}. Additionally, to also demonstrate the ability to extrapolate outside the training input function space, results for a within training distribution and out of training distribution initial condition are shown for all four PDEs.

\begin{figure}[h!]
        \centering
        {{\includegraphics[width=0.45\textwidth]{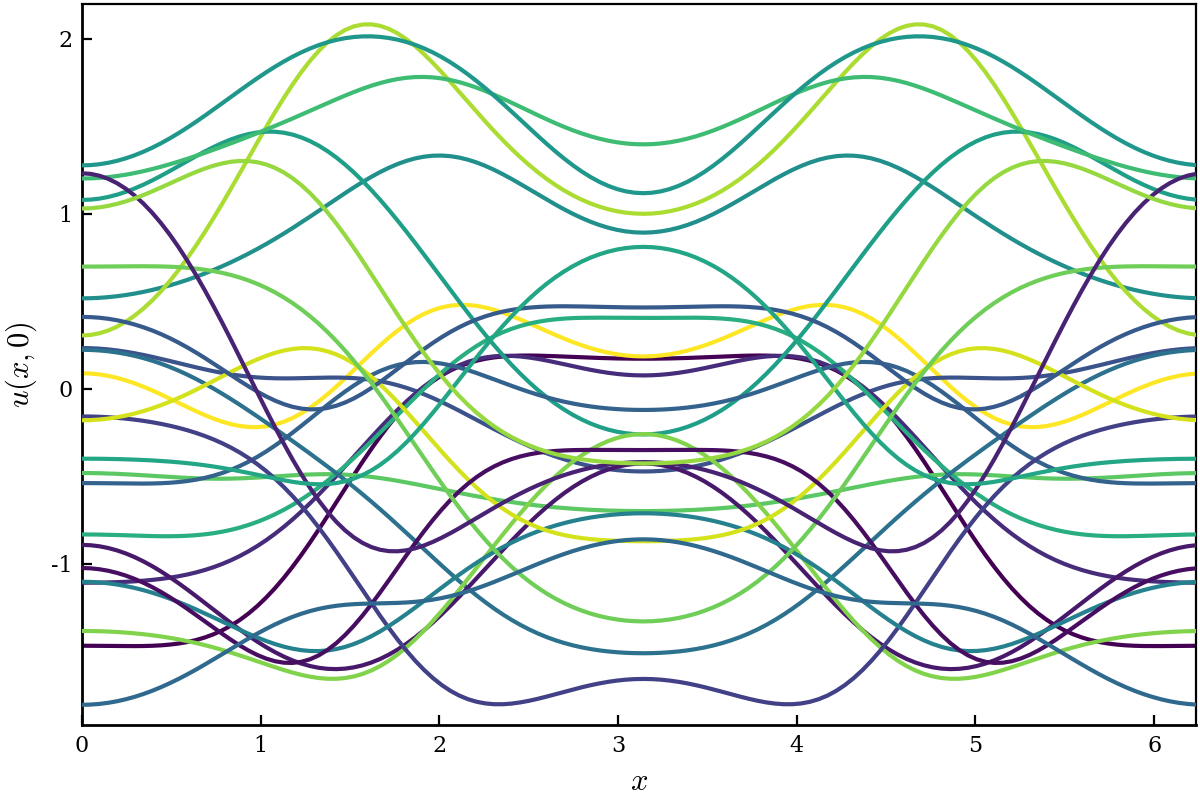} }}
        \caption{Twenty-five example initial conditions sampled from the Gaussian random field.}
        \label{fig:init_conds}
    \end{figure}

\begin{figure}[h!]
        \centering
        {{\includegraphics[width=0.45\textwidth]{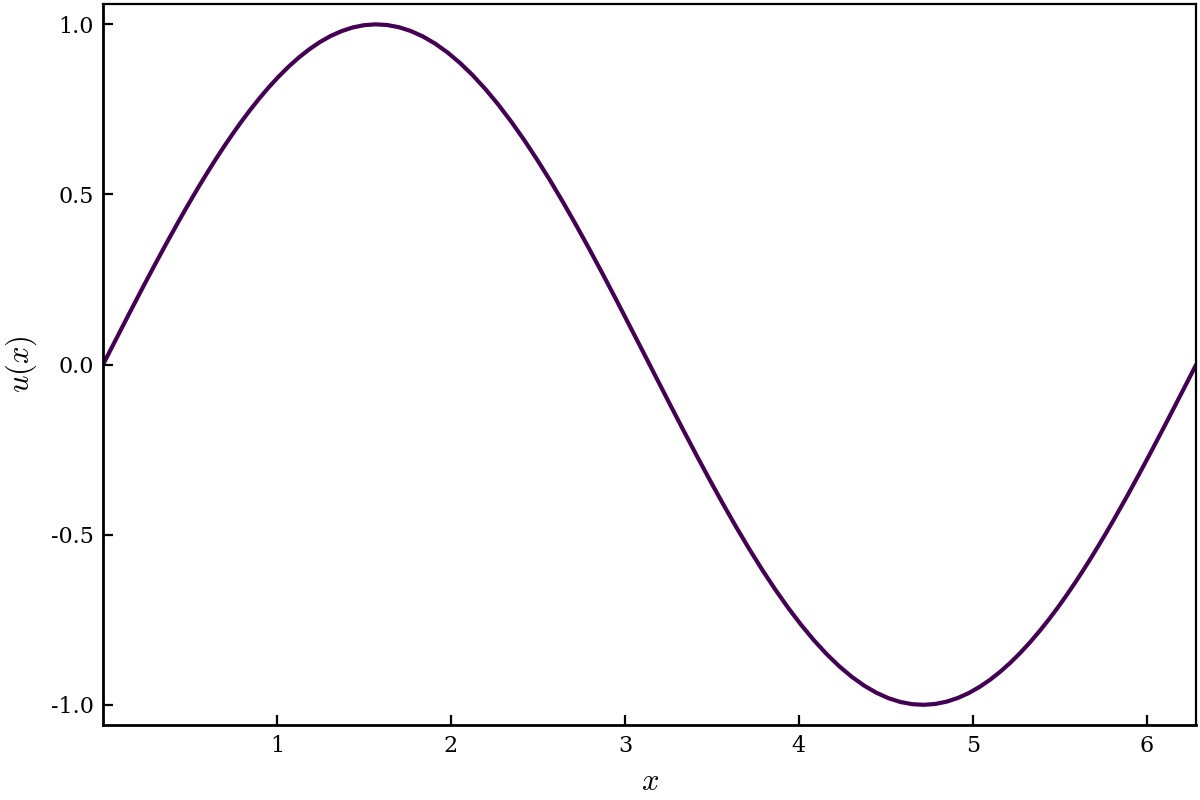} }}
        \caption{Out of distribtuion test initial condition, $u(x,0) =\sin(x)$, for each PDE.}
        \label{fig:init_out_dist}
    \end{figure}

The DeepONets that underlie our construction were trained using the network parameters presented in Table \ref{tbl:network_params}. All branch layers include bias, and the output layer of the branch network does not apply an activation function. Additionally, the network width was constant for all layers. We sampled $f(\sin^2(x/2))$ from a mean zero Gaussian random field with covariance kernel,
\begin{equation}\label{eqn:grf_dist}
    \kappa_l(x_1,x_2)=e^{\frac{-\Vert x_1 - x_2\Vert ^2}{2l^2}}, \;\; l = 0.5,
\end{equation}
to generate the random training and testing intial conditions. Twenty-five example initial conditions are shown in Figure \ref{fig:init_conds}. The in-distribution initial conditions for each PDE are presented in their respective sections below (see Figures \ref{fig:init_in_dist_advection}, \ref{fig:init_in_dist_advection_diffusion},  \ref{fig:init_in_dist_viscous_burgers}, and  \ref{fig:init_in_dist_inviscid_burgers}). The out of distribution initial condition tested for each of the example PDEs is shown in Figure \ref{fig:init_out_dist}. To indicate the effects of random initialization of the neural networks, the mean testing error and corresponding standard deviation based on three training runs are presented in Table \ref{tbl:network_errors} for each example PDE. 

The ground truth data for the advection, advection-diffusion, and viscous Burgers equations was generated by writing the solutions in terms of $M = 128$ Fourier modes,
\begin{equation}\label{eqn:fourier_expansion}
        u_\text{G}^M(t,x) = \sum_{k = -M/2}^{M/2-1} \hat{u}_k(t)e^{ikx},
\end{equation}
with the most negative mode set to zero to avoid asymmetry between the positive and negative modes \cite{canuto2012spectral}. The resulting systems of differential equations were solved for $t \in [0,1]$ using a Runge-Kutta-Dormand-Prince integrator with adaptive step size, relative error tolerance $10^{-10}$, and absolute tolerance $10^{-14}$ \cite{differential_equations_2017}. The solution was saved at times values $10^{-3}$ apart for the linear and $10^{-4}$ for the nonlinear PDE examples. The convolution sum that results from using \eqref{eqn:fourier_expansion} for solving the viscous Burgers equation was evaluated by padding the Fourier solution using the $3/2$-rule for de-aliasing \cite{canuto2012spectral}, transforming the solution to real space, and then computing the fast Fourier transform (FFT) of the product of the real space solution with itself. 

    \begin{table}[h!]
        \begin{center}
            \begin{tabular}{ | l || c |}
                \hline
                Parameter & Setting \\
                \hhline{|=#=|}
                Activation functions & Tanh \\
                \hline
                Optimizer & Adam \\
                \hline
                Error & Mean squared error \\
                \hline
                Learning rate & 0.00001 \\
                \hline
                Number of training epochs & 50000 \\
                \hline
                Number of sensors & 128 \\
                \hline
                Number of solution locations & 100 \\
                \hline
                Number of training initial conditions & 500 \\
                \hline
                Number of testing initial conditions & 1000 \\
                \hline
                Branch net depth & 2 \\
                \hline
                Branch net width & 128 \\
                \hline
                Trunk net depth & 3 \\
                \hline
                Trunk net width & 128 \\
                \hline
                Weight initialization & Glorot uniform \\
                \hline
                Bias initialization & Zero \\
                \hline
                Mini-batch size & 100 \\
                \hline
            \end{tabular}
            \caption{Parameter settings for training the DeepONets.}
            \label{tbl:network_params}
        \end{center}
    \end{table}

The ground truth data for the inviscid Burgers equation was generated by utilizing a MUSCL (monotonic upwind scheme for conservation laws) scheme with a second-order Roe scheme for the flux and a minmod slope limiter \cite{Pletcher_2012,Yang_1992}. The spatial domain was discretized using 4096 points and the resulting equations were solved for $t \in [0,3.5]$ using a Bogacki--Shampine $3/2$ integrator with adaptive step size, relative error tolerance $10^{-6}$, and absolute tolerance $10^{-8}$ \cite{differential_equations_2017}. To train the DeepONet, the solution was saved $10^{-4}$ time units apart and down sampled from the 4096 solution locations to the 128 uniformly spaced sensor locations.

    \begin{table}[h!]
        \begin{center}
            \begin{tabular}{ | l || c | c |}
                \hline
                PDE & Mean & Standard deviation \\
                \hhline{|=#=|=|}
                Advection & $3.37\times 10^{-5}$ & $3.59\times 10^{-6}$ \\
                \hline
                Advection-diffusion & $6.90\times 10^{-5}$ & $2.69\times 10^{-5}$ \\
                \hline
                Viscous Burgers & $2.36\times 10^{-3}$ & $1.79\times 10^{-4}$ \\
                \hline
                Inviscid Burgers & $4.08\times 10^{-2}$ & $6.47\times 10^{-4}$ \\
                \hline
            \end{tabular}
            \caption{DeepONet mean testing errors and standard deviation for all the example PDEs based on three training runs each.}
            \label{tbl:network_errors}
        \end{center}
    \end{table}

For each PDE, the custom bases were constructed as in Section \ref{sec:basis_functions_operator_nets} with $L=127$ (see \eqref{LegProj}) and with the trunk net functions sampled at $t = 0$. Recall that if the singular value corresponding to a basis function is too small, the function itself is predominantly noise. Since we utilize $b$ expansion coefficient(s) to enforce the boundary condition(s), the basis functions that appear later in the hierarchy need to be well behaved. As a result, the number $r$ of basis functions used to solve the PDE is determined empirically and specified to be larger than $10^{-12}$, but as near this lower bound as possible while still allowing for stable evolution in time. This can be seen as a form of model reduction without the memory term. The treatment of the boundary conditions when using the tau-method leads to the sensitivity associated with specifying the number of basis functions (see Appendix \ref{app:custom_basis_evolution}). We are also pursuing an alternative approach for the boundary conditions in the spirit of Discontinuous Galerkin methods \cite{lesaint1974finite,cockburn1998local,shu2001different} which will appear in a future publication.

All the examples presented are solved on a 128 node Gauss--Legendre quadrature grid to facilitate the usage of this highly accurate scheme for the calculation of inner products. The differentiation of basis functions was performed using automatic differentiation. The system of ODEs for the expansion coefficients was integrated in time using the classical fourth order Runge-Kutta scheme. Fixed time stepsizes of $10^{-3}$ and $10^{-4}$ were used for the linear and nonlinear PDEs respectively.

The error presented for each PDE is a relative Euclidean two-norm error defined by
\begin{equation}
    E_2(t) = \frac{\norm{u^r(t,\cdot) - u^M_\text{G}(t,\cdot)}_2} {\norm{u^M_\text{G}(t,\cdot)}_2}, \label{ErrorDef}
\end{equation}
where $u^M_\text{G}$ is the ground truth solution computed using a Fourier expansion or MUSCL scheme, and $u^r$ is the custom basis function solution using $r$ basis functions. When comparing against a Fourier solution, the $M$ = 128 Fourier expansion coefficients and \eqref{eqn:fourier_expansion} are used to approximate the solution at the non-uniform quadrature nodes. In the case of MUSCL, the solution computed on the spatial grid consisting of 4096 discretization points is interpolated using piecewise linear interpolation to approximate the solution at each of the non-uniform quadrature nodes. The average error over $[0,T]$ is defined as
\begin{equation}
    \overline{E}_2 = \frac{1}{T}\int_0^T E_2(t)dt
\end{equation}

All work was performed using the Julia language \cite{julia_2017} with the DeepONets implemented using the Flux machine learning library \cite{flux_2018}. The library of code used to generate the presented results is available upon request. 


\subsection{Advection equation} \label{sec:advection_results}

\begin{figure}[h!]
        \centering
        \subfloat[The singular values]{{\includegraphics[width=0.45\textwidth]{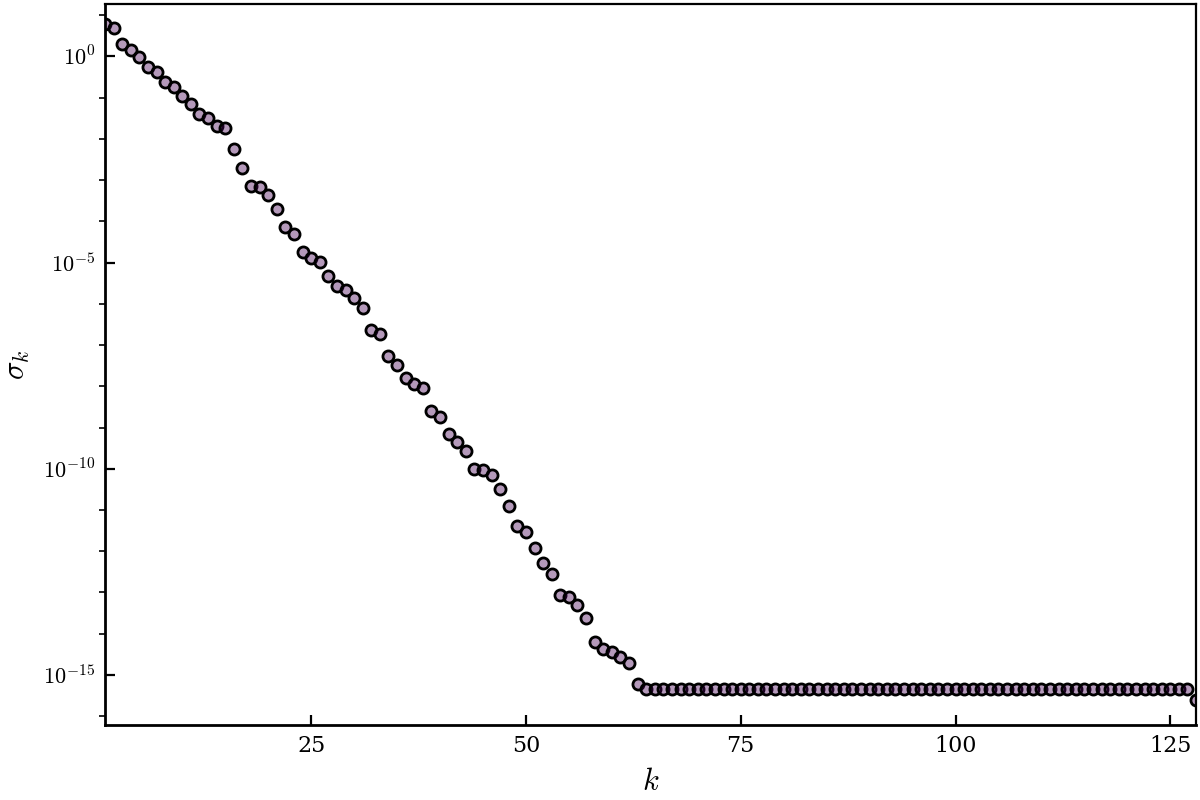} }}
        \qquad
        \subfloat[The first 40 custom basis functions]{{\includegraphics[width=0.45\textwidth]{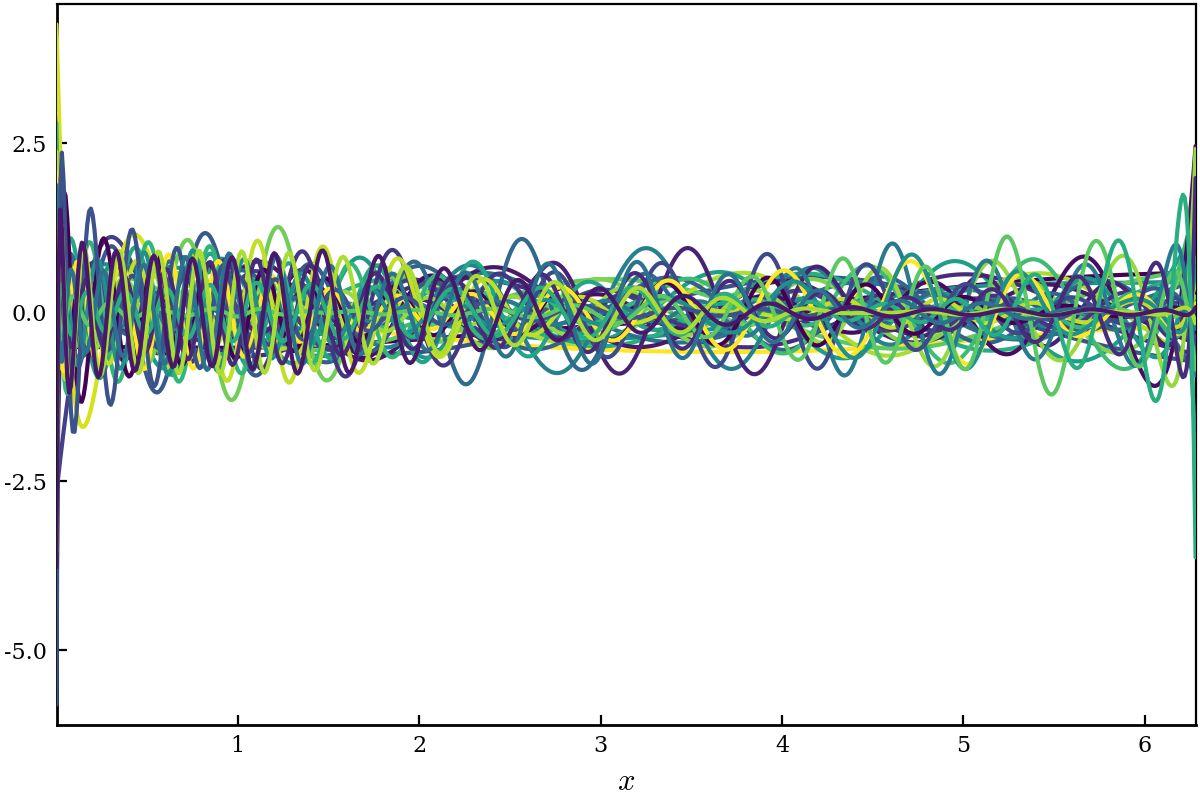} }}
        \qquad
        \caption{Singular values and custom basis functions for the periodic advection equation.}
        \label{fig:sv_basis_advection}
    \end{figure}

The one-dimensional advection equation with periodic boundary conditions is given by
\begin{equation}\label{eqn:advection}
    \frac{\partial u}{\partial t} + \frac{\partial u}{\partial x} = 0, \;\; x \in [0,2\pi].
\end{equation}
The singular value spectrum and the basis functions satisfying the singular value threshold, $\sigma_k > 10^{-9}$, that results from the construction presented in Section \ref{sec:basis_functions_operator_nets} are shown in Figure \ref{fig:sv_basis_advection}.

Using the procedure outlined in Appendix \ref{app:custom_basis_evolution}, we obtain
\begin{equation} 
        \frac{da_m(t)}{dt} = - \sum_{k=1}^r a_k(t) \langle \phi_{m},\phi'_{k} \rangle, \quad m=1,2,\dots,r-b,
\end{equation}
where $r = 40$ and the primes denote differentiation in space.

\begin{figure}[h!]
        \centering
        {{\includegraphics[width=0.45\textwidth]{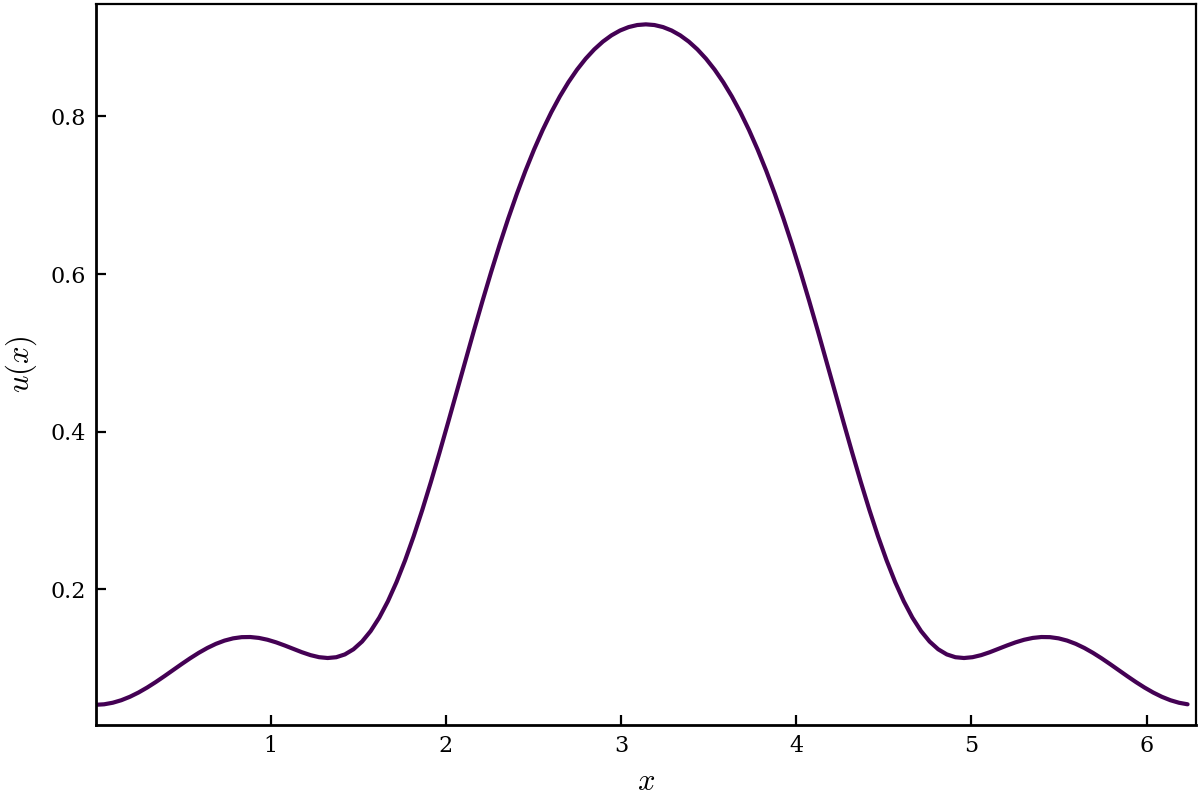} }}
        \caption{Random in-distribution test initial condition for the advection equation.}
        \label{fig:init_in_dist_advection}
    \end{figure}

\begin{figure}[h!]
        \centering
        \subfloat[$M=128$ Fourier solution]{{\includegraphics[width=0.45\textwidth]{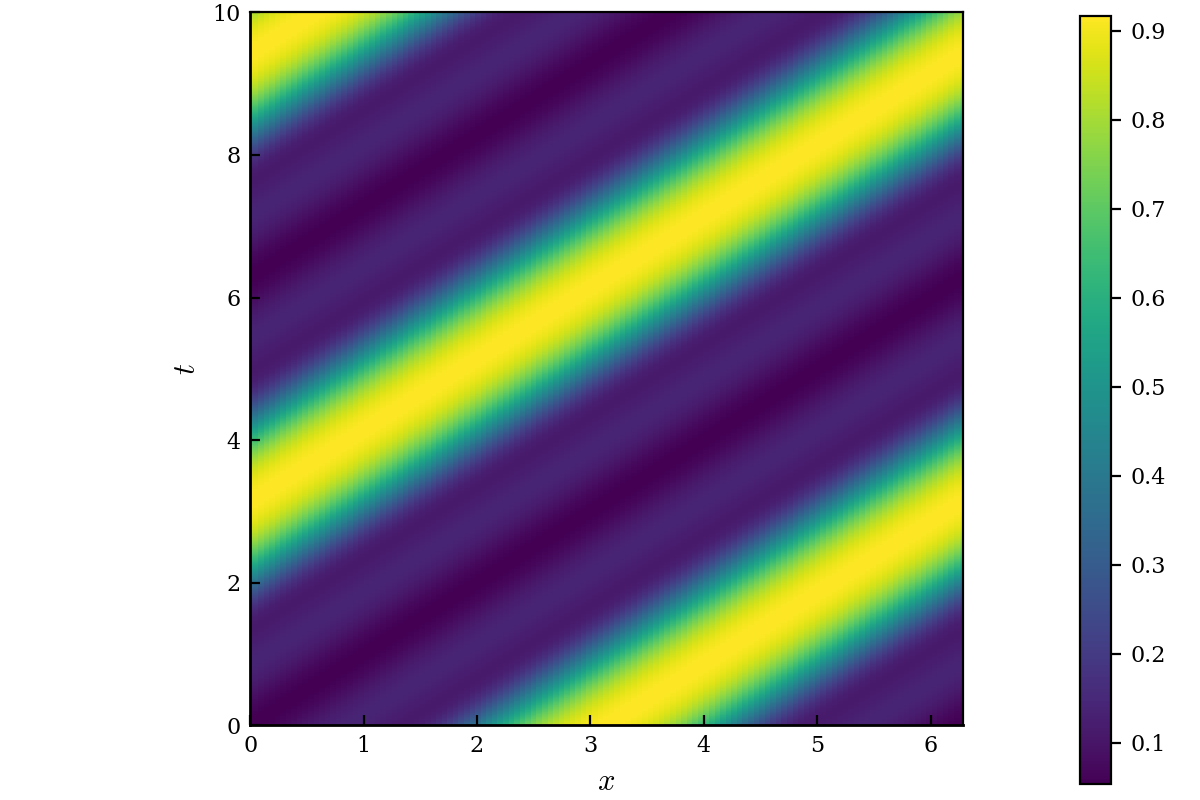} }}
        \qquad
        \subfloat[$r=40$ custom basis function solution]{{\includegraphics[width=0.45\textwidth]{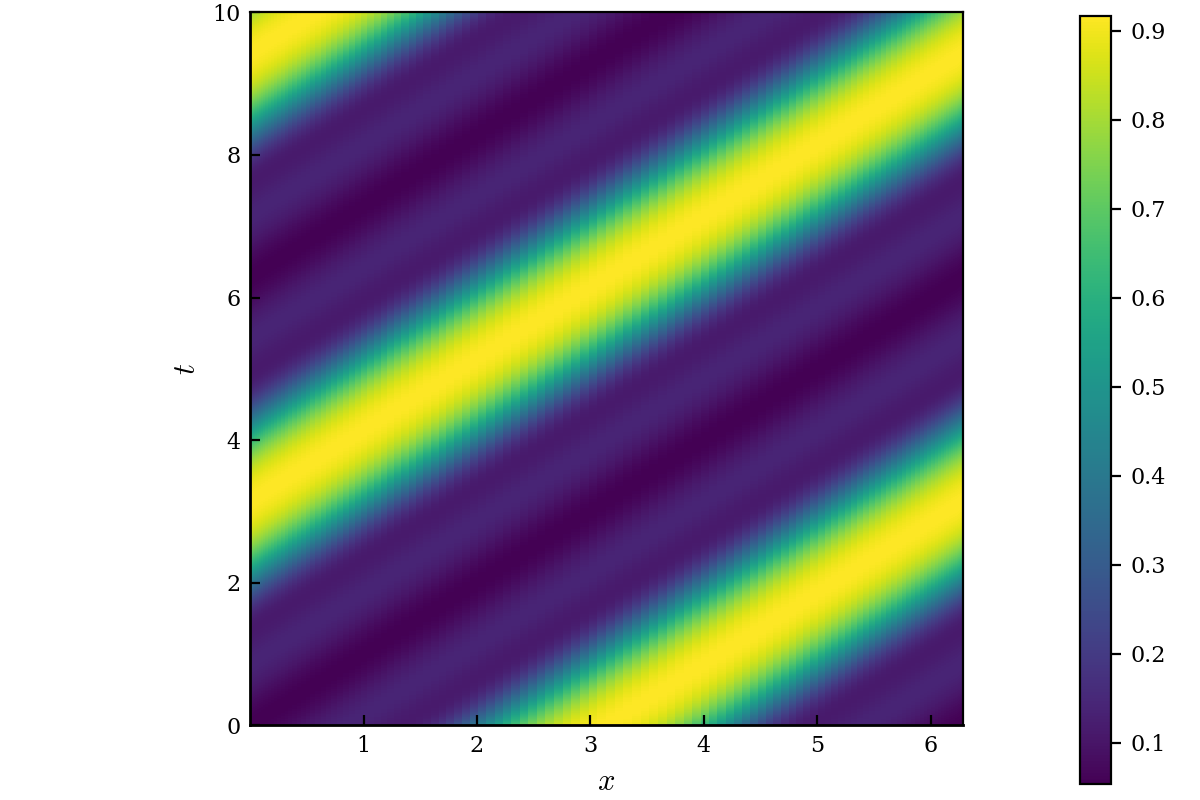} }}
        \qquad
        \caption{Results for the advection equation for $t \in [0,10]$ and the random in training distribution initial condition.}
        \label{fig:advection_comparo_in}
    \end{figure}

\begin{figure}[h!]
        \centering
        \subfloat[$M=128$ Fourier solution]{{\includegraphics[width=0.45\textwidth]{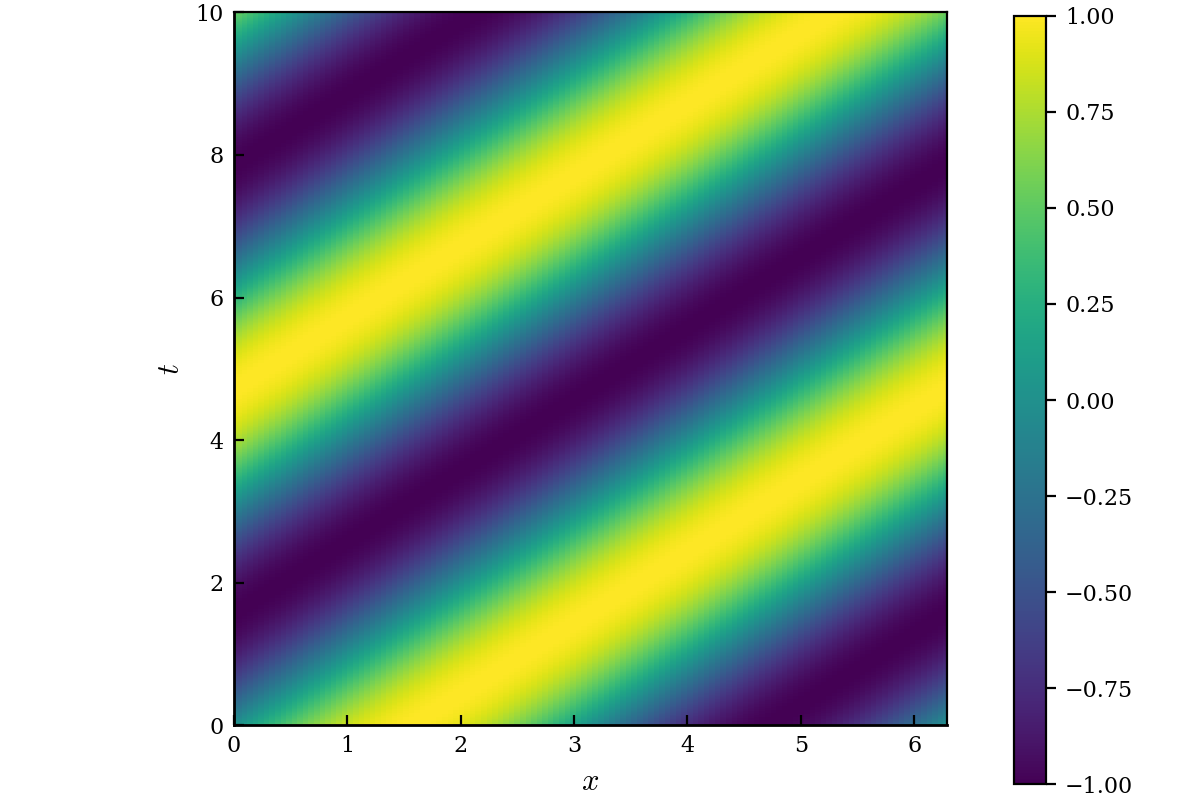} }}
        \qquad
        \subfloat[$r=40$ custom basis function solution]{{\includegraphics[width=0.45\textwidth]{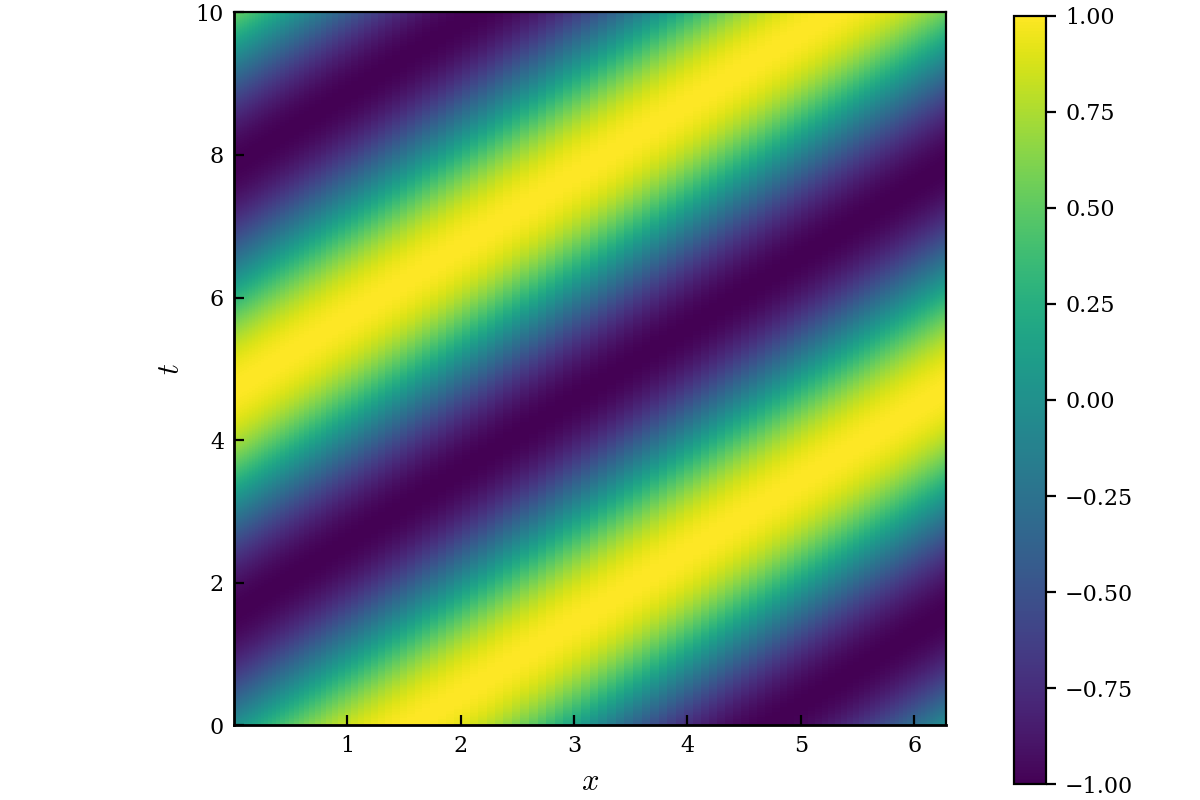} }}
        \qquad
        \caption{Results for the advection equation for $t \in [0,10]$ and the out of training distribution initial condition $u(0,x) = \sin(x)$.}
        \label{fig:advection_comparo_out}
    \end{figure}

\begin{figure}[h!]
        \centering
        \subfloat[In-distribution random initial condition]{{\includegraphics[width=0.45\textwidth]{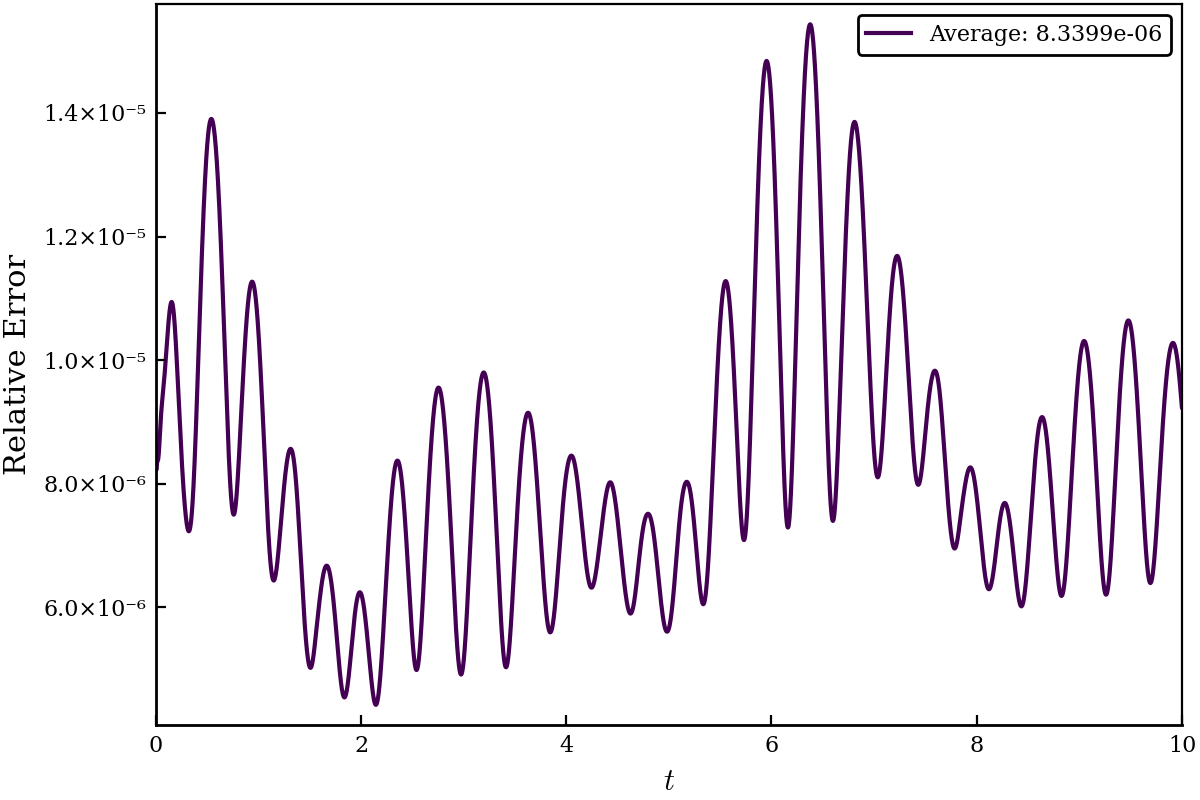} }}
        \qquad
        \subfloat[Out of distribution initial condition: $u(0,x) = \sin(x)$]{{\includegraphics[width=0.45\textwidth]{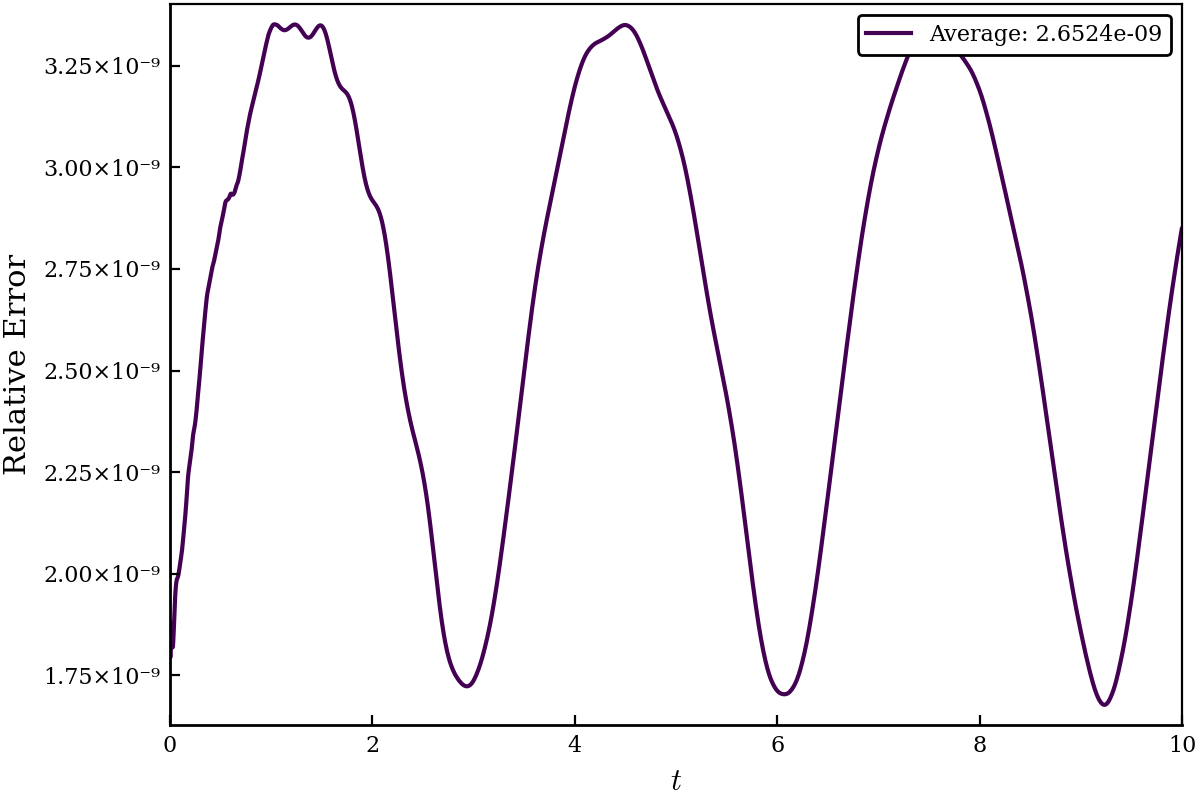} }}
        \qquad
        \caption{Relative errors for the advection equation for $t \in [0,10]$.}
        \label{fig:advection_errors}
    \end{figure}

Results using a $M = 128$ mode Fourier expansion and the custom basis function expansion for the evolution of the advection equation for $t \in [0,10]$ and the in-distribution test initial condition are shown in Figure \ref{fig:advection_comparo_in}. Similarly, results for the out of distribution test initial condition are shown in Figure \ref{fig:advection_comparo_out}. Refer to Figure \ref{fig:init_in_dist_advection} for a plot of the in-distribution and to Figure \ref{fig:init_out_dist} for a plot of the out of distribution test initial conditions. The relative and average errors compared to the $M=128$ Fourier expansion for the two test cases are presented in Figure \ref{fig:advection_errors}.

\subsection{Advection-diffusion equation}

\begin{figure}[h!]
        \centering
        \subfloat[The singular values]{{\includegraphics[width=0.45\textwidth]{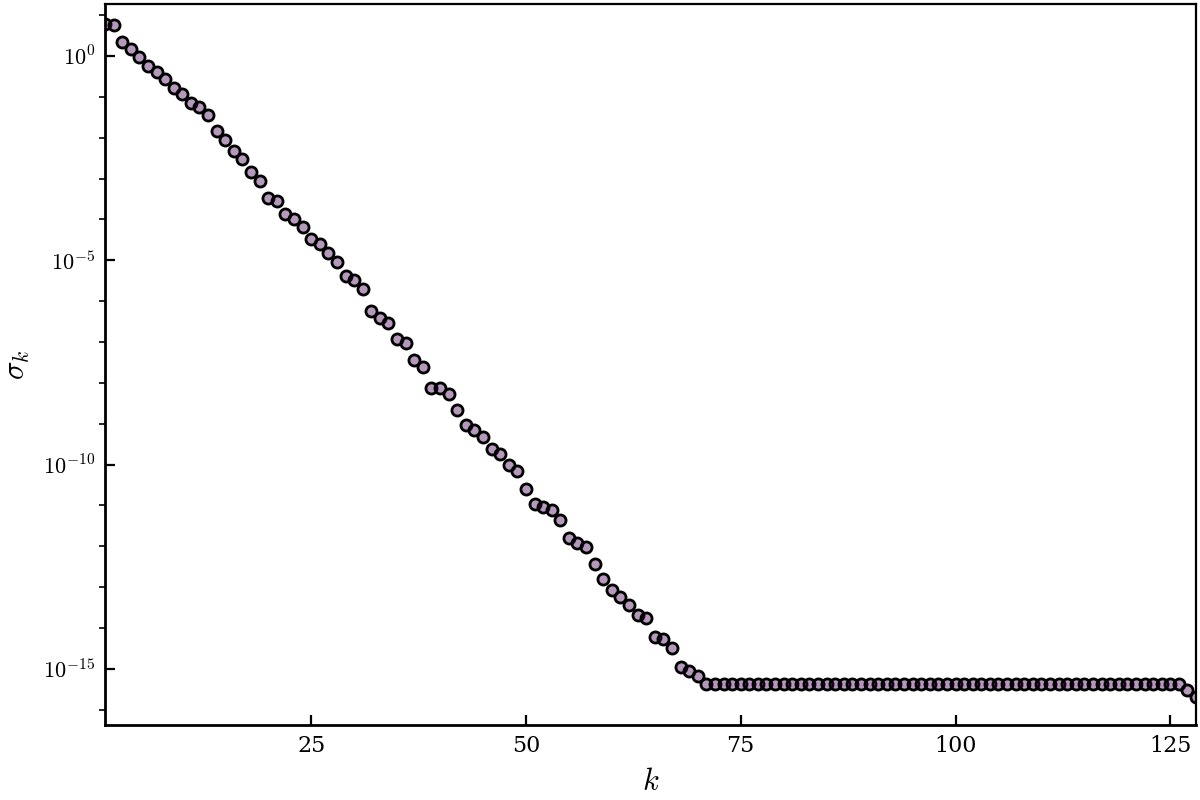} }}
        \qquad
        \subfloat[The first 34 custom basis functions]{{\includegraphics[width=0.45\textwidth]{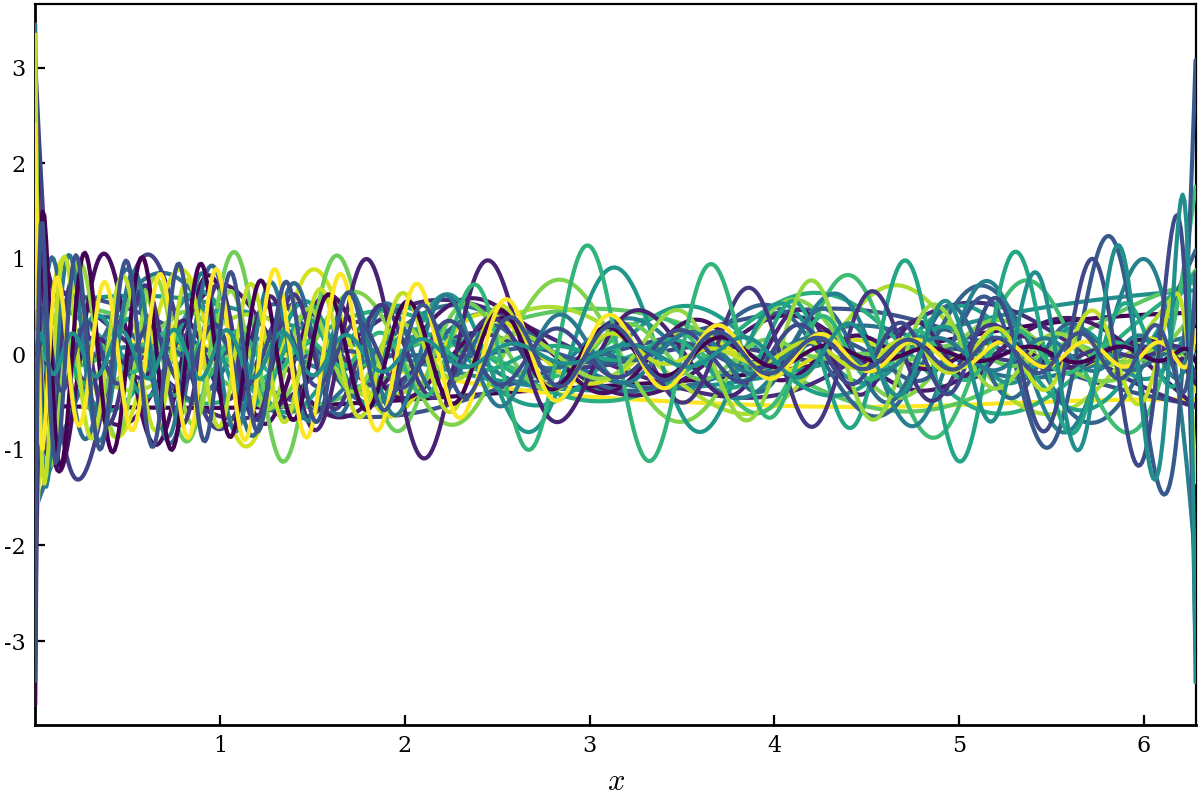} }}
        \qquad
        \caption{Singular values and custom basis functions for the advection-diffusion equation.}
        \label{fig:sv_basis_advection_diffusion}
    \end{figure}

The one-dimensional advection-diffusion equation with periodic boundary conditions is given by
\begin{equation}\label{eqn:advection_diffusion}
    \frac{\partial u}{\partial t} + \frac{\partial u}{\partial x} - \nu \frac{\partial^2 u}{\partial x^2} = 0, \quad x \in [0,2\pi],
\end{equation}
where we take $\nu = 0.1$ throughout. The singular value spectrum and the basis functions satisfying the specified singular value threshold, $\sigma_k > 10^{-7}$, are shown in Figure \ref{fig:sv_basis_advection_diffusion}. The equations for the evolution of the expansion coefficients, $a_m(t), \; m=1,2,\dots,r-b$ are given by,
\begin{equation} 
        \frac{da_m(t)}{dt} = - \sum_{k=1}^r a_k(t) \langle \phi_{m} , \phi'_{k}  \rangle + \nu \sum_{k=1}^r a_k(t) \langle \phi_{m},\phi''_{k} \rangle,
\end{equation}
where $r = 34$ and the primes denote differentiation in space.

\begin{figure}[h!]
        \centering
        {{\includegraphics[width=0.45\textwidth]{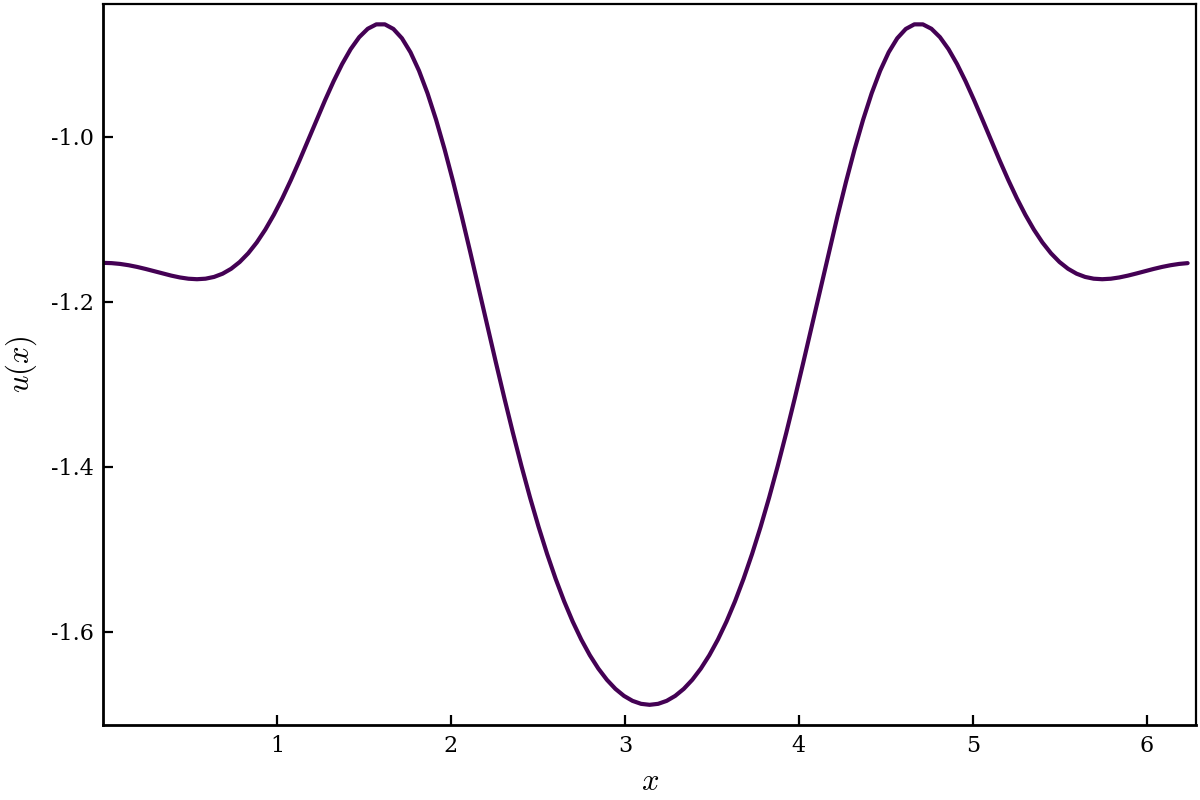} }}
        \caption{Random in-distribution test initial condition for the advection-diffusion equation.}
        \label{fig:init_in_dist_advection_diffusion}
    \end{figure}

\begin{figure}[h!]
        \centering
        \subfloat[$M=128$ Fourier solution]{{\includegraphics[width=0.45\textwidth]{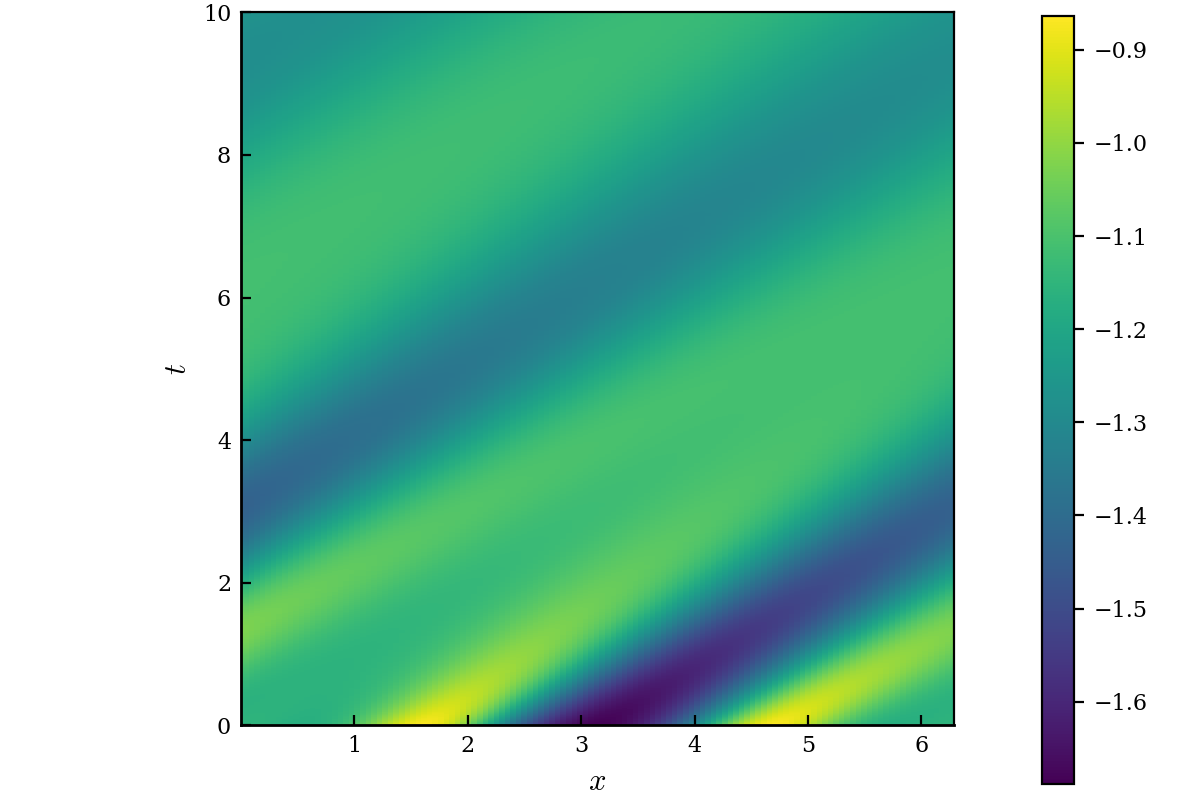} }}
        \qquad
        \subfloat[$r=34$ custom basis function solution]{{\includegraphics[width=0.45\textwidth]{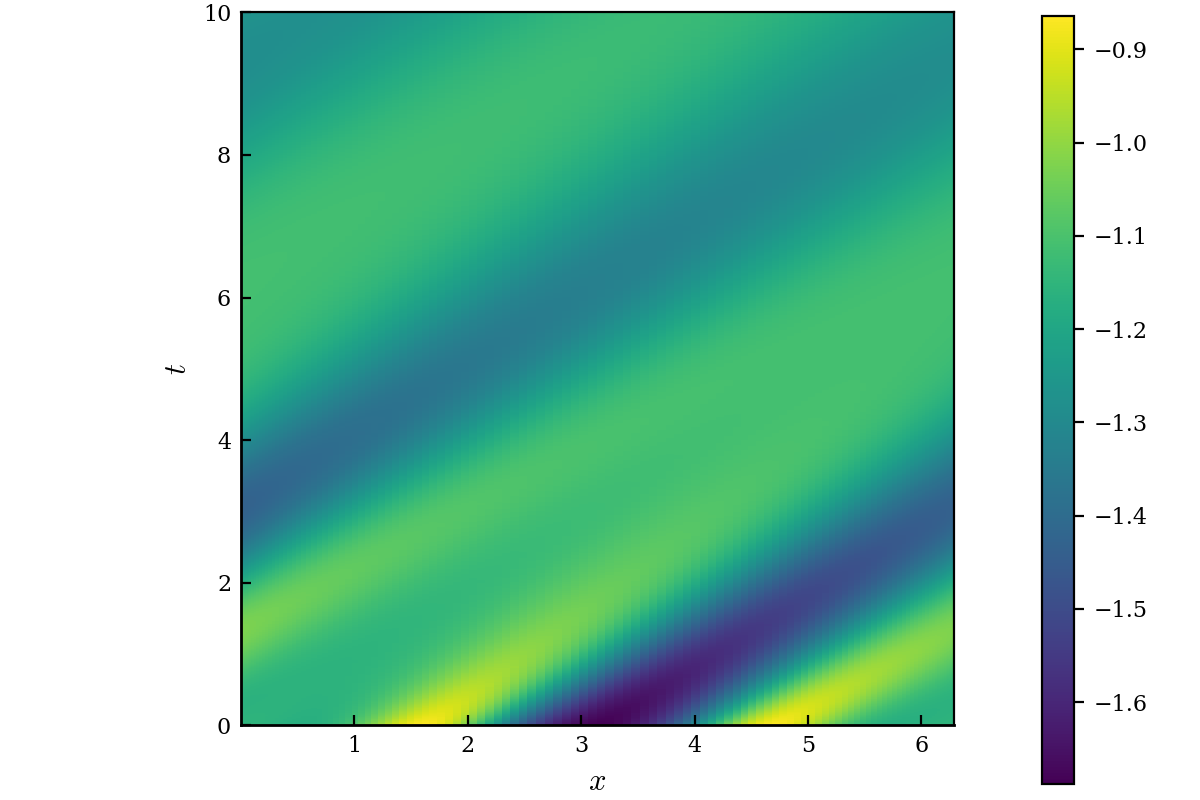} }}
        \qquad
        \caption{Results for the advection-diffusion equation for $t \in [0,10]$ and the random in training distribution initial condition.}
        \label{fig:advection_diffusion_comparo_in}
    \end{figure}

\begin{figure}[h!]
        \centering
        \subfloat[$M=128$ Fourier solution]{{\includegraphics[width=0.45\textwidth]{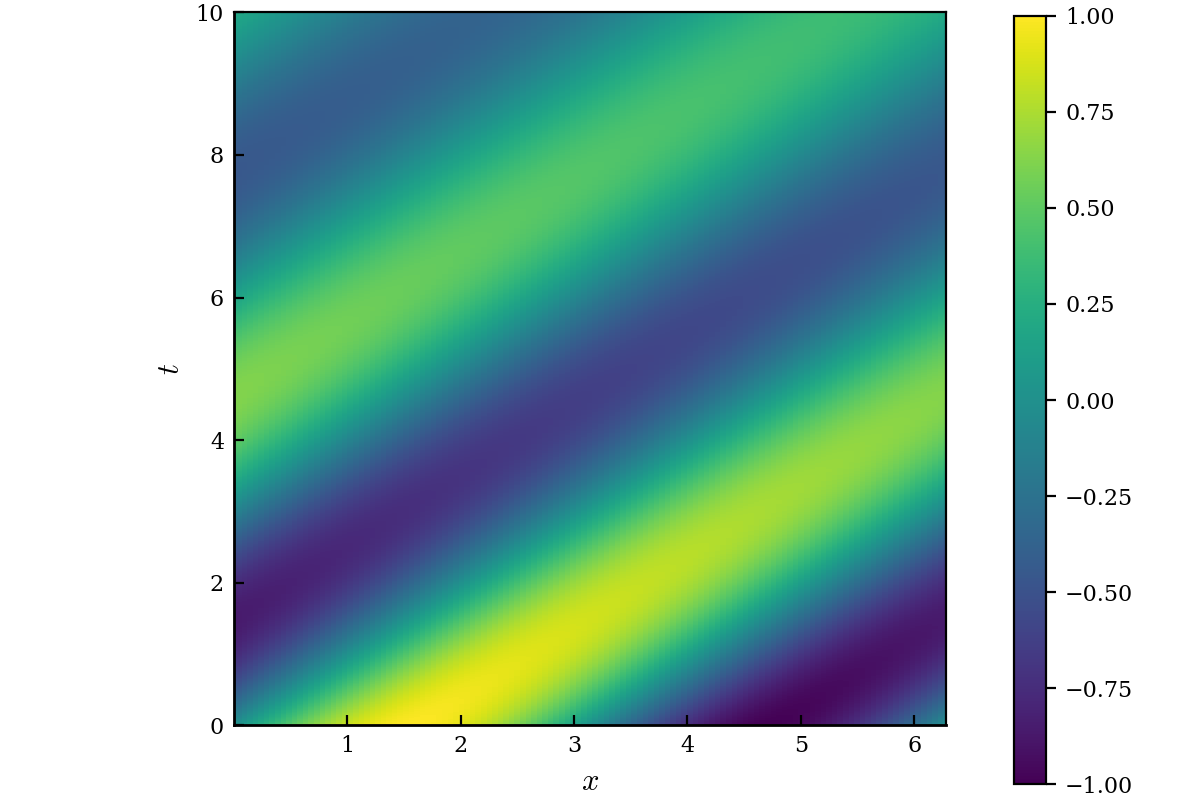} }}
        \qquad
        \subfloat[$r=34$ custom basis function solution]{{\includegraphics[width=0.45\textwidth]{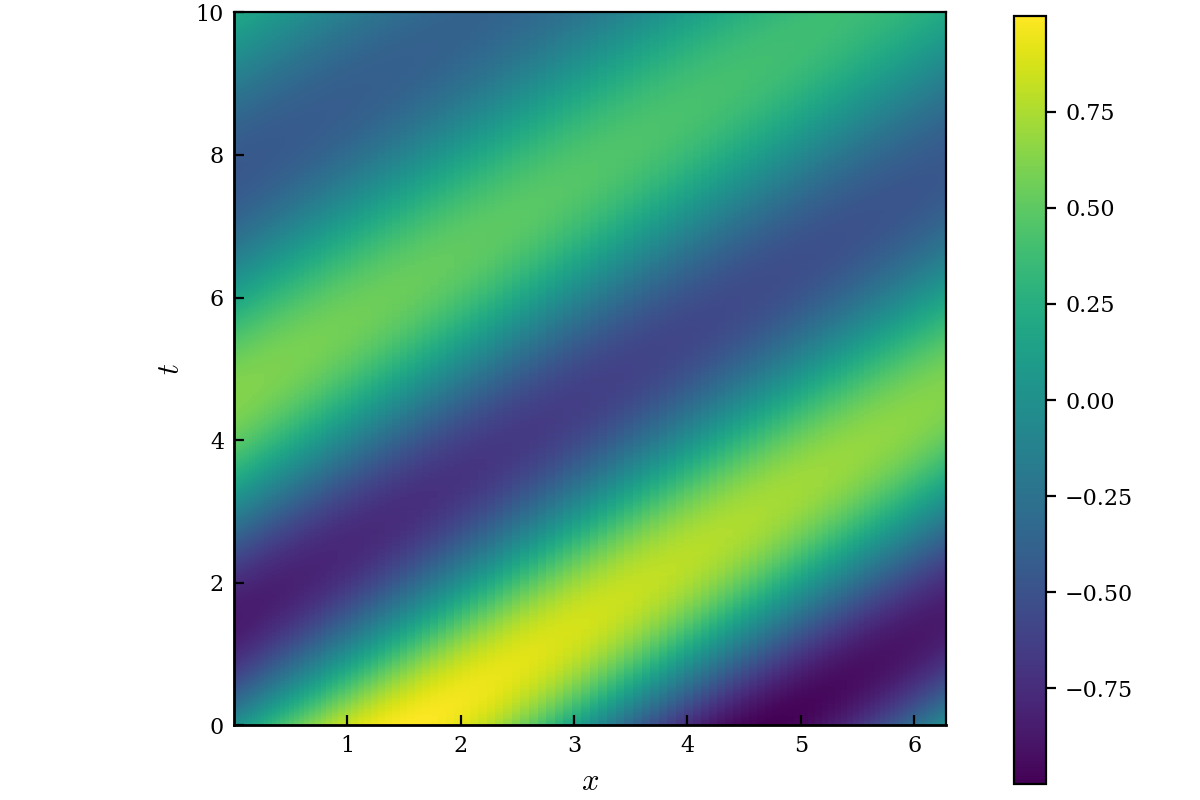} }}
        \qquad
        \caption{Results for the advection-diffusion equation for $t \in [0,10]$ and the out of training distribution initial condition $u(0,x) = \sin(x)$.}
        \label{fig:advection_diffusion_comparo_out}
    \end{figure}

\begin{figure}[h!]
        \centering
        \subfloat[In-distribution random initial condition]{{\includegraphics[width=0.45\textwidth]{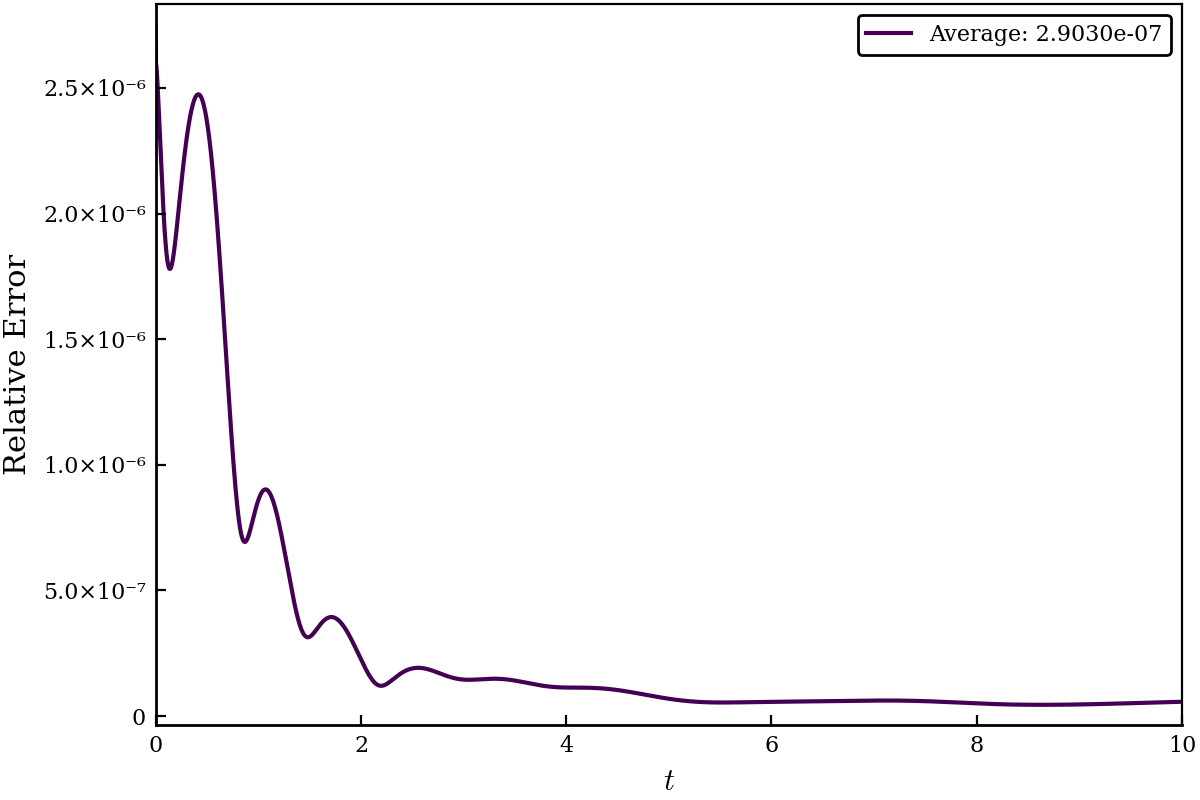} }}
        \qquad
        \subfloat[Out of distribution initial condition: $u(0,x) = \sin(x)$]{{\includegraphics[width=0.45\textwidth]{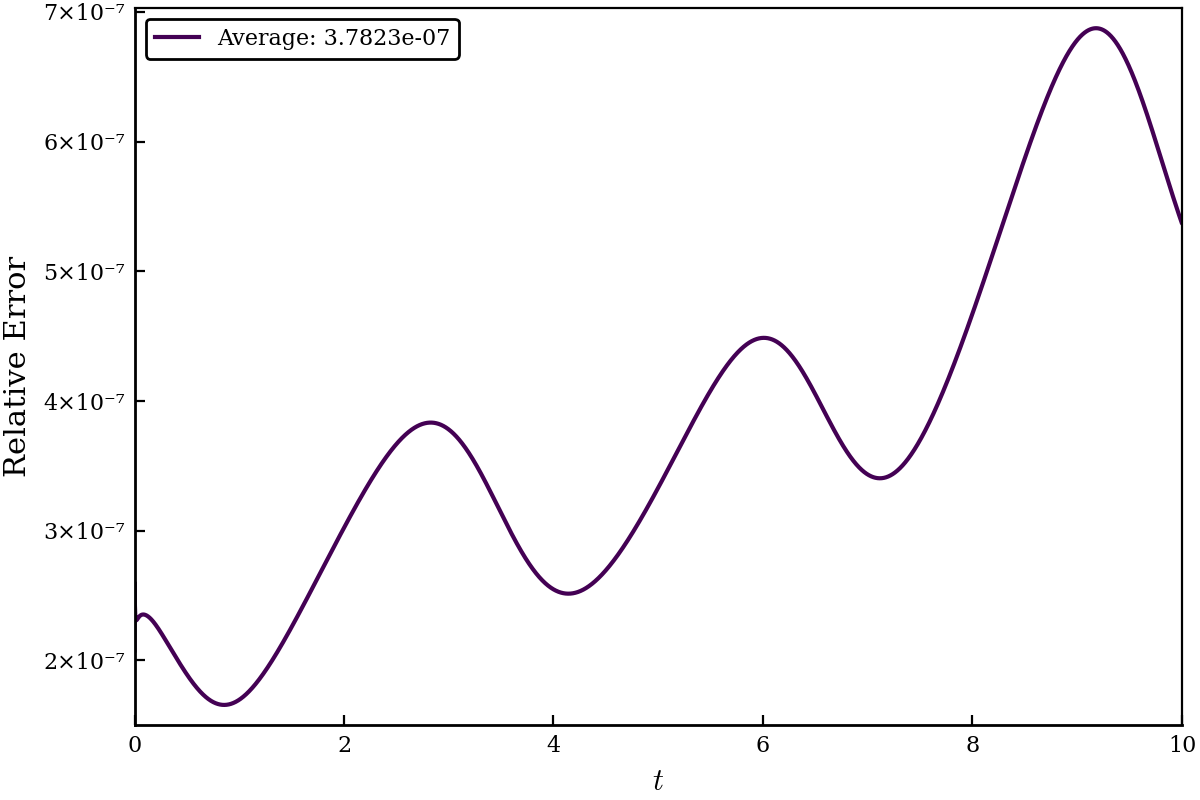} }}
        \qquad
        \caption{Relative errors for the advection-diffusion equation for $t \in [0,10]$.}
        \label{fig:advection_diffusion_errors}
    \end{figure}

Results using a $M = 128$ mode Fourier expansion and the custom basis function expansion for the evolution of the advection-diffusion equation for $t \in [0,10]$ and the in-distribution test initial condition are shown in Figure \ref{fig:advection_diffusion_comparo_in}. See Figure \ref{fig:init_in_dist_advection_diffusion} for a plot of the in-distribution test initial condition. Similarly, results for the out of distribution test initial condition are shown in Figure \ref{fig:advection_diffusion_comparo_out}. The relative and average errors compared to the $M=128$ Fourier expansion for the two test cases are presented in Figure \ref{fig:advection_diffusion_errors}.

\subsection{Viscous Burgers equation}

\begin{figure}[h!]
        \centering
        \subfloat[The singular values]{{\includegraphics[width=0.45\textwidth]{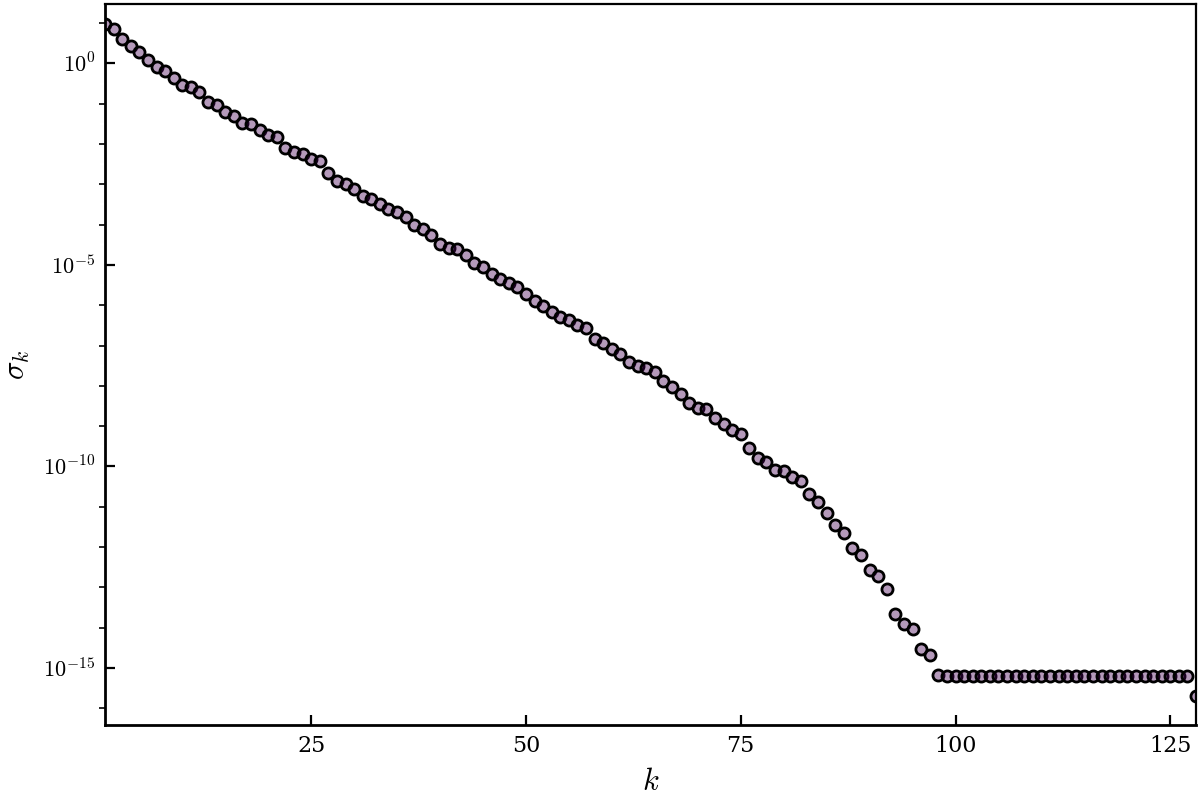} }}
        \qquad
        \subfloat[The first 58 custom basis functions]{{\includegraphics[width=0.45\textwidth]{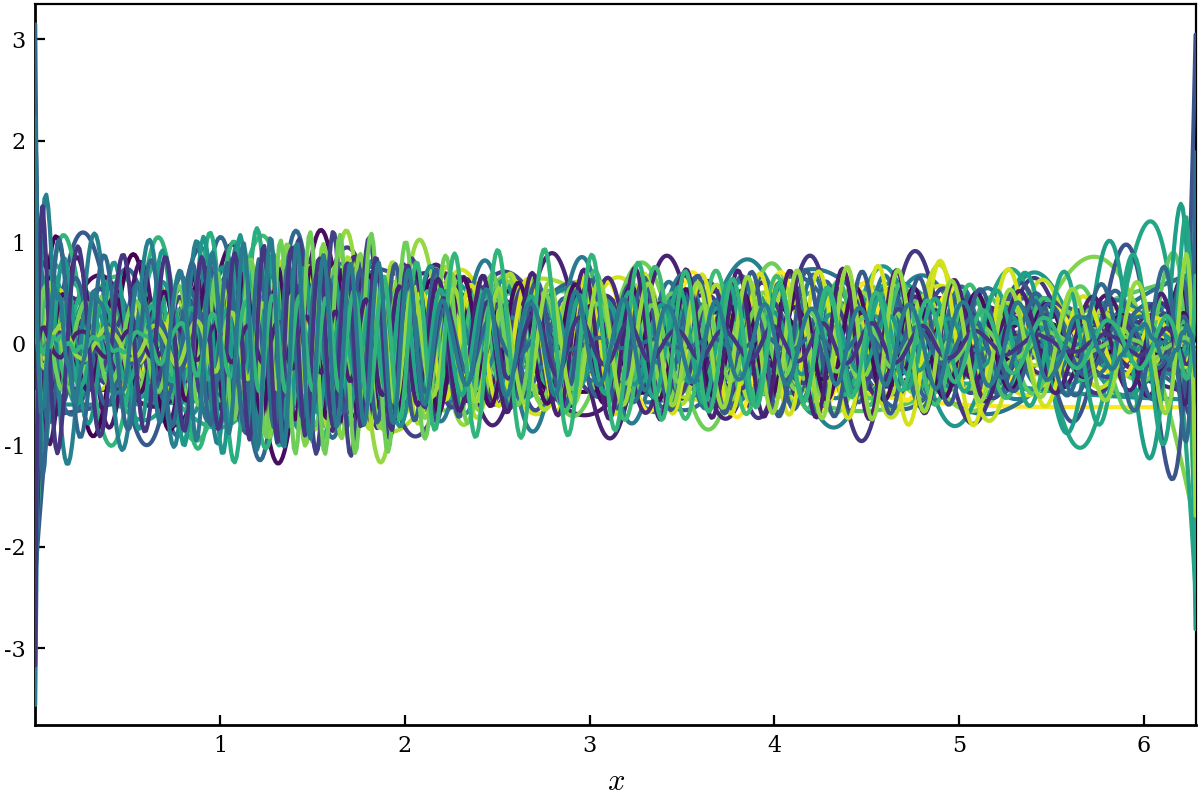} }}
        \qquad
        \caption{Singular values and custom basis functions for the viscous Burgers equation.}
        \label{fig:sv_basis_viscous_burgers}
    \end{figure}

The one-dimensional viscous Burgers equation with periodic boundary conditions is given by
\begin{equation}\label{eqn:viscous_burgers}
    \frac{\partial u}{\partial t} + u \frac{\partial u}{\partial x} - \nu \frac{\partial^2 u}{\partial x^2} = 0, \quad x \in [0,2\pi], 
\end{equation}
where we again set $\nu = 0.1$ throughout. The singular value spectrum and the basis functions satisfying the specified singular value threshold, $\sigma_k > 10^{-7}$, are shown in Figure \ref{fig:sv_basis_viscous_burgers}.
The equations for the evolution of the expansion coefficients, $a_m(t), \; m=1,2,\dots,r-b$ are given by,
\begin{equation}\label{eqn:custom_viscous_burgers}
        \frac{da_m(t)}{dt} = - \sum_{k=1}^r \sum_{l=1}^r a_k(t) a_l(t) \langle \phi_{m}, \phi_{k} \phi'_{l}\rangle + \nu \sum_{k=1}^r a_k(t) \langle \phi_{m},\phi''_{k} \rangle,
\end{equation}
where $r = 58$ and the primes denote differentiation in space.

Direct computation of the triple product integral $\langle \phi_{m}, \phi_{k}\phi'_{l} \rangle$ and a pseudospectral transform were tested for computing the nonlinear term in \eqref{eqn:custom_viscous_burgers}. For the pseudospectral transform, the nonlinear term is computed by transforming the solution from the custom basis function space to real space, performing the multiplication of the two terms, and transforming back to the custom basis function space. 
For all results, the triple product integral was computed directly for the nonlinear term. In Section \ref{sec:discussion_future_work} we provide additional comments on using DeepONets to develop a fast custom basis function inverse transform, which would accelerate the pseudospectral transform. 

\begin{figure}[h!]
        \centering
        {{\includegraphics[width=0.45\textwidth]{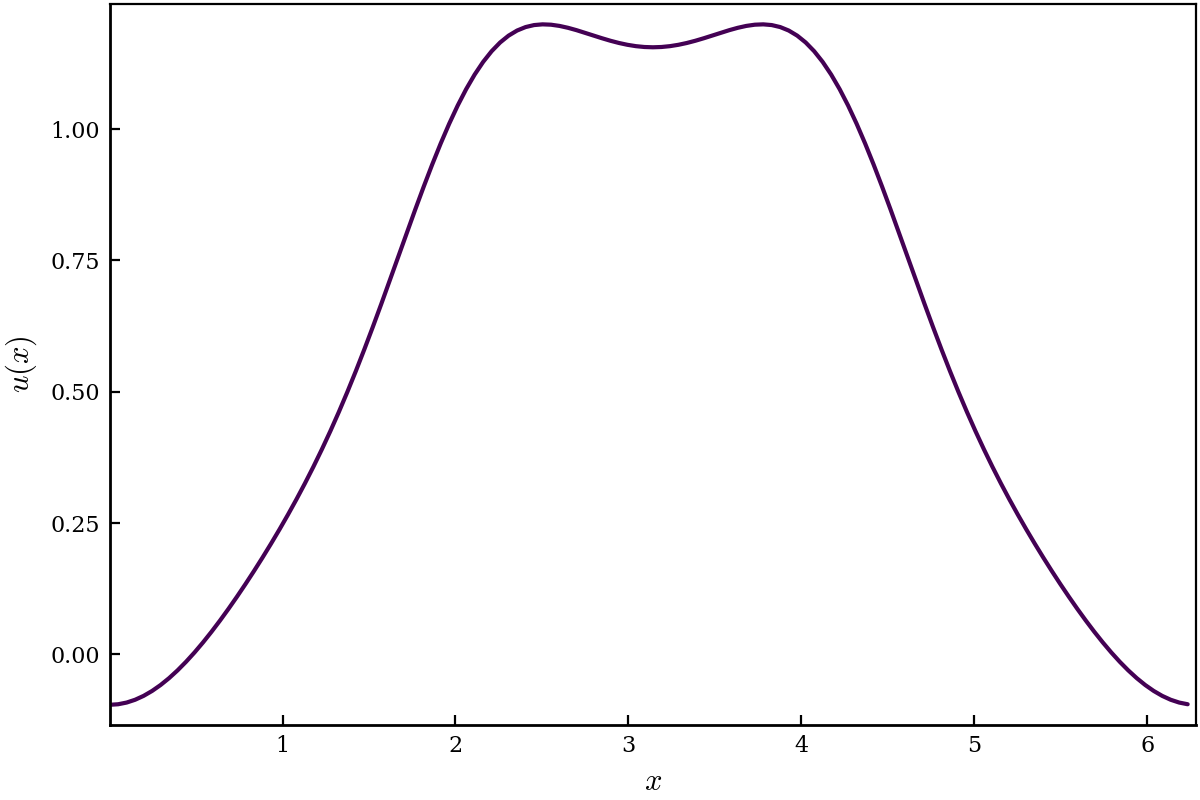} }}
        \caption{Random in-distribution test initial condition for the viscous Burgers equation.}
        \label{fig:init_in_dist_viscous_burgers}
    \end{figure}

\begin{figure}[h!]
        \centering
        \subfloat[$M=128$ Fourier solution]{{\includegraphics[width=0.45\textwidth]{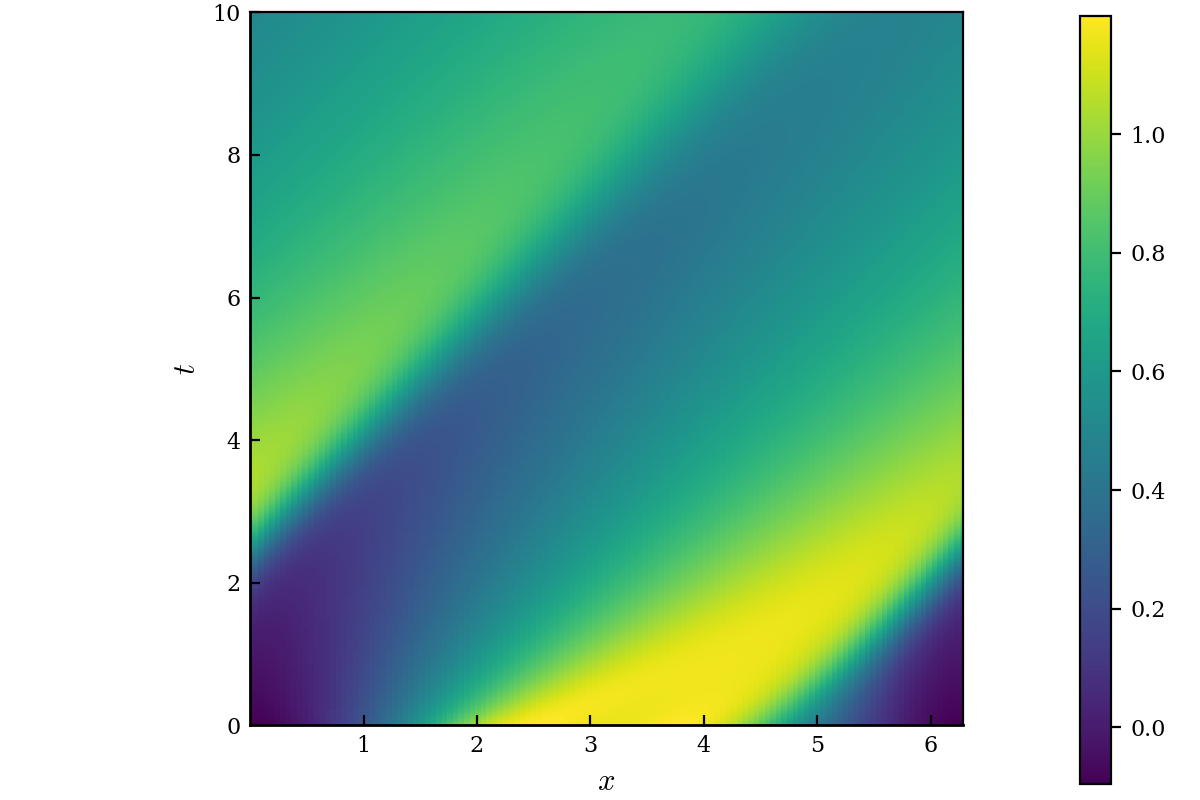} }}
        \qquad
        \subfloat[$r=58$ custom basis function solution]{{\includegraphics[width=0.45\textwidth]{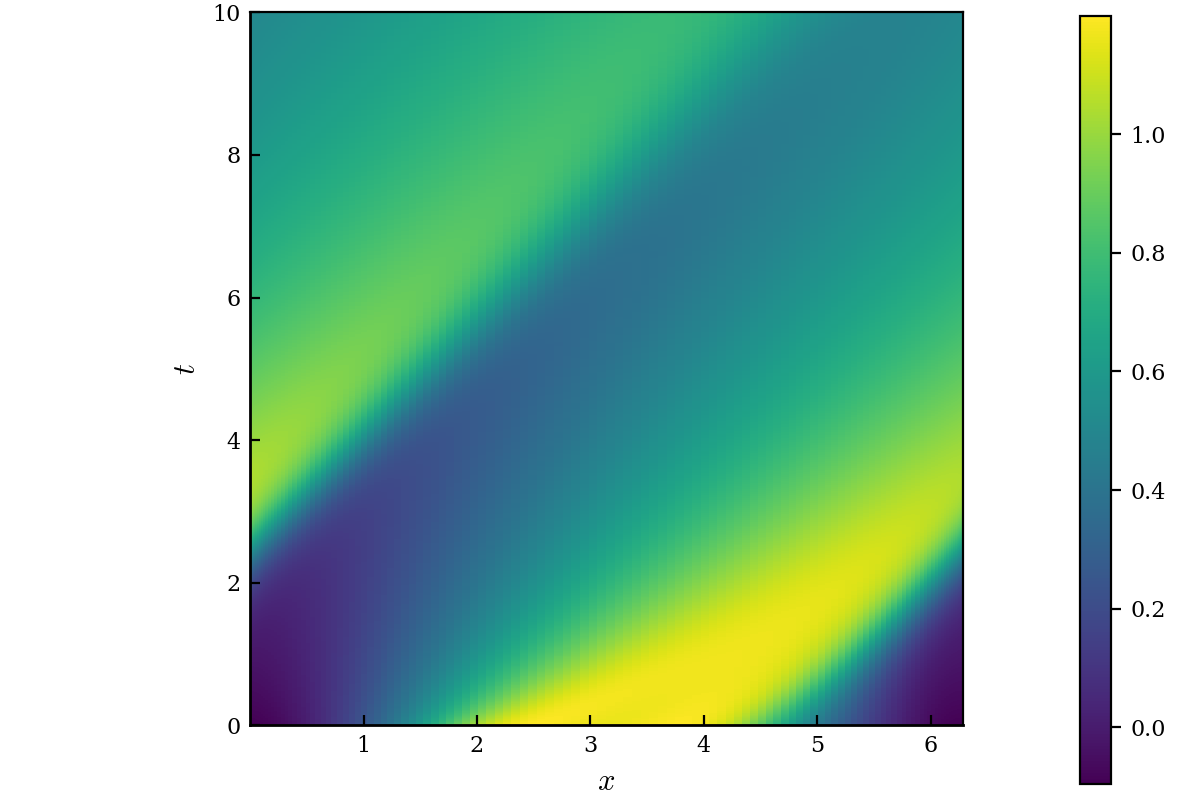} }}
        \qquad
        \caption{Results for the viscous Burgers equation for $t \in [0,10]$ and the random in training distribution initial condition.}
        \label{fig:viscous_burgers_comparo_in}
    \end{figure}

\begin{figure}[h!]
        \centering
        \subfloat[$M=128$ Fourier solution]{{\includegraphics[width=0.45\textwidth]{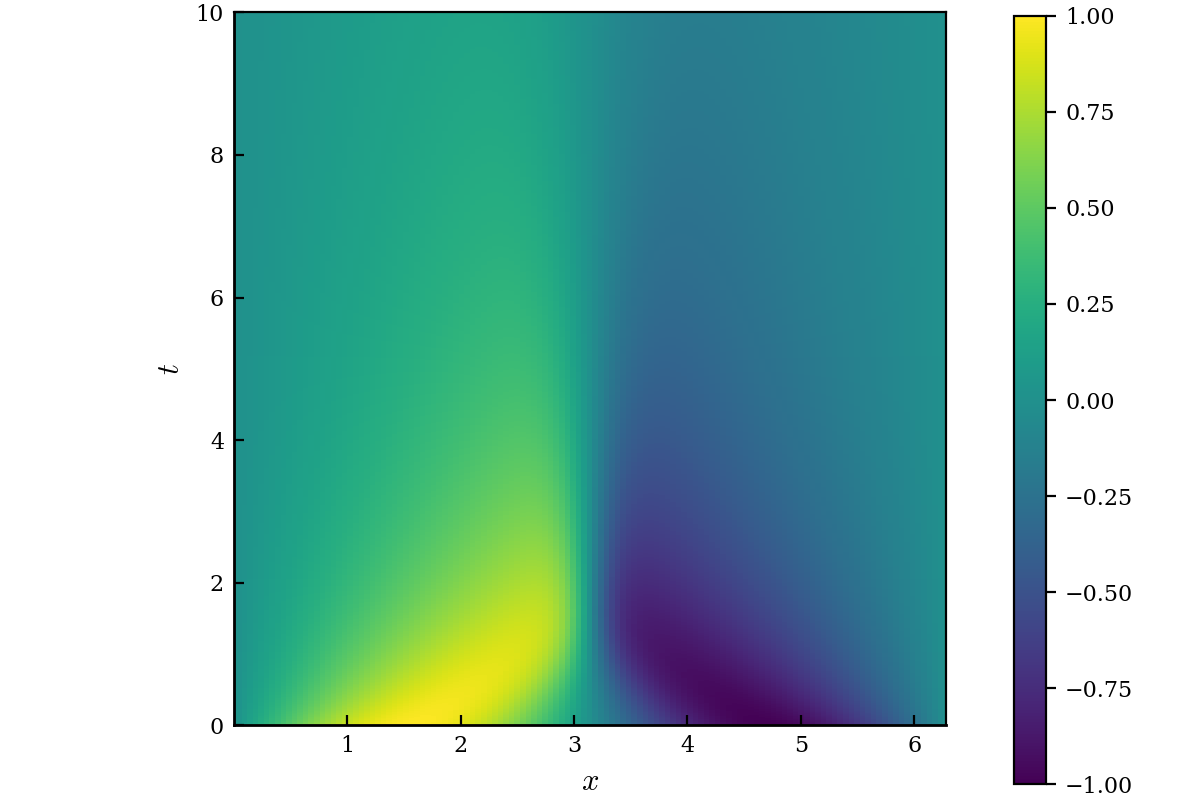} }}
        \qquad
        \subfloat[$r=58$ custom basis function solution]{{\includegraphics[width=0.45\textwidth]{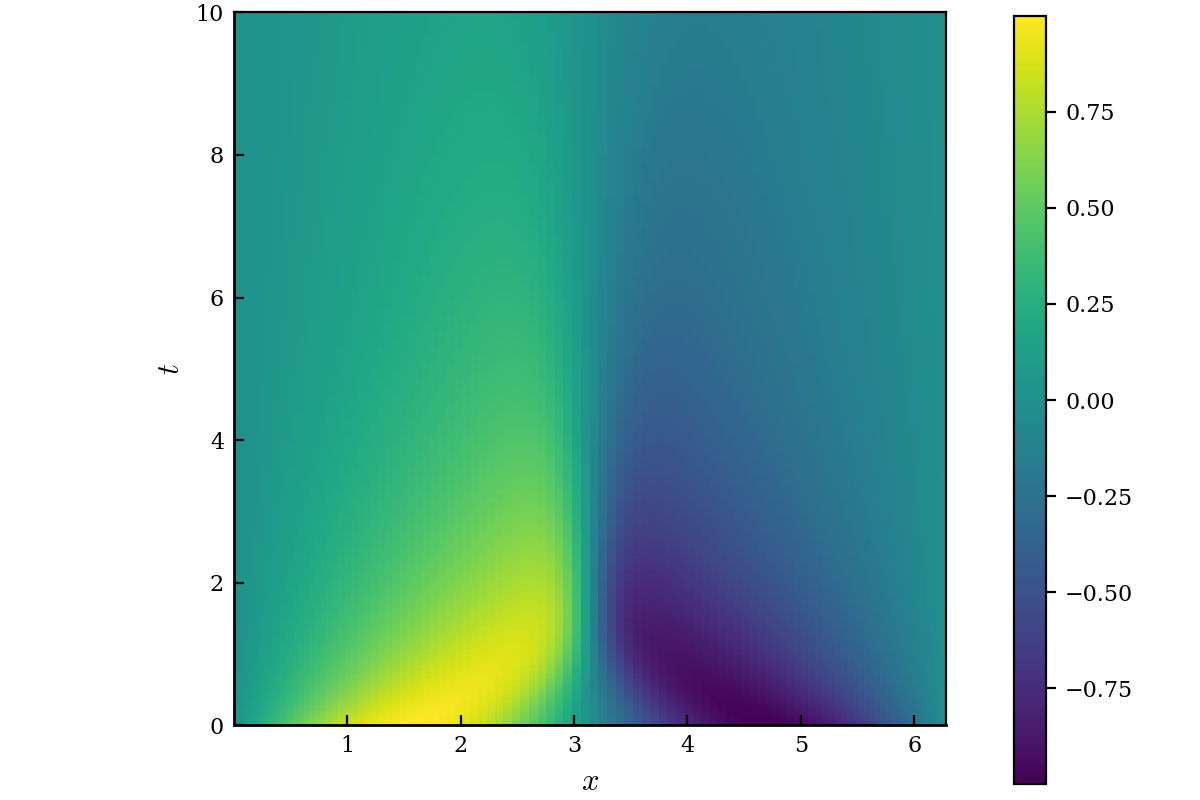} }}
        \qquad
        \caption{Results for the viscous Burgers equation for $t \in [0,10]$ and the out of training distribution initial condition $u(0,x) = \sin(x)$.}
        \label{fig:viscous_burgers_comparo_out}
    \end{figure}

\begin{figure}[h!]
        \centering
        \subfloat[ In-distribution random initial condition]{{\includegraphics[width=0.45\textwidth]{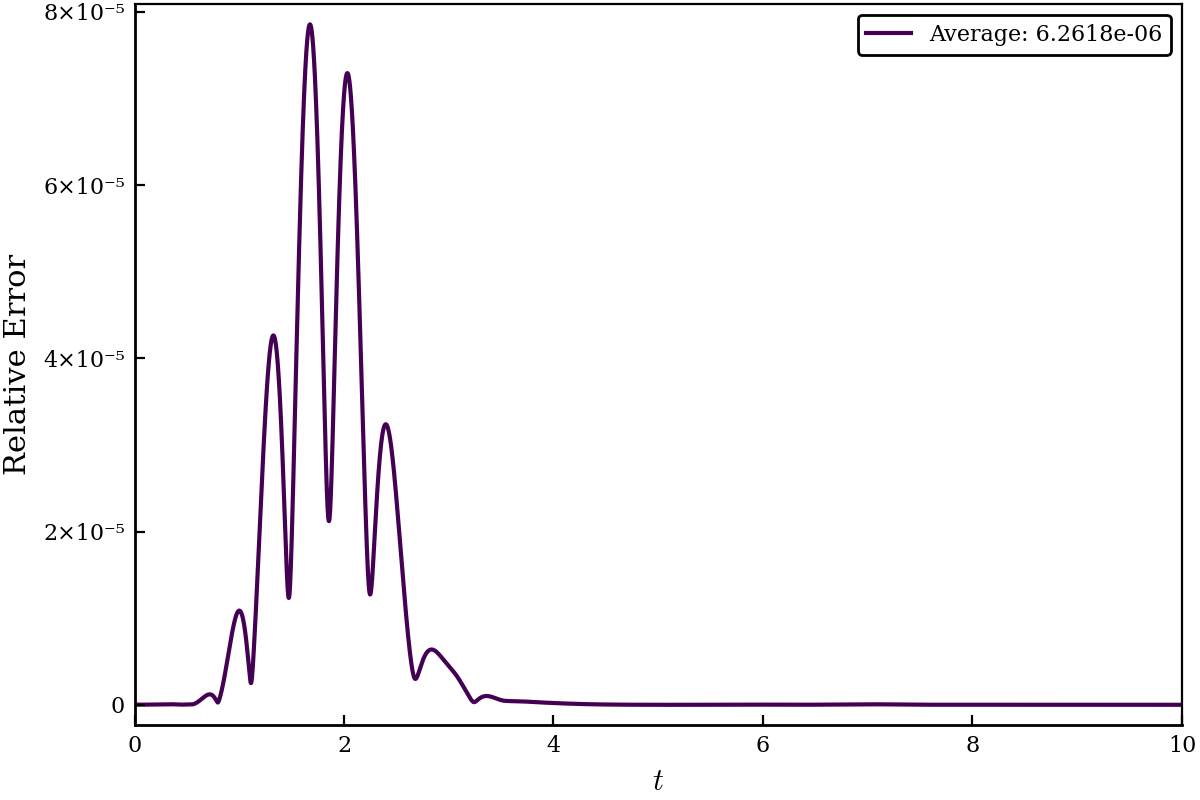} }}
        \qquad
        \subfloat[Out of distribution initial condition: $u(0,x) = \sin(x)$]{{\includegraphics[width=0.45\textwidth]{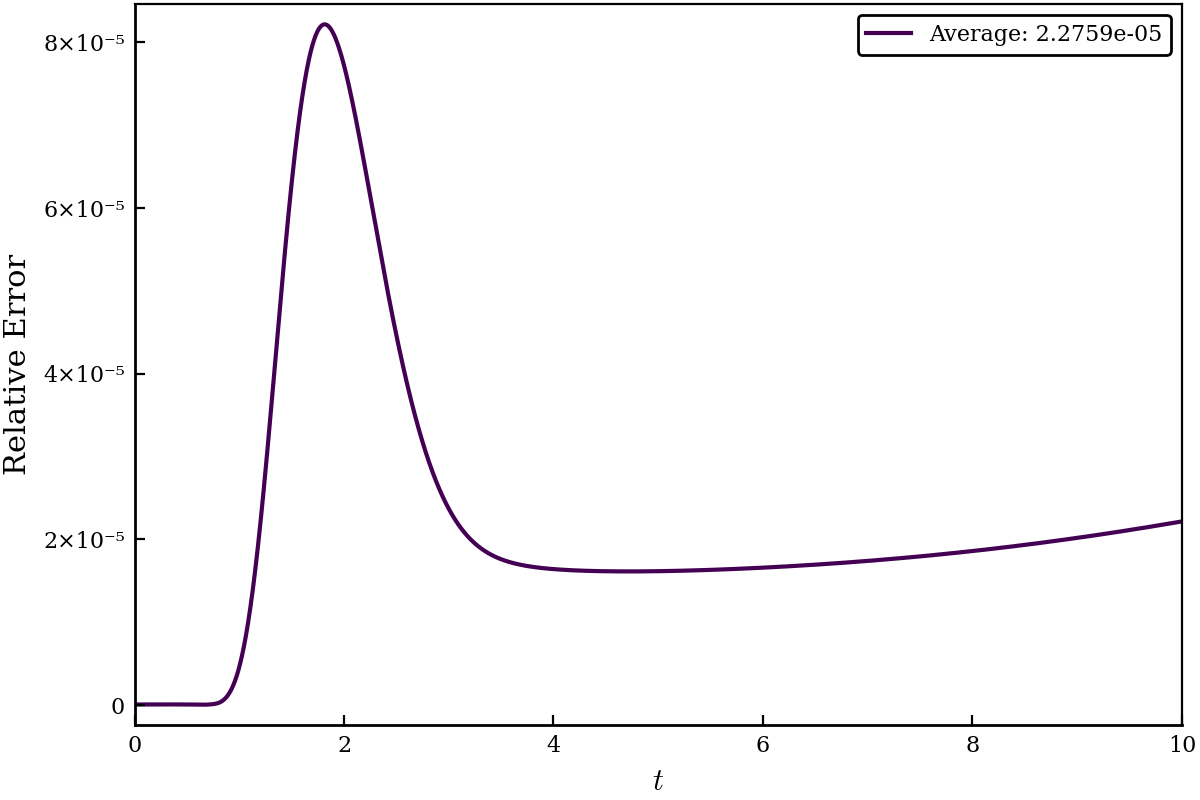} }}
        \qquad
        \caption{Relative errors for the viscous Burgers equation for $t \in [0,10]$.}
        \label{fig:viscous_burgers_errors}
    \end{figure}

Results using a $M = 128$ mode Fourier expansion and the custom basis function expansion for the evolution of the viscous Burgers equation for $t \in [0,10]$ and the in-distribution test initial condition are shown in Figure \ref{fig:viscous_burgers_comparo_in}. See Figure \ref{fig:init_in_dist_viscous_burgers} for a plot of the in-distribution test initial condition. Similarly, results for the out of distribution test initial condition are shown in Figure \ref{fig:viscous_burgers_comparo_out}. The relative and average errors compared to the $M=128$ Fourier expansion for the two test cases are presented in Figure \ref{fig:viscous_burgers_errors}.

\subsection{Inviscid Burgers equation}

The one-dimensional inviscid Burgers equation with periodic boundary conditions is given by
\begin{equation}\label{eqn:inviscid_burgers}
    \frac{\partial u}{\partial t} + u \frac{\partial u}{\partial x} = 0, \;\; x \in [0,2\pi].
\end{equation}

\begin{figure}[h!]
        \centering
        \subfloat[The singular values]{{\includegraphics[width=0.45\textwidth]{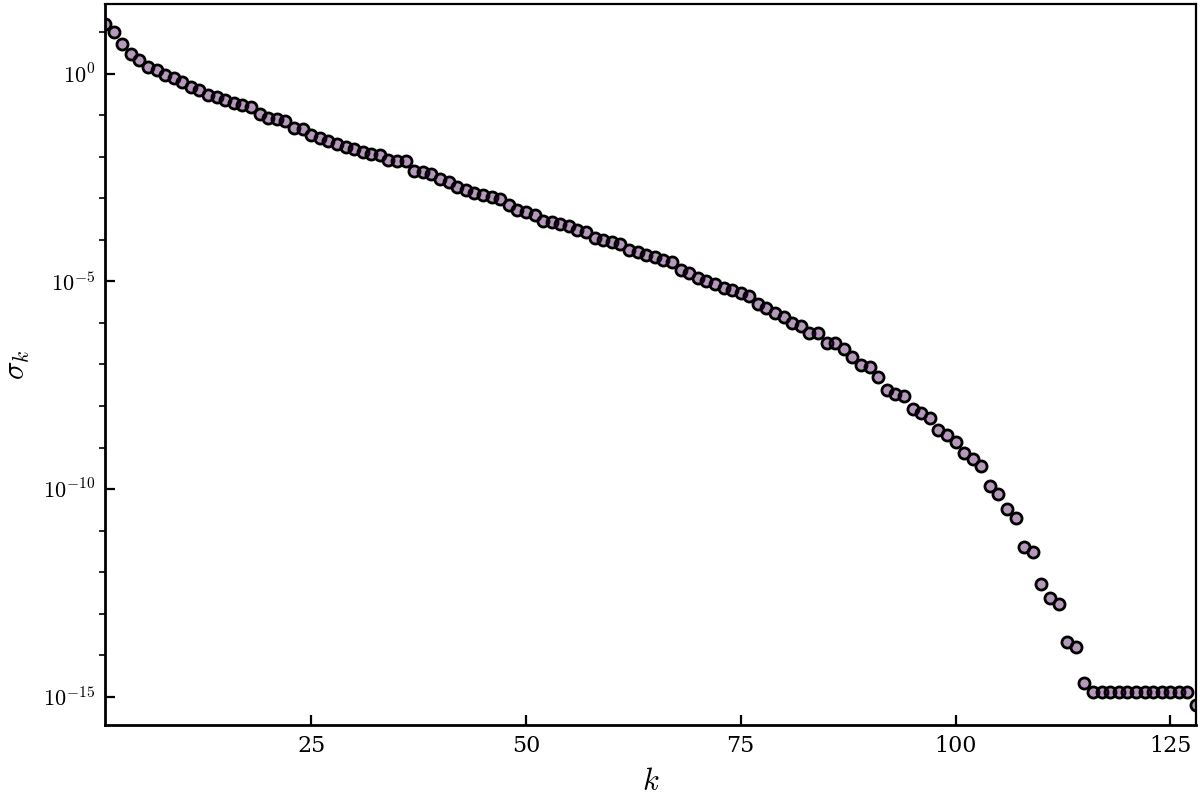} }}
        \qquad
        \subfloat[The first 108 custom basis functions]{{\includegraphics[width=0.45\textwidth]{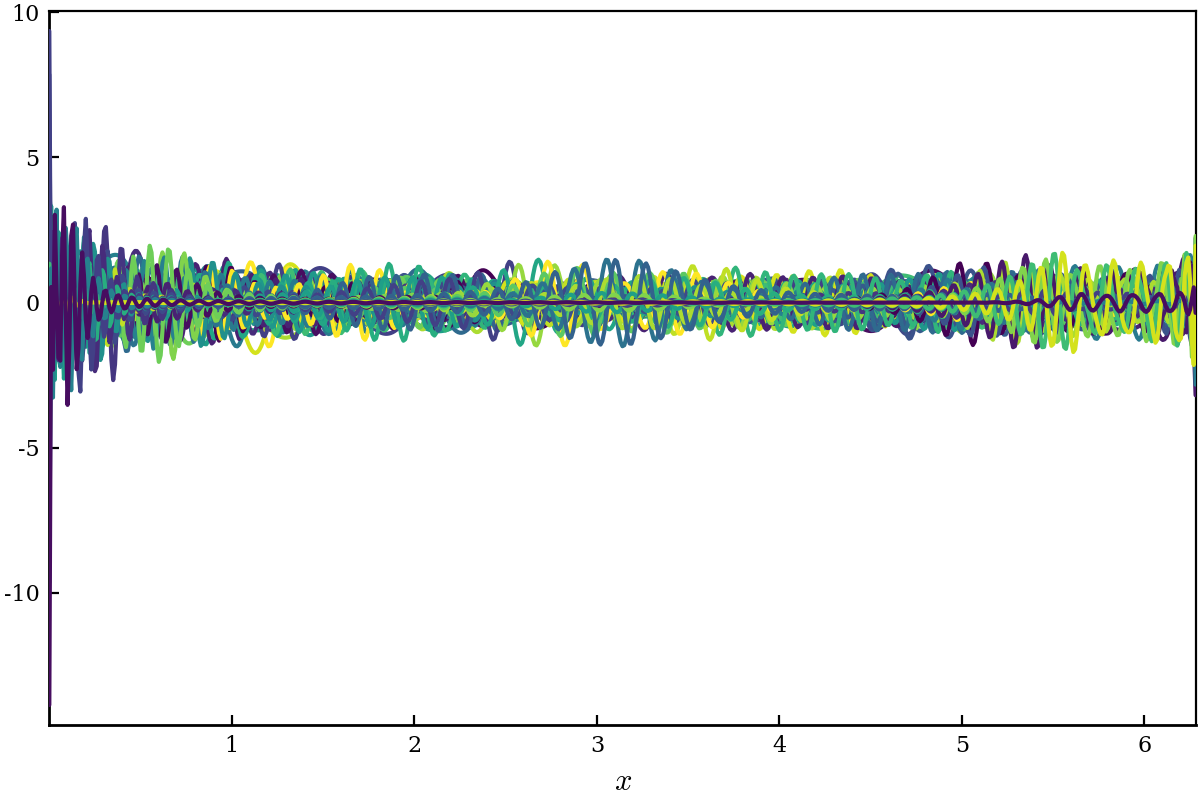} }}
        \qquad
        \caption{Singular values and custom basis functions for the inviscid Burgers equation.}
        \label{fig:sv_basis_inviscid_burgers}
    \end{figure}

The equations for the evolution of the expansion coefficients then are $m=1,2,\dots,r-b$ are given by,
\begin{equation} 
        \frac{da_m(t)}{dt} = - \sum_{k=1}^r \sum_{l=1}^r a_k(t) a_l(t) \langle \phi_{m},\phi_{k}\phi'_{l} \rangle,
\end{equation}
where $r = 108$ and the primes denote differentiation in space. The singular value spectrum and the basis functions satisfying the specified singular value threshold, $\sigma_k > 10^{-12}$, are shown in Figure \ref{fig:sv_basis_inviscid_burgers}.

\begin{figure}[h!]
        \centering
        {{\includegraphics[width=0.45\textwidth]{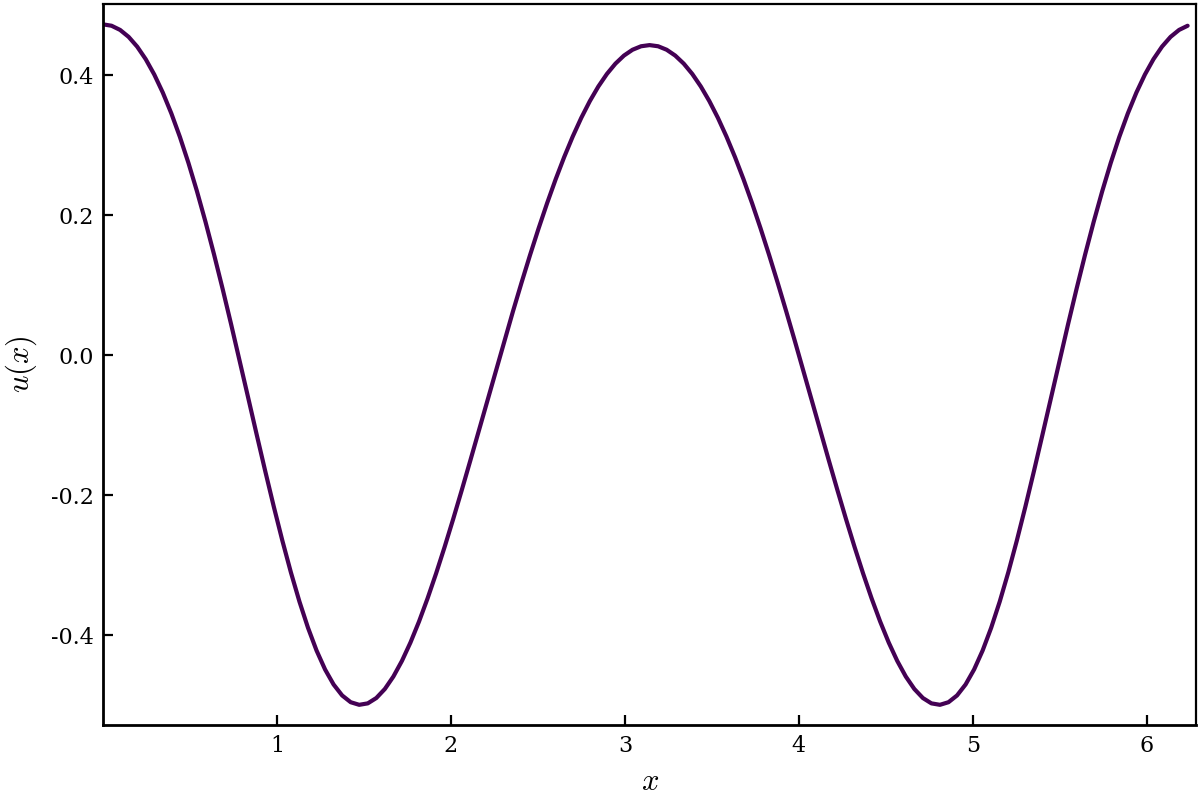} }}
        \caption{Random in-distribution test initial condition for the inviscid Burgers equation.}
        \label{fig:init_in_dist_inviscid_burgers}
    \end{figure}

\begin{figure}[h!]
        \centering
        \subfloat[The MUSCL solution]{{\includegraphics[width=0.45\textwidth]{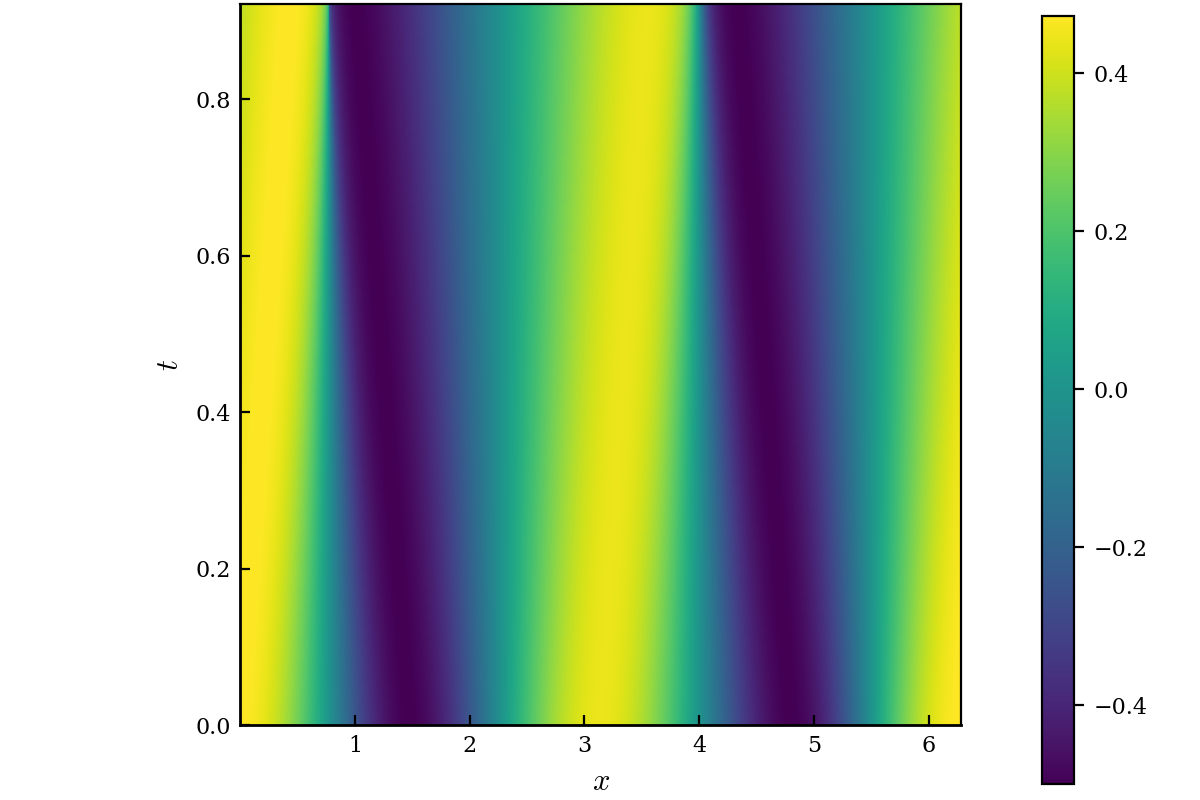} }}
        \qquad
        \subfloat[$r=108$ custom basis function solution]{{\includegraphics[width=0.45\textwidth]{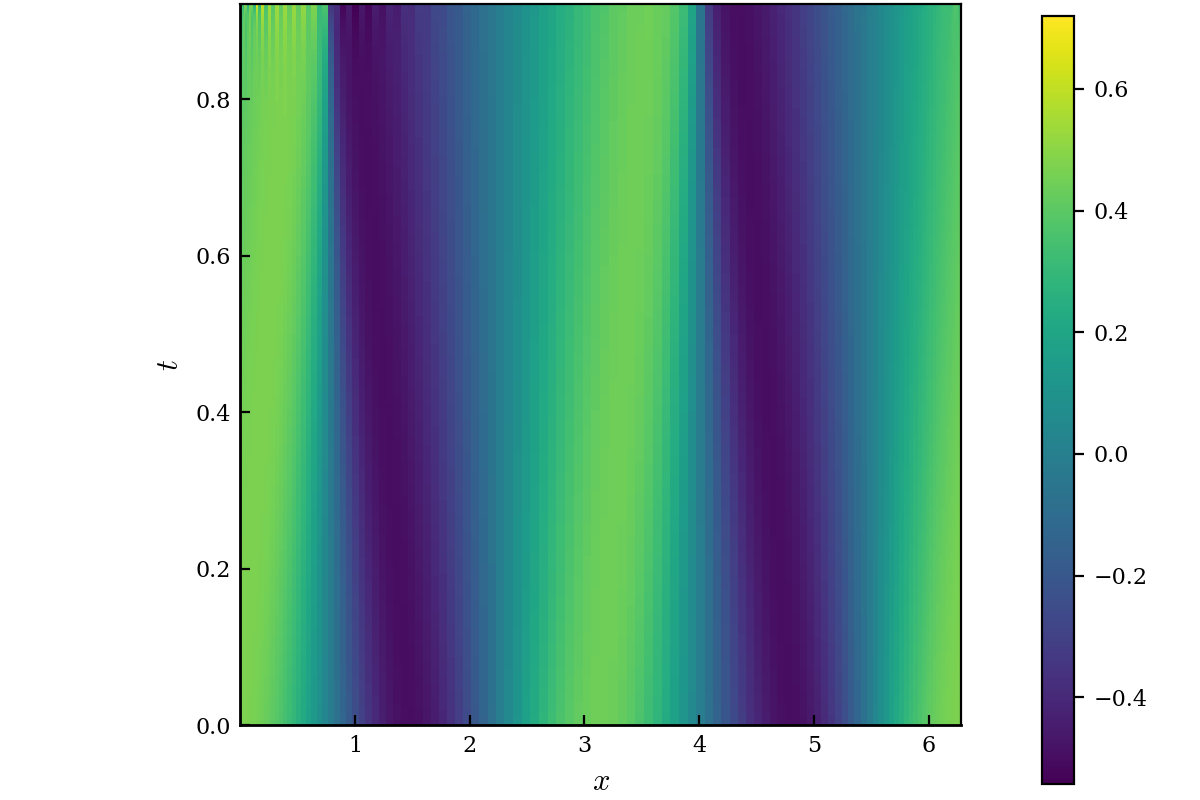} }}
        \qquad
        \caption{Results for the inviscid Burgers equation for $t \in [0,0.9221]$ and the random in training distribution initial condition.}
        \label{fig:inviscid_burgers_comparo_in}
    \end{figure}

\begin{figure}[h!]
        \centering
        \subfloat[The MUSCL solution]{{\includegraphics[width=0.45\textwidth]{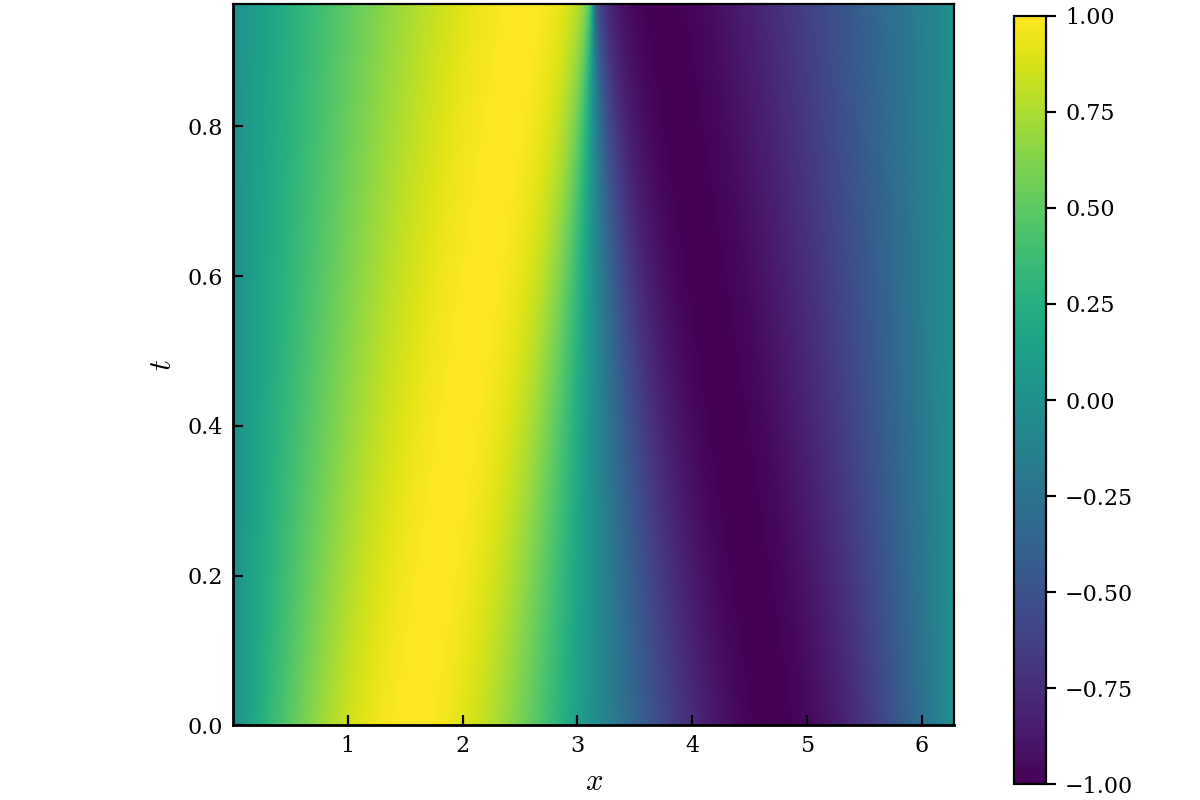} }}
        \qquad
        \subfloat[$r=108$ custom basis function solution]{{\includegraphics[width=0.45\textwidth]{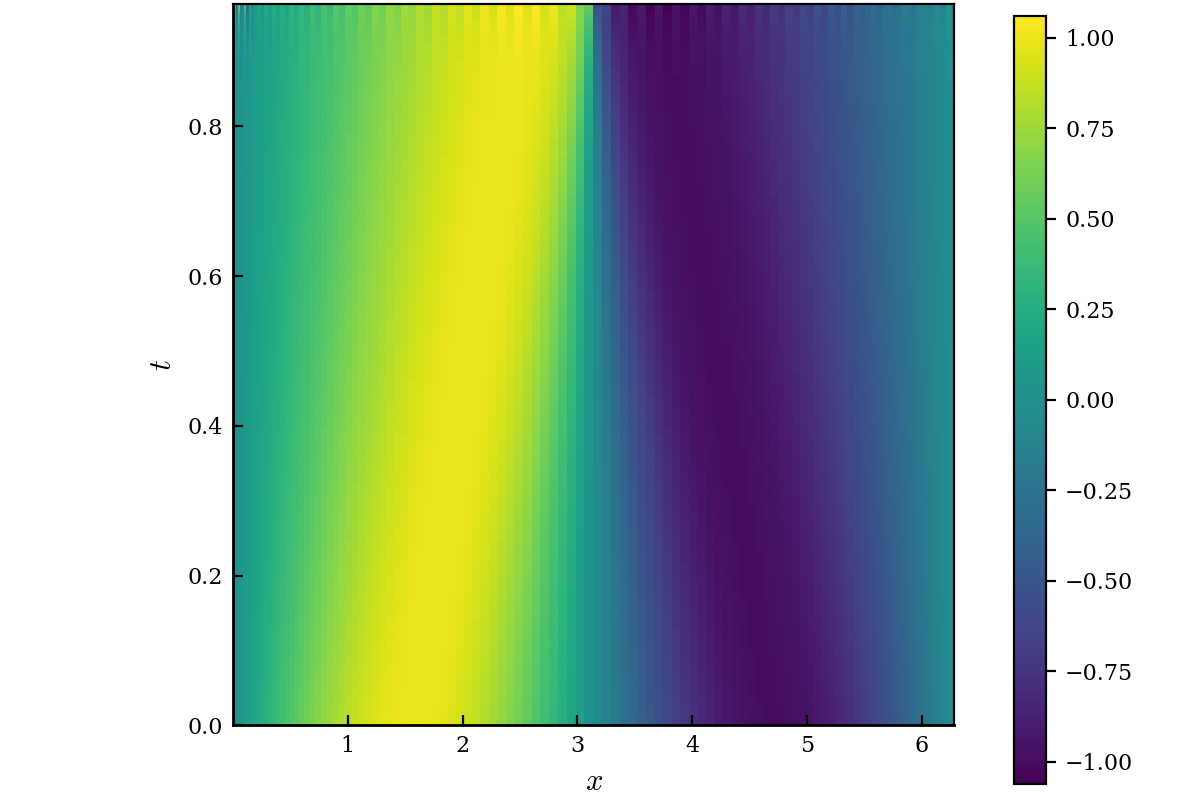} }}
        \qquad
        \caption{Results for the inviscid Burgers equation for $t \in [0,0.9636]$ and the out of training distribution initial condition $u(0,x) = \sin(x)$}
        \label{fig:inviscid_burgers_comparo_out}
    \end{figure}

Results using a MUSCL scheme with 4096 spatial discretization points and the custom basis function expansion for the evolution of the inviscid Burgers equation for $t \in [0,0.9221]$ and the in-distribution test initial condition are shown in Figure \ref{fig:inviscid_burgers_comparo_in}. See Figure \ref{fig:init_in_dist_inviscid_burgers} for a plot of the in-distribution test initial condition. Similarly, results for $t \in [0,0.9636]$ and the out of distribution test initial condition are shown in Figure \ref{fig:inviscid_burgers_comparo_out}. The final solution time presented for both test initial conditions corresponds to the time at which the solution becomes unstable and is specified as the time at which the energy in the system at a given solution time step exceeds 102.5\% of initial energy present in the system.

\begin{figure}[h!]
        \centering
        \subfloat[In-distribution random initial condition]{{\includegraphics[width=0.45\textwidth]{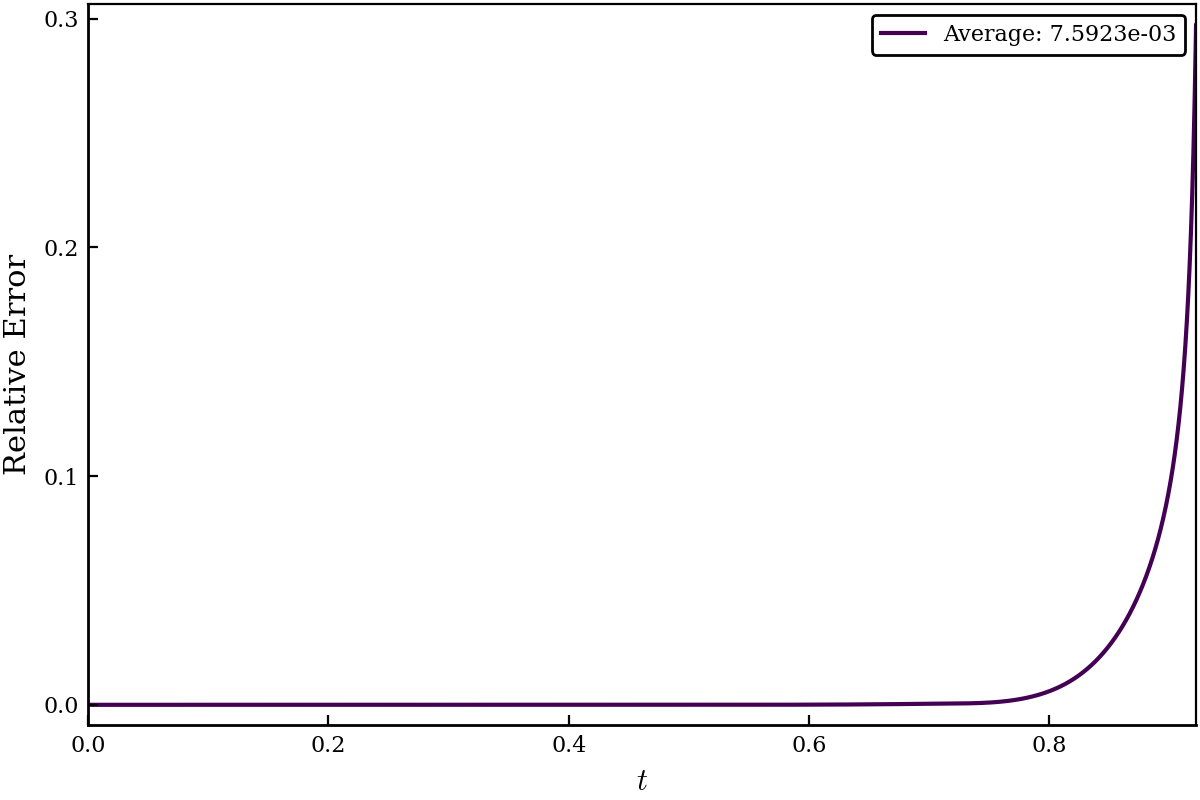} }}
        \qquad
        \subfloat[Out of distribution initial condition: $u(0,x) = \sin(x)$]{{\includegraphics[width=0.45\textwidth]{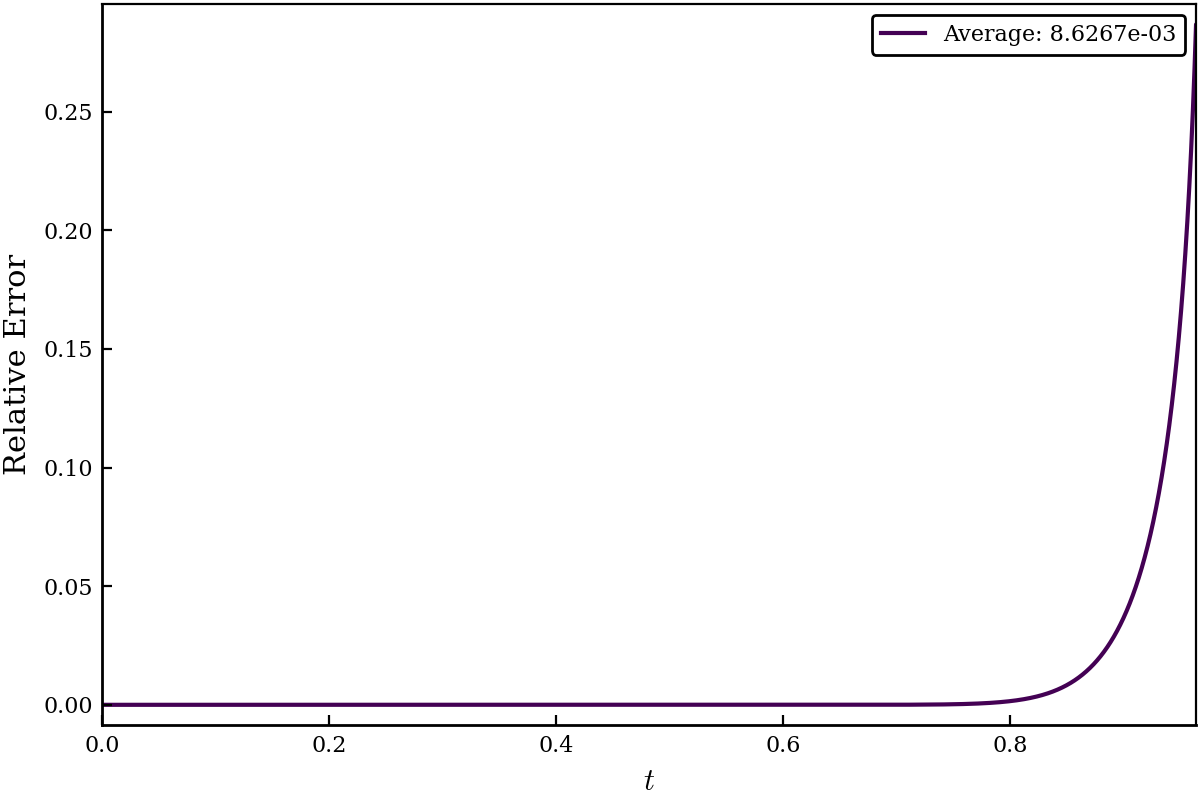} }}
        \qquad
        \caption{Relative errors for the invsicid Burgers equation.}
        \label{fig:inviscid_burgers_errors}
    \end{figure}

The relative and average errors compared to the MUSCL solution for the two test cases are presented in Figure \ref{fig:inviscid_burgers_errors}. For all results the triple product integral was computed directly for the nonlinear term.

\section{Discussion and future work}\label{sec:discussion_future_work}

We have presented a general framework for using DeepONets to identify spatial functions which can be transformed into a hierarchical orthonormal basis and subsequently used to solve PDEs. We illustrated this framework and its interpolation and extrapolation capabilities using four one-dimensional PDEs, two of which that are linear and two that are nonlinear. 

The results for the advection, advection-diffusion, and viscous Burgers equations show strong agreement with the $M=128$ Fourier mode solution over the entire temporal domain. In particular, that errors remain low for time values well beyond the temporal training interval of the DeepONet demonstrates the temporal extrapolation capabilities of the presented framework. This also illustrates the effectiveness of scientific machine learning techniques \cite{baker2019workshop} as the presented framework consists of embedding the information gleaned from a neural network, which is purely data driven, into the PDE and solving it using conventional techniques.

Results were also presented for the inviscid Burgers equation, which, unlike the other three examples whose solutions remain smooth over time when initialized from a smooth initial condition, can develop shocks in finite time. As a result, a temporal interpolation interval shorter than the training domain of $t \in [0,3.5]$ was used. For the time values before the shock, we obtain strong agreement between the custom basis function solution and the ground truth MUSCL solution. However, as evidenced by Figure \ref{fig:inviscid_burgers_errors}, around the time instant when the shock forms, the approximate solution becomes more inaccurate and is eventually no longer stable. This increasing level of inaccuracy should not be seen as a shortcoming of the presented framework; instead, this is an issue commonly encountered when using spectral methods for the evolution of singular PDEs \cite{hesthaven2007spectral}. This fact motivated the use of a MUSCL solution to generate the ground truth for training the DeepONets as the use of a Fourier expansion also provides inaccurate, albeit, stable results. The instability occurs due to the unavailability of a mechanism to eject the energy that is being consumed by the shock. To account for the ejection of energy and to accurately capture the evolution of the energy in time, we need to augment the system with a memory term (e.g. \cite{PriceJacob2021Orom}). In the case of the inviscid Burgers equation, the inclusion of a memory term allows for energy to be drained from the scales resolved by the simulation \cite{Stinis_2013}. Combining the presented framework with the methods developed in \cite{PriceJacob2021Orom} is an active area of investigation and will appear in a future publication.

For all four test PDEs, results were shown for two different initial conditions, one that was randomly selected from within the training distribution, and one that was outside the training distribution, $\sin(x)$. Referencing Figures \ref{fig:advection_errors}, \ref{fig:advection_diffusion_errors}, \ref{fig:viscous_burgers_errors}, and \ref{fig:inviscid_burgers_errors}, strong agreement is shown with the $M=128$ mode Fourier or MUSCL solution for both of initial conditions (in advance of the shock in the case of inviscid Burgers). With the exception of the advection equation results, we find an increase in the average error over the temporal interval for the out of distribution initial condition compared to the in-distribution initial condition; however, the presented results demonstrate the opportunity to extrapolate not only temporally, but also in terms of the input function space when utilizing the presented framework.

The presented general framework provides many interesting future research directions in addition to those already noted in this section. First, we need to perform meticulous optimization of DeepONet parameters to improve the quality of the custom basis functions. Second on the list is the development of a fast custom basis function inverse transform. For the results presented, the triple product integral was computed directly for all nonlinear terms which is computationally expensive and scales poorly, and hence is impractical for large scale computation \cite{canuto2012spectral}. Preliminary work is underway to develop a fast inverse custom basis function transform using DeepONets. These networks takes as the inputs to the branch and trunk nets the expansion coefficients and spatial locations respectively. Once trained, they will approximate the functions corresponding to the expansion coefficients. This will enable the use of a fast pseudospectral transform technique that provides a function on which auto differentiation may be performed so that nonlinear terms can be computed efficiently in real space. Third, as we move to problems on complex domains in higher dimensions, obvious generalizations of Legendre expansions are not available. As discussed in Section \ref{sec:basis_functions_operator_nets}, the development of alternative approaches is an active area of investigation. Fourth, a more thorough investigation of the time-sampling approach is also required, including an assessment of its performance on PDEs more complicated than the advection problem (see Appendix \ref{app:time-sampled}). Fifth, it is interesting to investigate if the custom-made basis functions developed for one PDE can be used to expand accurately the solution of another PDE (see Appendix \ref{app:use_another} for preliminary results).

Another interesting avenue for exploration is analyzing the basis functions obtained from DeepONets trained on time-independent problems. Our machinery can be deployed on solution operators for static equations that map, e.g., boundary data, forcing terms, or diffusivity coefficients to the solutions to yield promising custom bases. The elimination of the temporal dimension implies that, along with a possible reduction in the network training cost, the ambiguity associated with freezing the trunk net functions would be eradicated. 

Finally, the framework presented was initially based off the DeepONet architecture \cite{Lu_2021}, which is why we explicitly reference the trunk functions; however, there is reason to believe that this framework could be readily extended to other operator neural network architectures, e.g., \cite{li2020fourier,li2020neural}, and is the subject of further investigation.

\section{Acknowledgements}
We would like to thank Dr. Nathaniel Trask, Prof. George Karniadakis, and Prof. Lu Lu for helpful discussions and comments. The work of PS is supported by the DOE-ASCR-funded "Collaboratory on Mathematics and Physics-Informed Learning Machines for Multiscale and Multiphysics Problems (PhILMs)". Pacific Northwest National Laboratory is operated by Battelle Memorial Institute for DOE under Contract DE-AC05-76RL01830.

\bibliographystyle{plain}
\bibliography{bibliography}

\newpage

\appendix 

\renewcommand\thefigure{\thesection.\arabic{figure}}    
\setcounter{figure}{0}  
\renewcommand\thetable{\thesection.\arabic{table}}    
\setcounter{table}{0} 

\section{Solving partial differential equations using custom basis functions}\label{app:custom_basis_evolution}

We describe the general framework utilized for solving time-dependent partial differential equations (PDEs) using the custom basis functions presented in Section \ref{sec:basis_functions_operator_nets}. We begin by writing the approximate solution as a linear combination of the custom basis functions $\{\phi_k\}$
\begin{equation} \label{eqn:approx}
        u^r(t,x) = \sum_{k=1}^r a_k(t) \phi_k(x),
\end{equation}
where $r$ is the number of basis functions utilized. The test functions are specified to also be the custom basis functions; recall that the orthonormality condition
\begin{equation} \label{eqn:orthonormality}
    \langle \phi_k, \phi_m \rangle = \int_{\Omega} \overline{\phi_k(x)} \phi_m(x) dx = \delta_{km},
\end{equation}
holds for these functions. It should be noted, that based on the DeepONet architecture considered, the $\{\tau_k\}$ and resulting $\{\phi_k\}$ are always real-valued; however, the complex conjugation is shown for completeness.

Consider the general time-dependent PDE
\begin{equation}\label{eqn:general_pde}
    \frac{\partial u}{\partial t} = \tb{N} \left(u,\frac{\partial u}{\partial x},\dots,x,t \right),
\end{equation}
where \tb{N} is a, in general, nonlinear operator. To obtain a system of ordinary differential equations (ODEs) in terms of the expansion coefficients, we begin by inserting the approximate solution \eqref{eqn:approx} into \eqref{eqn:general_pde},
\begin{equation}
    \sum_{k=1}^r  \phi_k(x) \frac{da_k(t)}{dt} = \tb{N} \left(\sum_{k=1}^r a_k(t) \phi_k(x), \sum_{k=1}^r a_k(t) \phi'_k(x) ,\dots,x,t  \right),
\end{equation}
where the primes denote differentiation in space. Applying the Galerkin condition and using \eqref{eqn:orthonormality} yields the desired system of ODEs
\begin{equation}\label{eqn:system_odes}
    \frac{da_m(t)}{dt} = \ip{\phi_m , \text{\textbf{N}} \left(\sum_{k=1}^r a_k(t) \phi_k, \sum_{k=1}^r a_k(t) \phi'_k ,\dots,\cdot,t  \right) }.
\end{equation}
The initial condition for this system of equations is inferred from the initial condition $u(0,x)$ by
\begin{equation} \label{eqn:init_coeffs}
    a_m(0) = \ip{\phi_m,u(0,\cdot)}, \quad m = 1,\hdots,r.
\end{equation}

Note that the trial functions may not necessarily satisfy the boundary conditions of the PDE. To ensure compliance with the boundary conditions, the system of ODEs for the expansion coefficients \eqref{eqn:system_odes} are evolved for $m = 1,2,\dots,r-b$ where $b$ is the number of boundary conditions required for the PDE of interest. The $b$ number of boundary conditions are approximated using \eqref{eqn:approx} and enforced at each time step to determine the remaining expansion coefficients in an implementation similar to the tau-method for Chebyshev or Legendre polynomials (see e.g. \cite{Boyd_2001,canuto2012spectral} for a detailed discussion).

\section{Performance of time-sampled basis functions}\label{app:time-sampled} 
\begin{figure}[h!]
        \centering
        \subfloat[In-distribution random initial condition]{{\includegraphics[width=0.45\textwidth]{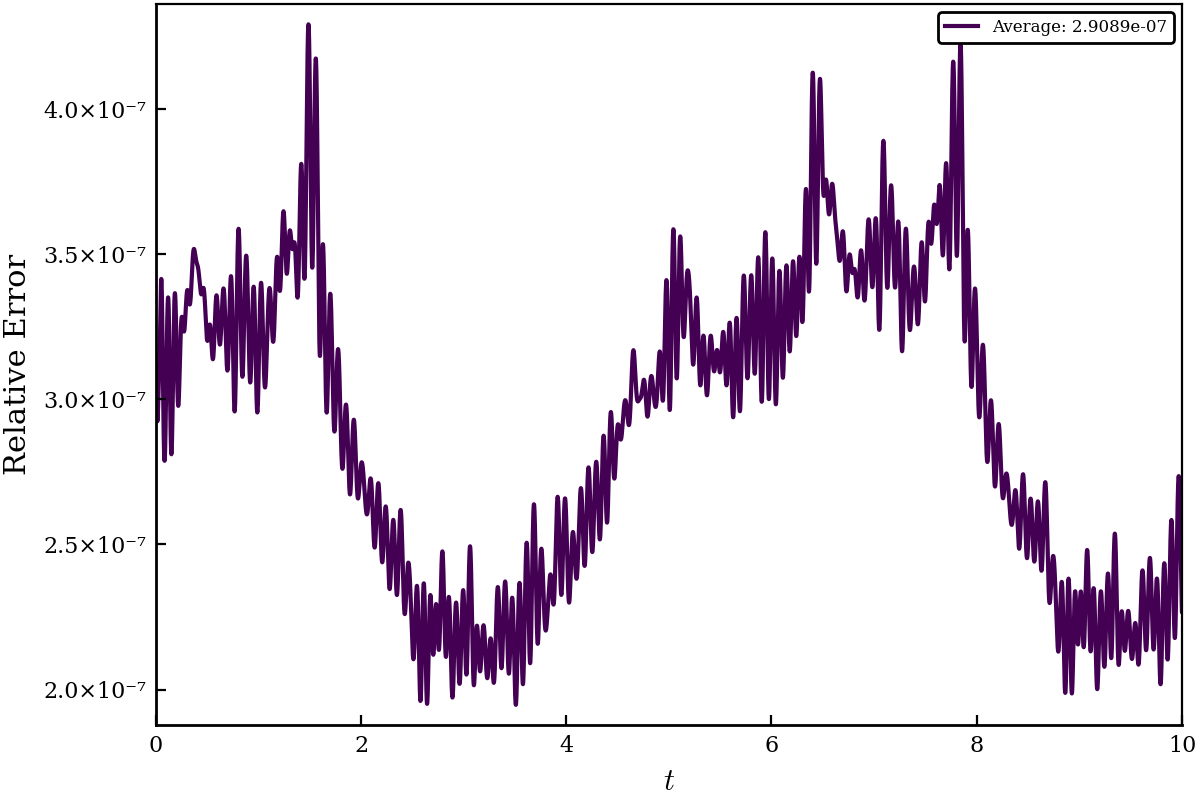} }}
        \qquad
        \subfloat[Out of distribution initial condition: $u(0,x) = \sin(x)$]{{\includegraphics[width=0.45\textwidth]{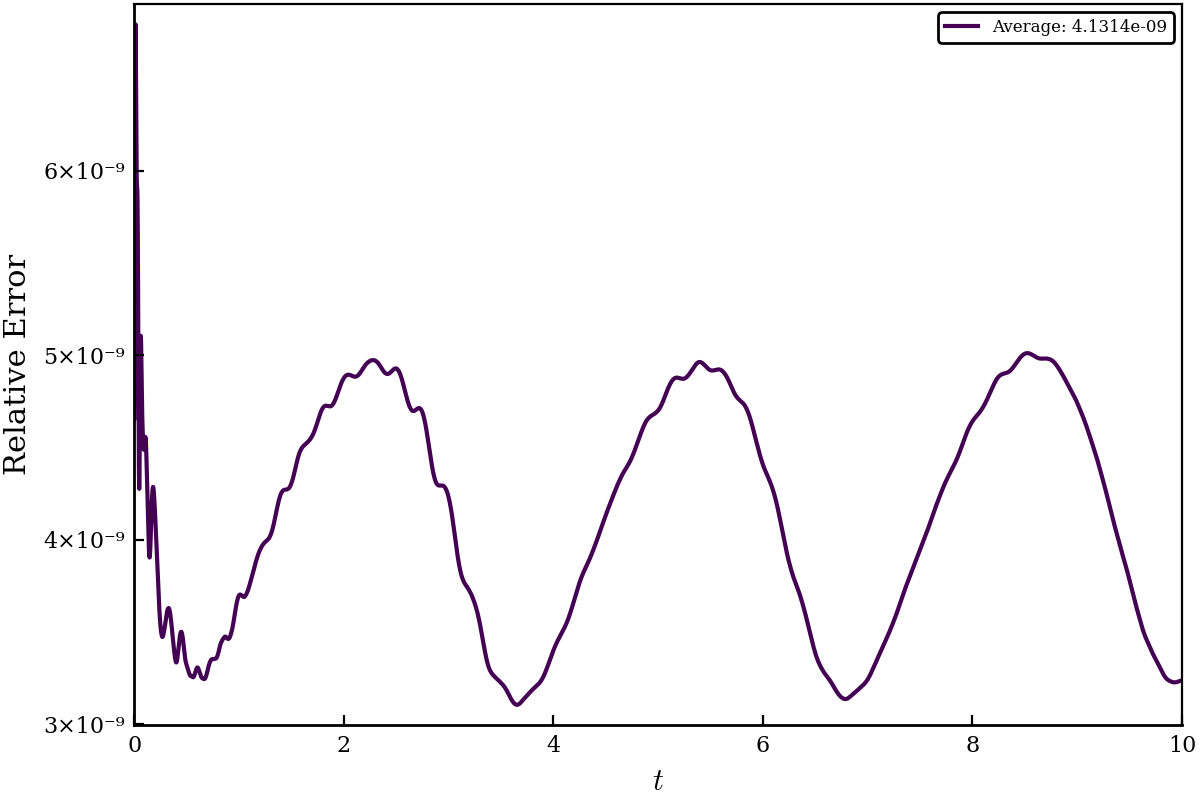} }}
        \qquad
        \caption{Relative errors for the advection equation for $t \in [0,10]$ and using a time-sampled set of candidate basis functions ($r = 56$).}
        \label{fig:advection_errors_sampling}
    \end{figure}

As was seen in Section \ref{sec:approx_capability}, the approximation capability of the custom basis is significantly improved by sampling the trunk net functions at various points in the temporal domain. Figure \ref{fig:advection_errors_sampling} shows representative results if this is applied to the advection equation using, as before, a sampling gap $\Delta t = 0.05$ for $t \in [0,1]$. Keeping all other evolution parameters constant, the algorithm finds 56 custom basis functions with a corresponding singular value greater than $10^{-7}$. The usage of time-sampling and the increased set of basis functions ($56$ vs. $40$) for solving the PDE results in a decrease in the average error over the testing interval for the in-distribution initial condition, while for the out of distribution problem, we find a slight increase in the average error over the temporal domain. 

\begin{figure}[h!]
        \centering
        \subfloat[In-distribution random initial condition]{{\includegraphics[width=0.45\textwidth]{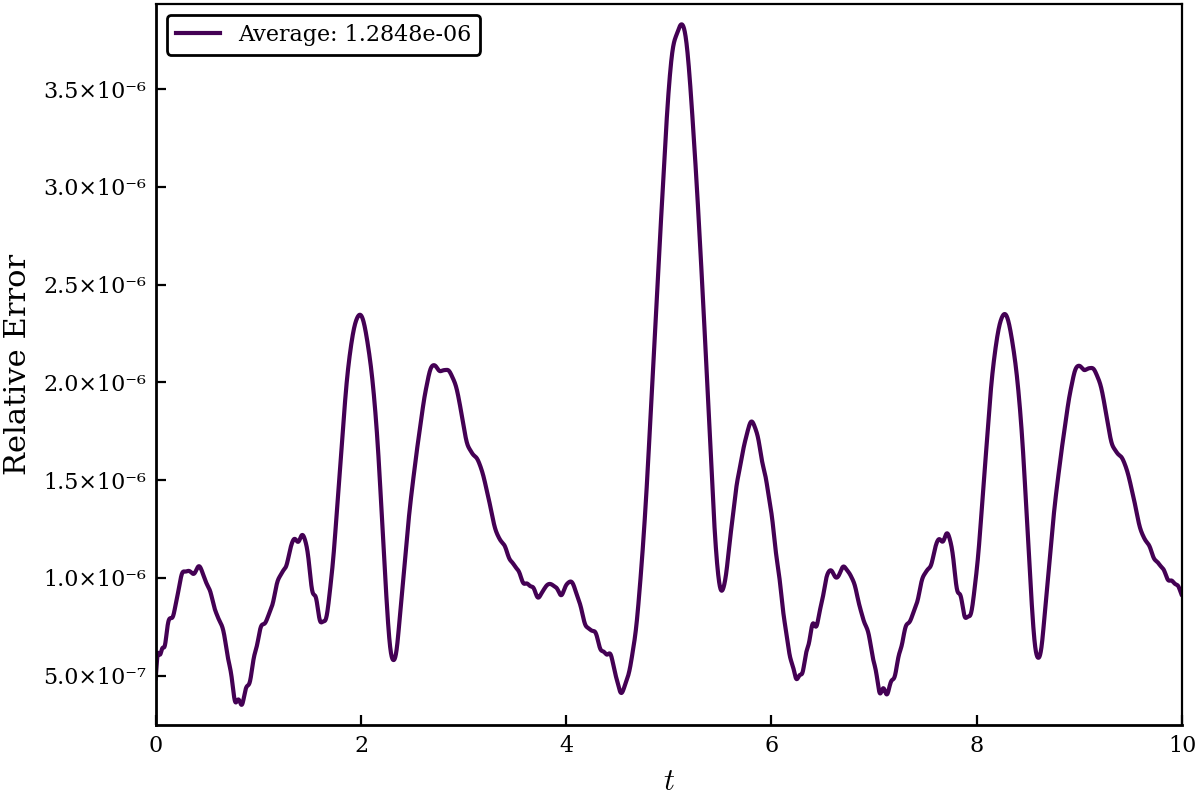} }}
        \qquad
        \subfloat[Out of distribution initial condition: $u(0,x) = \sin(x)$]{{\includegraphics[width=0.45\textwidth]{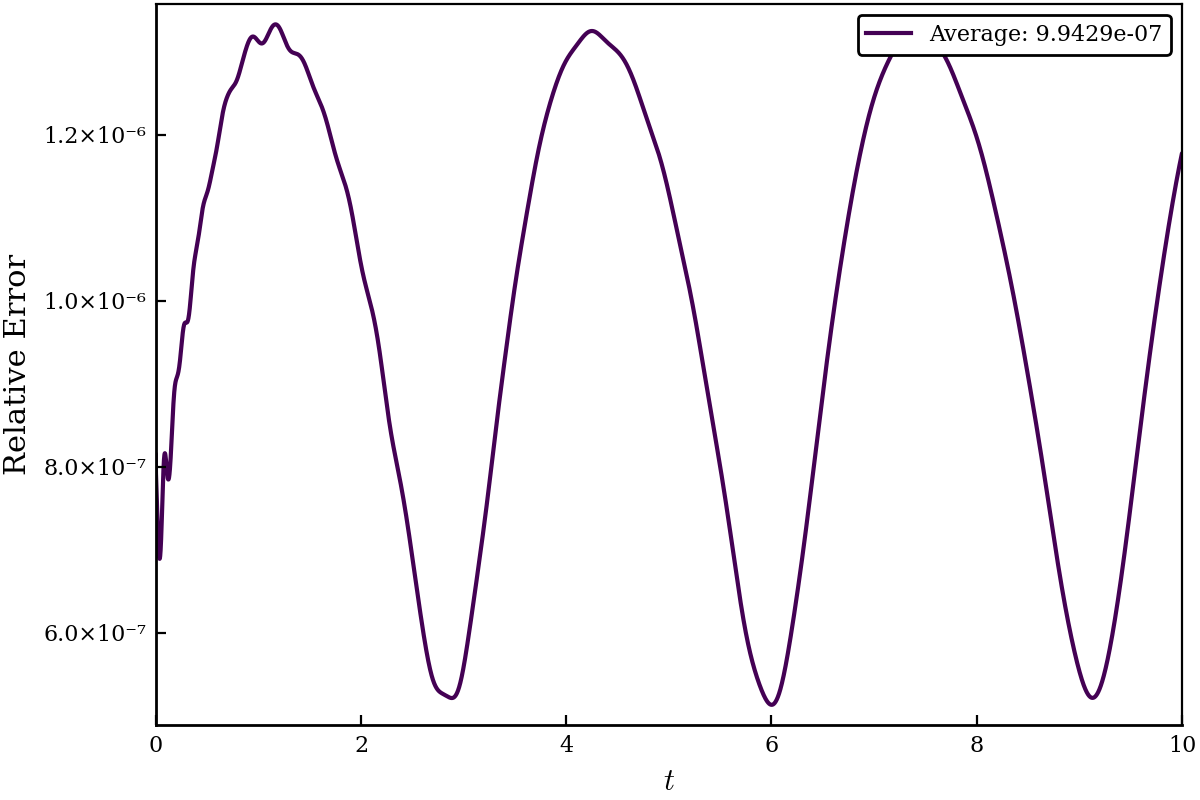} }}
        \qquad
        \caption{Relative errors for the advection equation for $t \in [0,10]$ and using a time-sampled set of candidate basis functions ($r = 40$). }
        \label{fig:advection_errors_sampling_40}
    \end{figure}

Figure \ref{fig:advection_errors_sampling_40} shows the comparison when using the same number of custom basis functions ($r = 40$) taken from $t = 0$ and the time-sampled set of functions. In this case, we find that the latter approach results in a slight decrease in the average error for the in-distribution initial condition, while for the out of distribution problem, we find an increase in the average error over the temporal domain. 

\begin{figure}[h!]
        \centering
        \subfloat[In-distribution random initial condition]{{\includegraphics[width=0.45\textwidth]{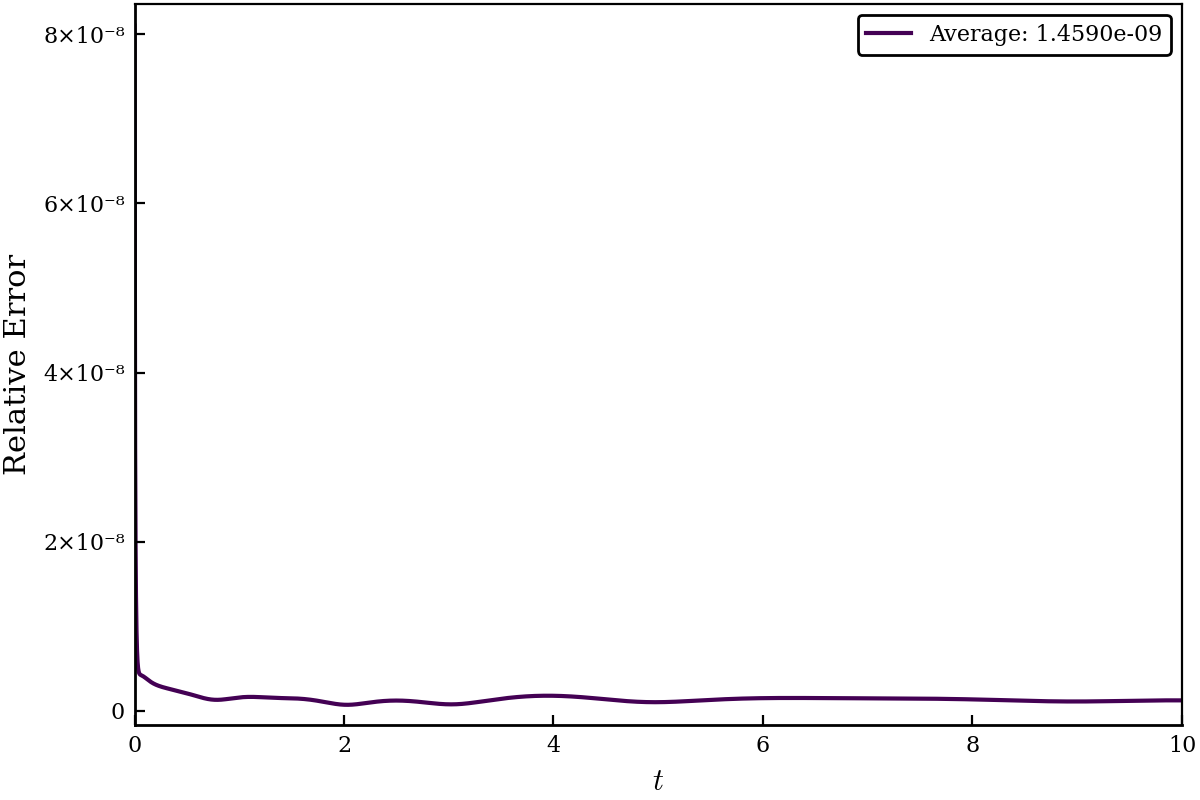} }}
        \qquad
        \subfloat[Out of distribution initial condition: $u(0,x) = \sin(x)$]{{\includegraphics[width=0.45\textwidth]{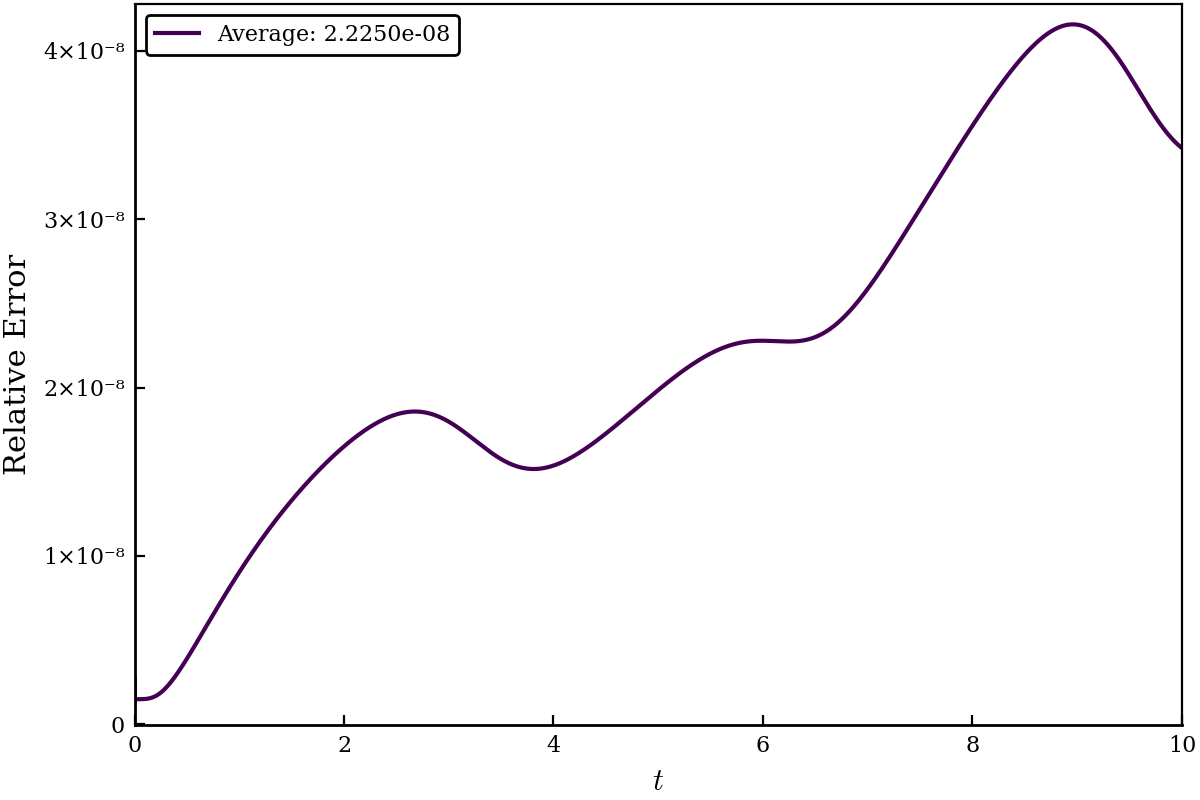} }}
        \qquad
        \caption{Relative errors for the advection-diffusion equation for $t \in [0,10]$ and using the time-sampled set of custom basis functions ($r = 50$).}
        \label{fig:advection_diffusion_errors_sampling}
    \end{figure}

\begin{figure}[h!]
        \centering
        \subfloat[In-distribution random initial condition]{{\includegraphics[width=0.45\textwidth]{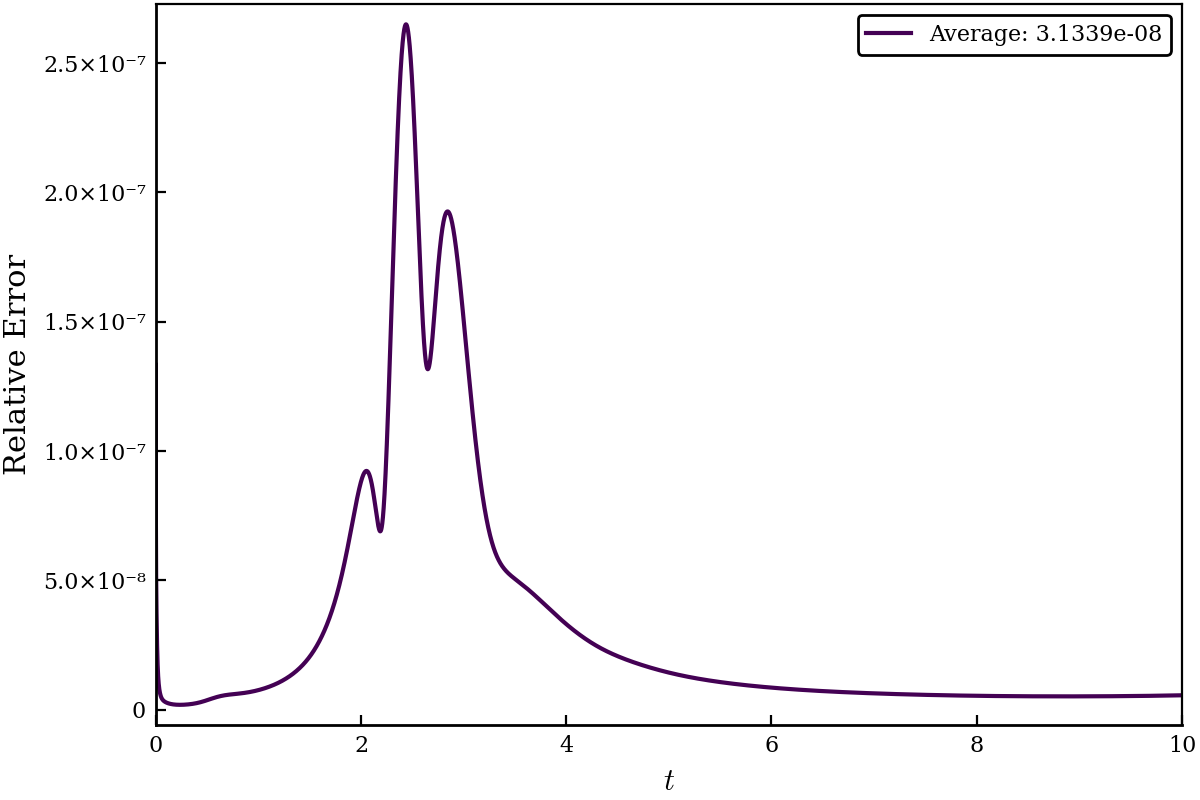} }}
        \qquad
        \subfloat[Out of distribution initial condition: $u(0,x) = \sin(x)$]{{\includegraphics[width=0.45\textwidth]{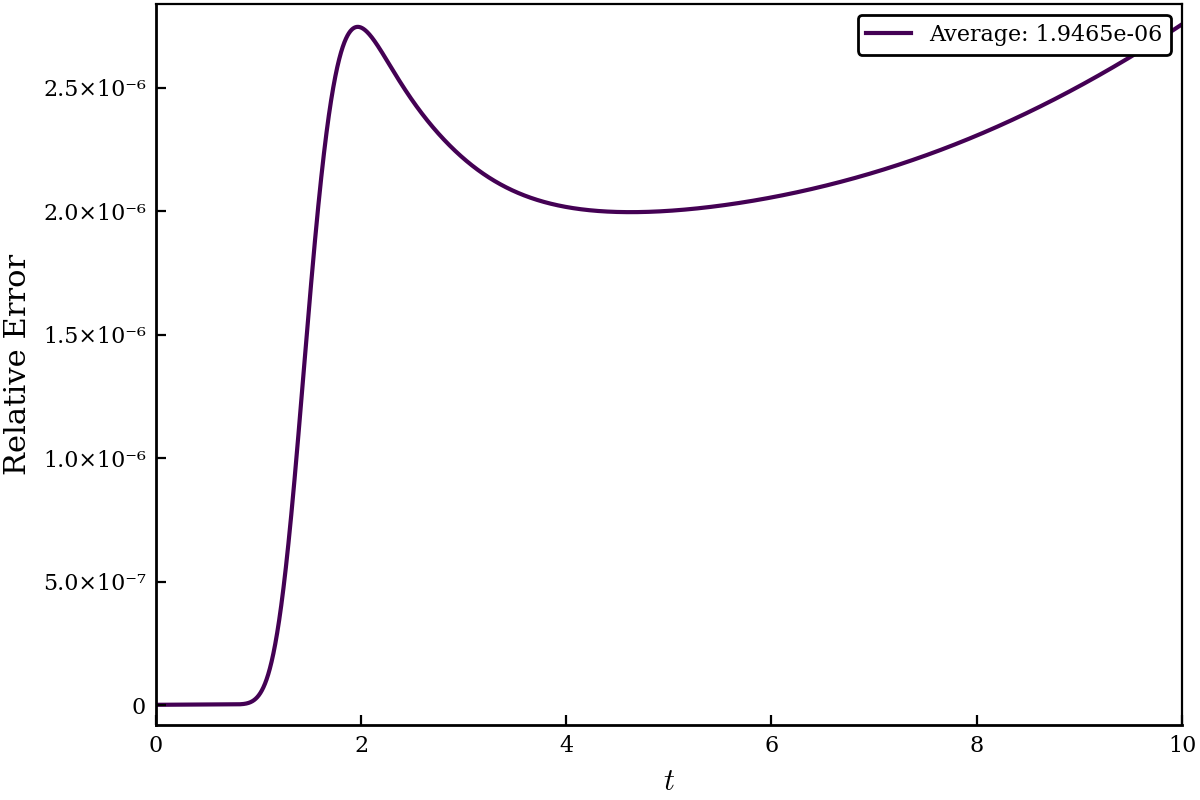} }}
        \qquad
        \caption{Relative errors for the viscous Burgers' equation for $t \in [0,10]$ and using the time-sampled set of custom basis functions ($r = 98$).}
        \label{fig:viscous_burgers_errors_sampling}
    \end{figure}

\begin{figure}[h!]
        \centering
        \subfloat[In-distribution random initial condition, {$t \in [0,0.9231]$}]{{\includegraphics[width=0.45\textwidth]{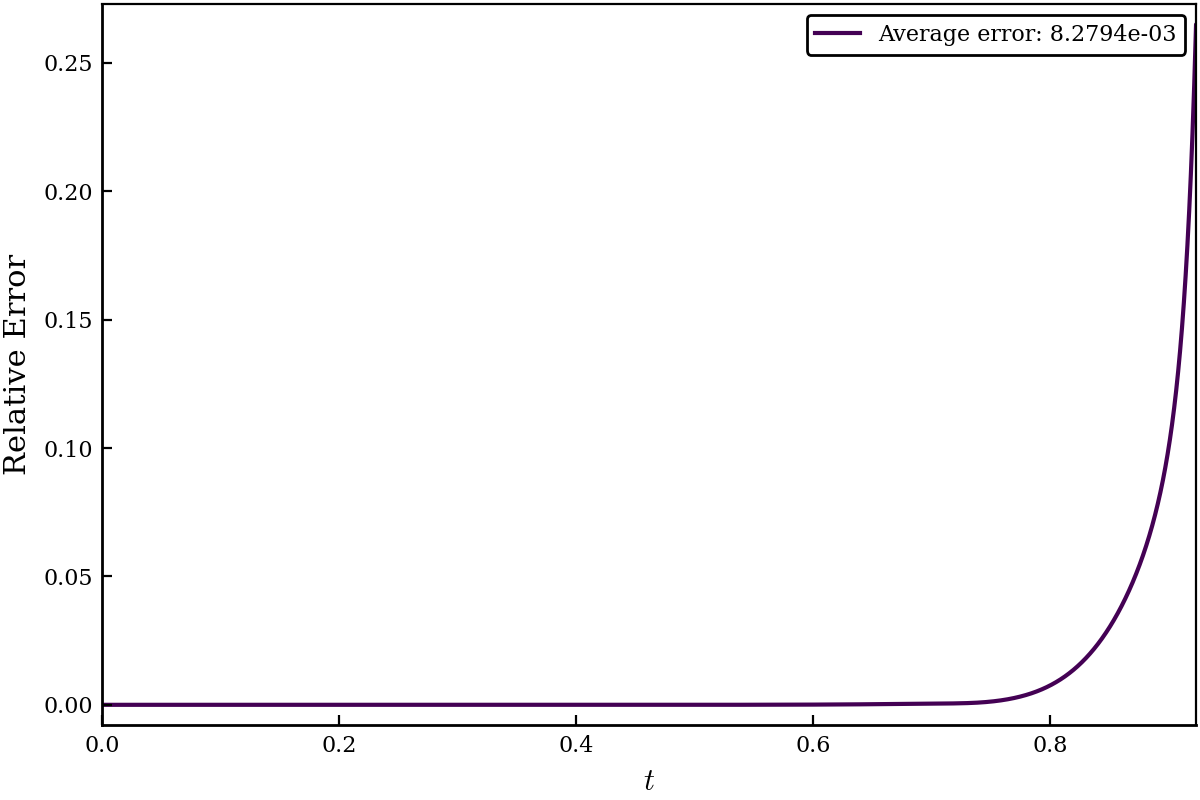} }}
        \qquad
        \subfloat[Out of distribution initial condition: $u(0,x) = \sin(x)$, {$t \in [0,0.9463]$}]{{\includegraphics[width=0.45\textwidth]{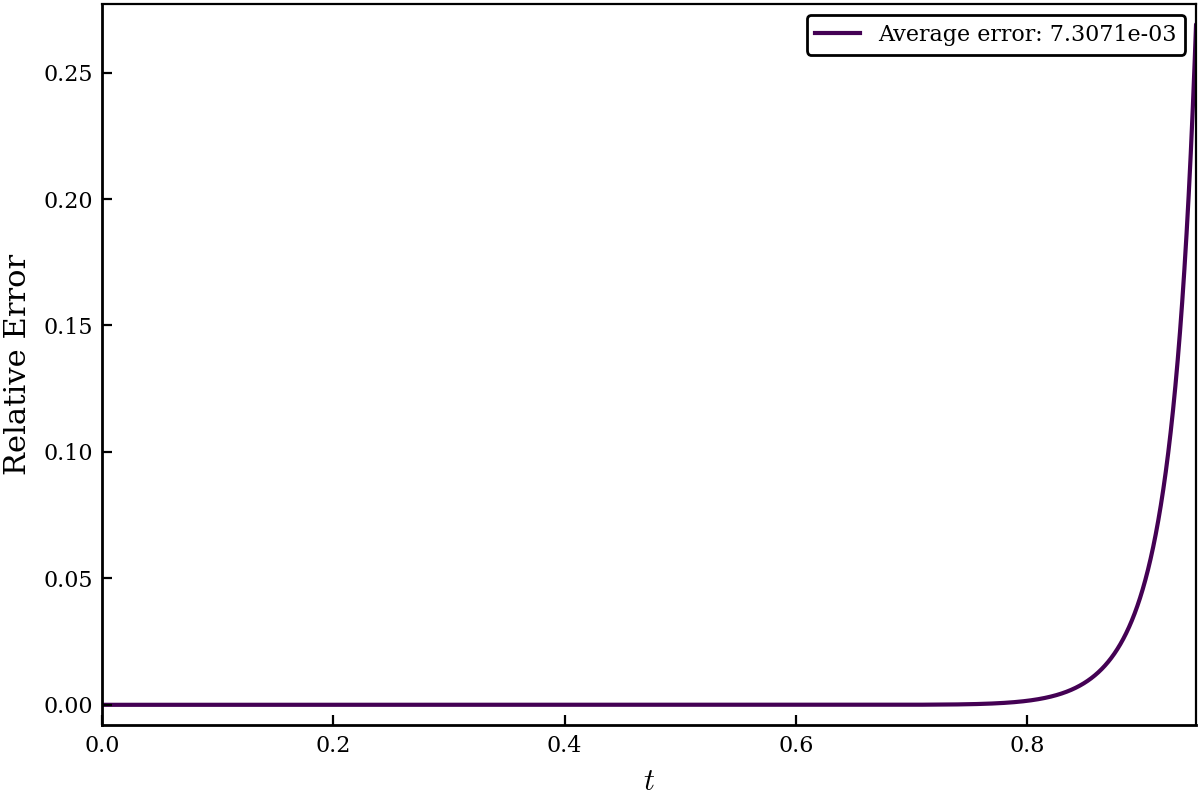} }}
        \qquad
        \caption{Relative errors for the inviscid Burgers' equation using the time-sampled set of custom basis functions ($r = 106$).}
        \label{fig:inviscid_burgers_errors_sampling}
    \end{figure}

Representative results for the advection-diffusion equation, viscous Burgers equation, and inviscid Burgers equation are shown on Figure \ref{fig:advection_diffusion_errors_sampling}, Figure \ref{fig:viscous_burgers_errors_sampling}, and Figure \ref{fig:inviscid_burgers_errors_sampling} respectively. For all of these results a sampling gap $\Delta t = 0.05$ for $t \in [0,1]$ was used in conjunction with a singular value threshold of $10^{-7}$. For the advection-diffusion and viscous Burgers equation, the use of time-sampled basis functions leads to an increase in accuracy relative to using basis functions taken from $t=0.$ For the case of inviscid Burgers there is no change in accuracy relative to using basis functions taken from $t=0.$ 

Sampling the trunk function at various times throughout the temporal training domain is the subject of continued investigation and will appear in a future publication. 

\section{Use the custom-made basis functions from one PDE to expand the solution of another PDE}\label{app:use_another}

\begin{figure}[h!]
        \centering
        \subfloat[In-distribution random initial condition]{{\includegraphics[width=0.45\textwidth]{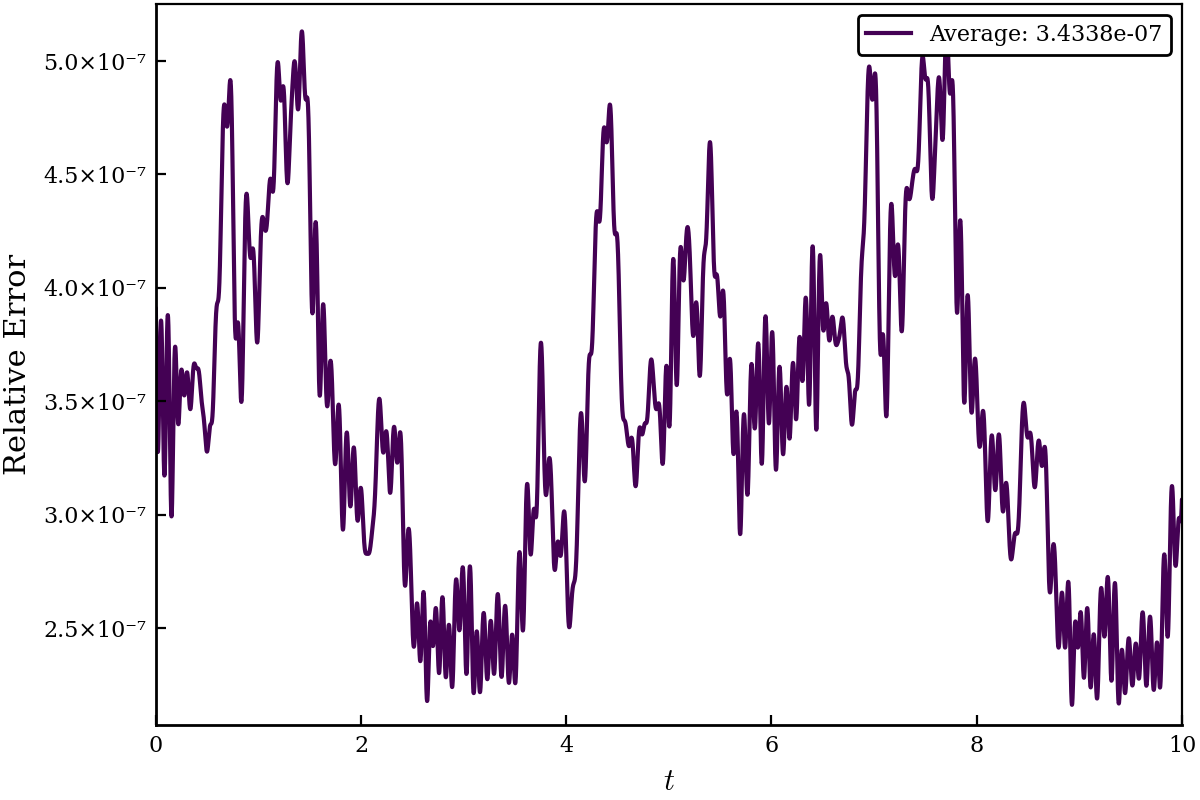} }}
        \qquad
        \subfloat[Out of distribution initial condition: $u(0,x) = \sin(x)$]{{\includegraphics[width=0.45\textwidth]{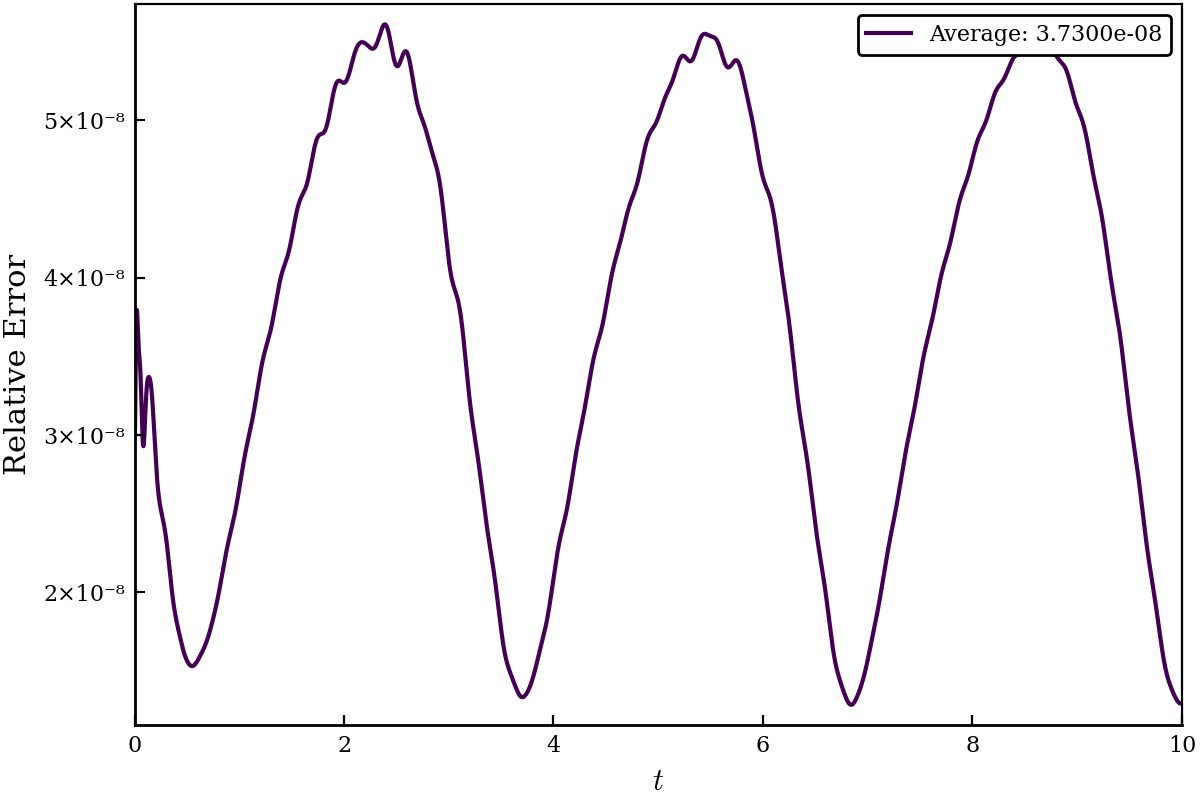} }}
        \qquad
        \caption{Relative errors for the advection equation for $t \in [0,10]$ and using the custom basis functions found for the viscous Burgers equation ($r = 58$).}
        \label{fig:advection_errors_viscous_burgers}
    \end{figure}

For each example PDE, a DeepONet was trained for that specific PDE; however, the natural question arises: can the custom basis functions identified for PDEs which exhibit more complex dynamics be utilized for simpler systems? Figure \ref{fig:advection_errors_viscous_burgers} shows results for the advection equation and for the two test initial conditions when using the custom basis functions identified for the viscous Burgers equation ($\nu=0.1$). From the test case shown in Figure \ref{fig:advection_errors_viscous_burgers}, we find a decrease in the average error for the in-distribution initial condition and an increase in the average error for the out of distribution initial condition. The concept of using basis functions trained for more complex PDEs for a general class of PDEs is also the subject of continued investigation and will appear in a future publication.

\end{document}